\documentclass[11pt]{report}
\usepackage{amssymb,amsmath,enumerate,latexsym,geometry}
\usepackage[all]{xy}
\begin{document}

\newtheorem{example}{Example}[section] 
\newtheorem{Def}[example]{Definition}  
\newtheorem{fact}[example]{Fact} 
\newtheorem{rem}[example]{Remark}  
\newtheorem{con}[example]{Conjecture} 
\newtheorem{prop}[example]{Proposition}
\newtheorem{cor}[example]{Corollary} 
\newtheorem{thm}[example]{Theorem} 
\newtheorem{lem}[example]{Lemma} 
\newenvironment{out}{\noindent {\bf Outline proof} }{\mbox{} \hfill $\Box$ \mbox{}\\} 
\newenvironment{proof}{\noindent {\bf Proof} }{\mbox{} \hfill $\Box$ \mbox{}\\} 

\newcommand{\id}{i\!d}
\newcommand{\ep}{\varepsilon}
\newcommand{\bA}{\mathsf{A}}
\newcommand{\bB}{\mathsf{B}}
\newcommand{\bC}{\mathsf{C}}
\newcommand{\bF}{\mathsf{F}}
\newcommand{\bI}{\mathsf{1}}
\newcommand{\bN}{\mathbb{N}}
\newcommand{\bP}{\mathsf{P}}
\newcommand{\bZ}{\mathbb{Z}}
\newcommand{\dd}{\delta}
\newcommand{\Xt}{\widetilde{X}}
\newcommand{\Gt}{\widetilde{G}}
\newcommand{\Rt}{\widetilde{R}}
\newcommand{\db}{\tilde{\delta}}
\newcommand{\Cb}{\widetilde{C}}
\newcommand{\lra}{\leftrightarrow}
\newcommand{\lan}{\langle}
\newcommand{\ran}{\rangle}
\newcommand{\arr}{\mathrm{Arr}}
\newcommand{\ob}{\mathrm{Ob}}
\newcommand{\sets}{\mathsf{Sets}}
\newcommand{\relb}{\mathrm{RelB}}
\newcommand{\cF}{\mathcal{F}}
\newcommand{\leadmon}{\mathtt{LM}}
\newcommand{\leadco}{\mathtt{LC}}
\newcommand{\reman}{\mathtt{rem}}
\newcommand{\spol}{\mathtt{spol}}
\newcommand{\tags}{\mathtt{tag}}
\pagestyle{empty}

\begin{titlepage}
\begin{center}
\textbf{\huge Applications of Rewriting Systems and \\
\vspace{0.2cm}
               Gr\"obner Bases to Computing Kan Extensions \\
\vspace{0.2cm}
               and Identities Among Relations}\\
{\Large
\vspace{1cm}
Thesis submitted to the University of Wales in support of \\
the application for the degree of Philosophi\ae   Doctor\\
\vspace{1cm}
by\\
\vspace{0.7cm}
A. Heyworth\\
\vspace{0.7cm}
supervised by\\
\vspace{0.7cm}
Prof. R. Brown and Dr. C. D. Wensley\\
\vspace{0.7cm}
October 1998\\
\vspace{0.7cm}
}
\end{center}
\vspace{1cm}

{\large
\underline{KEY WORDS:} 
          Presentation, Congruence, Category, Kan Extension, Rewrite 
System,\linebreak 
Gr\"obner Basis, Normal Form, Automaton, Regular Expression, 
Peiffer Relation,\linebreak 
Module, Crossed Module,  Identities Among Relations,
Covering Groupoid, \linebreak Contracting Homotopy, \ Crossed Resolution.\\

\vspace{2cm}

A. Heyworth,\\
School of Mathematics,\\
University of Wales,\\
Bangor,\\
Gwynedd LL57 1UT.\\}                            

\end{titlepage}

\pagestyle{plain}
\setcounter{page}{1}
\pagenumbering{roman}

\tableofcontents
\begin{center}
\textbf{\Large Summary}
\end{center}

This thesis concentrates on the development and application of Gr\"obner bases
methods to a range of combinatorial problems (involving groups, semigroups,
categories, category actions, algebras and $K$-categories) and the use of 
rewriting for calculating Kan extensions.

The first chapter gives a short introduction to 
presentations, rewrite systems, and completion.\\

Chapter Two contains the most important result, which is the application of 
Knuth-Bendix procedures to Kan extensions,
showing how rewriting provides a useful method for attempting to solve a variety 
of combinatorial problems which can be phrased in terms of Kan extensions.
A GAP3 program for Kan extensions is included in the appendix.\\


Chapter Three shows that the standard Knuth-Bendix algorithm is 
step-for-step a special case of Buchberger's algorithm. 
The one-sided cases and higher dimensions are considered, and the relations
between these are made precise.
The standard noncommutative Gr\"obner basis 
calculation may be expressed as a Kan extension over modules.
A noncommutative Gr\"obner bases program (in $\mathsf{GAP3}$) has been written.\\


Chapter Four relates rewrite systems, Gr\"obner bases and automata.
Automata which only accept irreducibles, and automata which output reduced
forms are discussed for presentations of Kan extensions.
Reduction machines for rewrite systems are identified with standard output 
automata and the reduction machines devised for algebras are expressed as
Petri nets.\\


Chapter Five uses the completion of a group rewriting system to algorithmically 
determine a contracting homotopy necessary in order to compute the set of
generators for the module of identities among relations using the covering 
groupoid methods devised by Brown and Razak Salleh \cite{BrSa97}. 
(The resulting algorithm has been implemented in GAP3). 
Reducing the resulting set of submodule generators is identified as a
Gr\"obner basis problem.\\

\begin{center}
\textbf{\Large Acknowledgements}
\end{center}

 I would first like to express my deepest appreciation to my parents. \\
 
 I would also like to thank the School of Mathematics at the University 
 of Wales in Bangor for giving me the opportunity to do the PhD and
 a friendly environment in which to work. I am particularly grateful to
 Ronnie Brown and Chris Wensley for their joint supervision and
 encouragement and to
 Tim Porter and Larry Lambe for their additional advice and inspirations.\\
 
 I am grateful to many true and kind friends who have encouraged me. 
 My wonderful brother Ben.
 My `little sisters' Angie and Nergiz;
 My good friend Tanveer; 
 The lovely people at Barnardo's -- especially Si\"an and Yvonne;
 My neighbours in Rachub -- especially Emma H and family;
 Helen, Emma M and Val who made University a nicer place to be.\\
 
 There is so much to learn \ldots

\setcounter{chapter}{0}
\setcounter{page}{1}
\pagestyle{plain}
\pagenumbering{arabic}

\chapter{Introduction}

\section{Presentations}

\subsection{Background}

A computational problem in group theory typically begins ``Given a group $G$, 
determine...''. Methods of solution of the problem depend on the way the
information about $G$ is given. The study of groups given by presentations is
called combinatorial group theory. 
Study of other algebraic objects (for example categories) through presentations
may be called combinatorics.
This section is an attempt to outline a little of the (controversial) history of
and motivation for the study of groups and in particular the use of group
presentations.

The origins of group theory might go back to 1600 BC.
Stone tablets remain as evidence that the Babylonians knew how to solve 
quadratic equations (though they had no algebraic notation). The solution (by
radicals) of a cubic equation was not discovered until the 16th century, and
published simultaneously with the method for solving quartics (by reducing to a 
cubic). Mathematicians such as Euler and Lagrange worked on the problem, and in
1824 Abel proved that there was no general solution by radicals of a quintic
equation. Work began on determining whether a given quintic was soluble, and
it is from \'Evariste Galois's paper ``On the Conditions of Solubility of 
Equations by Radicals'' (submitted and rejected in 1831) that group theory 
really began.
(That is not to say that group theoretic ideas did not exist before Galois 
(according to \cite{Sims}, they did) and a number of results were obtained 
before the definition of an abstract group reached its final form.)
The first formal development of group theory followed Galois's ideas and 
was limited almost entirely to finite groups. 
The idea of an abstract infinite group is included in Arthur Cayley's work 
(1854,
1878) on group axioms, but was not pursued at that time.
Finitely generated groups were defined by Dyck in 1882, and it is (disputedly) 
here that the first definition of a presentation by generators and relations 
was given.

Studying groups became important; groups of transformations came from 
symmetries and congruences in Euclidean Geometry, (semigroups come from partial
symmetries) automorphism groups were used in Klein's ``Erlangen Programme'', 
cyclic groups came from numbers and modular arithmetic and more groups from 
Gauss's composition of binary quadratic forms (groupoids from Brandt's 
generalisation of this problem). 
Abstract finite groups were defined by Weber in 1882, and it was in 1893 that he
published what we recognise as the modern definition of an arbitrary abstract 
group.

A major stimulus to the study of infinite discrete groups,
however, was the development of topology. In 1895 Poincar\'e introduced the notion of a 
fundamental group $\Pi_1(X,a)$ of closed paths of a space $X$ from a point $a$.
The properties of the fundamental group of a topological space correspond to 
some properties of the space. Interest in classifying the topological spaces 
generated interest in fundamental groups.
In 1911 Max Dehn, a student of Hilbert's, wrote a paper \cite{Dehn} which 
dealt with presentations of fundamental groups of closed, orientable surfaces, 
for which he formulated three fundamental decision problems: the word problem, 
the conjugacy problem, and the isomorphism problem. It is thought that by this 
time the idea of trying to determine properties of a group given by a finite 
presentation was already familiar. Anyway, some consider the problems 
to be part of what became known as ``Hilbert's Programme''.
Nielsen was also an important influence: his work led naturally to the study of 
groups presented through generators and relators.

There are certain advantages of presentations as a method for studying groups, 
or indeed other algebraic structures (monoids, categories, algebras). 
One advantage is that a presentation is compact as compared to (say) a Cayley 
table. An efficient presentation describes the group with the minimal amount of 
information.
By now there is a lot of theoretical machinery for working with presentations, 
this may be called computational group theory (or computational category theory,
etc), which really began with Turing and Newman's work at the end of World 
War II. 
Modern work in computational group theory may be found in Charles Sims's recent
book \cite{Sims}, and a lot of work developing computer programs for group 
theoretic computations continues at Warwick (KBMAG), St Andrews (GAP) and
Sydney (MAGMA) to name a few. The area has also broadened, problems with
monoids are more widely researched and now categories are
coming into the picture. Computational category theory is one relatively new
field of computer algebra which has considerable prospects.\\

Rewriting systems are sets of directed equations or rules which are useful in
computations. Rewrite rules specify the repeated replacement of subterms of a 
given formula with equivalent terms.
Rewriting theory was introduced as a method of solving the word problem.
The original word problem was expressed by Axel Thue in 1914:\\
\hspace{1cm}
``Suppose one has a set of objects, and a set of transformations (rules)
that when applied to these objects yield objects in the same set.
Given two objects $x$ and $y$ in the set, can $x$ be transformed into $y$, or
is there perhaps a third object $z$ such that both $x$ and $y$ can be 
transformed into $z$?''.\\

Thue established some preliminary results about strings of symbols (i.e. 
elements of a free monoid) and suggested that the approach might extend to more
structured combinatorial objects (at about this time Dehn was working on the 
beginnings of combinatorial group theory). Thue wanted to develop a ``calculus'' 
to decide the word problem, that is a set of procedures or algorithms that could 
be applied to the given objects to obtain the correct answer. He wanted 
a general algorithm to solve the word problem in a variety of different 
settings.\\

Apparently Thue's work was disregarded until the 1930's when logicians were 
seeking formal definitions of concepts like ``algorithm'' and ``effective 
procedure''.
In the mid 1950's and 60's notions of semi-Thue systems became important in 
mathematical linguistics. Work on formal language theory used semi-Thue systems
as mathematical models for phrase-structure grammars. 
At the same time technology was improving to the extent where mathematicians
began to consider mechanical theorem proving, and in the 1960's automated 
deduction quickly developed. As a form of computer program, rewriting systems 
made their debut in 1967 in a paper by Gorn.
A particularly influential role was played by a paper written by
Knuth and Bendix in 1970 \cite{KnBe}. They described an automatic procedure for solving word problems in abstract algebras.\\ 

In the 1970's term-rewriting systems took an important role in the study of 
automated deduction, which was still a rapidly developing area.
However, it was not really until the 1980's that Thue systems became popular. 
A book which contains the most fundamental results of the 1980's is \cite{BoOt}.
Since then, rewriting systems have continued to be of increasing interest, being
investigated for different properties and applied to a widening range of areas. 
The computational aspect is particularly important. Many modern programs for 
symbolic manipulation continue to use rewrite rules in an ad hoc manner, and 
there is now much work on the more formal use of rewriting systems in programming
(in particular see \cite{Holt1}\cite{Holt2}\cite{Sims}).

\subsection{Monoid and Group Presentations}

It is assumed that the reader is familiar with monoids and groups. The
terms and definitions for presentations are given in the following paragraphs to
fix the notation.\\

  Let $X$ be a set.
  The \textbf{free semigroup} $X^\dagger$ on $X$ consists of all nonempty sequences (strings) of 
elements of $X$. Composition is defined by concatenation of the strings.
  The \textbf{free monoid} $PX$ (sometimes denoted $X^*$) on $X$ consists of  all strings of elements of $X$, including the empty string. Composition is defined by string concatenation with the empty string acting as identity.\\

A \textbf{set of relations} $R$ for a monoid 
generated by $X$ is a subset of $PX \times PX$.
A \textbf{congruence} $=_S$ on a monoid $A$ is an equivalence relation on
$A$ such that, for all $u,v \in A$, if $l =_S r$ then $ulv =_S urv$.  
The \textbf{congruence $=_R$ generated by $R$ on $PX$}, where $R$ is a set of 
relations, is given by $x =_R y$ if and only if there is a system of equations 
\begin{eqnarray*} x & = & u_1 l_1 v_1 \\
        u_1 r_1 v_1 & = & u_2 l_2 v_2 \\
              \cdots&\cdots& \cdots \\
        u_n r_n v_n & = & y 
\end{eqnarray*} 
where either $(l_i, r_i)$ or $(r_i, l_i) \in R$ for $i=1,\ldots,n$, $n\geq 1$.
This is equal to the smallest equivalence relation on $PX$ containing $R$ 
such that for all $u,v \in PX \ x=_Ry \Rightarrow uxv =_R uyv$ 
\cite{Cremanns}.
If $A$ is a monoid and $=_S$ a congruence on $A$ then the 
\textbf{factor monoid} $A/=_S$ is the monoid whose elements are the 
congruence classes of $=_S$ on $A$ and whose composition is induced by that on 
$A$. The congruence class of an element $a\in A$ with respect to $S$ will be 
denoted $[a]_S$.\\

  A \textbf{monoid presentation} is a pair $mon\lan X | R \ran$, where $X$ is a 
set and $R\subseteq PX\times PX$ is a set of relations.  
The monoid it presents is the factor monoid $PX / =_R$. We say 
$mon\lan X | R \ran$ is a monoid presentation of $M$ if $M \cong PX / =_R$.
  The \textbf{free group} 
on $X$ is the group $F(X)$ with monoid presentation
$mon\lan \bar{X} | R_0 \ran$ where $\bar{X}:=\{x^+,x^- : x \in X \}$ and
$R_0:=\{(x^+x^-,\id),(x^-x^+,\id):x\in X\}$. 
  A \textbf{group presentation} 
is a pair $grp\lan X | R\ran $ where $X$ is a set and $R \subseteq F(X)$ 
(the group \emph{relators}). The group it presents is defined as the monoid 
that is presented by $mon\lan \bar{X} | \bar{R} \ran$ where 
$\bar{R}:=R_0 \cup \{(r,\id):r\in R\}$.
(To verify that this is a group note that any element has the form 
$[{x_1}^{\ep_1}\ldots{x_n}^{\ep_n}]_{\bar{R}}$ where $x_1,\ldots,x_n \in X, \,
\ep_1,\ldots,\ep_n \in \{+,-\}$ and so has inverse 
$[{x_n}^{-\ep_n}\ldots{x_1}^{-\ep_1}]_{R'}$ where $-(+):=-,-(-):=+$.)\\

  A monoid is \textbf{finitely presented} 
if it has a presentation $mon\lan X \,|\, R\ran$ where $X$ and $R$ are 
finite sets (similarly for groups).
  Monoid presentations are often used to give all the information about
the monoid in a compact form. The main question, given a monoid
presentation, is known as the word problem.
  The \textbf{word problem for a monoid presentation} 
$mon\lan X|R\ran$ is as follows:
\begin{center}
\begin{tabular}{lll}
INPUT: & $u, v \in PX$ & (two elements in the free monoid),\\
QUESTION: & $u =_R v$? & (do they represent the same element in the monoid 
presented?)\\
\end{tabular}
\end{center}
  Rewriting systems (defined later) are one method of tackling this problem
(another being the Todd-Coxeter procedure). However, as is well known,
rewriting cannot solve the problem in general but only when the rewriting
system can be \emph{completed} (defined later). Fortunately there are a large 
number of interesting examples (all finite monoids, all abelian monoids - 
see later) for which rewriting systems are completable.

\subsection{Category and Groupoid Presentations}
  
  It is assumed that the reader is familiar with the general concepts
of category, functor and natural transformation. The following paragraphs 
fix the notation used and define presentations of categories and groupoids 
and the associated word problem.\\

  A \textbf{directed graph} $\Gamma$ consists of a set of objects 
$\ob\Gamma$, a set of arrows $\arr\Gamma$ and two functions 
$src,tgt:\arr\Gamma \to \ob\Gamma$.
(Throughout the text, unless otherwise specified, ``graph'' should be taken 
to mean such a directed graph. If a graph has only one object this will
be denoted $\bullet$.) A \textbf{morphism of graphs} $F:\Gamma\to\Delta$ 
consists of functions $\ob F:\ob\Gamma\to\ob\Delta$, 
$\arr F:\arr\Gamma \to \arr\Delta$ such that 
$src\circ \arr F=\ob F \circ src$ and $tgt\circ\arr F=   \ob F \circ tgt$. 
This gives the category $\mathsf{DirG}$ of directed graphs.\\ 
 
The forgetful functor $U:\mathsf{Cat}\to\mathsf{DirG}$ from the category of 
small categories to directed Graphs has a left adjoint which we write $P$, 
the \textbf{free category} on a graph. It is realised in the usual way: 
if $\Gamma$ is a graph then $\ob P\Gamma:=\ob\Gamma$, and the non-identity 
arrows $P\Gamma(A_1,A_2)$ consist of all paths $a_1\cdots a_n$, i.e. sequences 
$a_1,\ldots,a_n \in \Gamma$ such that 
$tgt(a_i)=src(a_{i+1})$ for $i=1,\ldots,n-1,\ n\geq 1$. 
The identity arrows are such that for all objects $A$ of the free category
$\id_A a = a$ for any path $a$ with source $A$ and
$c\,\id_A = c$ for any path $c$ with target $A$. 
Composition is defined by concatenation. 
Thus if $\Gamma$ has one object then
$P\Gamma$ can be identified with the free monoid on $\arr\Gamma$.\\

  A set of \textbf{relations} $R$ for a category $\bA$ is a subset of 
$\arr \bA \times \arr \bA$, every relation $(l,r) \in R$ must satisfy 
$src(l)=src(r)$, $~ tgt(l)=tgt(r)$.
  A \textbf{congruence} $=_S$ on a category $\bA$ is an equivalence 
relation on the set $\arr\bA$ which satisfies 
$l =_S r \Rightarrow src(l)=src(r),tgt(l)=tgt(r)$ 
and for all $u,v \in \arr\bA$, 
if $l =_S r$ then $ulv =_S urv$ 
when these products are defined.   
  The \textbf{congruence $=_R$ generated by $R$ on $P\Gamma$}, where $R$ is a 
set of relations, is given by
$x =_R y$ if there is a system of equations
\begin{eqnarray*} x & = & u_1 l_1 v_1 \\
        u_1 r_1 v_1 & = & u_2 l_2 v_2 \\
              \cdots&\cdots&\cdots \\
        u_n r_n v_n & = & y 
\end{eqnarray*} 
where either $(l_i, r_i)$ or $(r_i, l_i) \in R$ for $i=1,..,n$ and the products 
$u_il_iv_i$ and $u_ir_iv_i$ are defined.
If $\bA$ is a category and $=_S$ is a congruence on $\bA$ then the 
  \textbf{factor category} $\bA/=_S$ is the category whose objects are
$\ob\bA$ and whose arrows are the congruence classes with respect
to $=_S$ of $\arr\bA$ with composition induced by that on $A$. The congruence class of an arrow $a\in\bA$ with respect 
to $S$ will be denoted $[a]_S$. Congruent arrows have the same sources and 
targets
as each other, so $src,tgt$ are preserved.\\

  A \textbf{category presentation} is a pair $cat\lan \Gamma|R \ran $, where 
$\Gamma$ is a graph and $R\subset \arr P\Gamma\times\arr P\Gamma$ is a set of 
relations. 
The category it presents is the factor category $P\Gamma / =_R$. 
We say that $cat\lan \Gamma|R \ran $ is a category presentation for $\bC$ if 
$\bC \cong P\Gamma / =_R$.\\

  The \textbf{free groupoid} on $\Gamma$ is denoted $F(\Gamma)$. It is 
defined to be the free category $P\bar{\Gamma}$ factored by the relations $R_0$
where $\ob\bar{\Gamma}:=\ob\Gamma$, 
$\arr\bar{\Gamma}:=\{a^+,a^-:a\in \arr\Gamma\}$ with $src(a^+)=tgt(a^-)=src(a)$
and $tgt(a^+)=src(a^-)=tgt(a)$ and  
$R_0:=\{(a^+a^-,\id_{src(a)}),(a^-a^+,\id_{tgt(a)}) : a \in \arr\Gamma \}$.
  A \textbf{groupoid presentation} is a pair $gpd \lan \Gamma | R \ran$ where
$\Gamma$ is a graph and $R$ is a subset of the disjoint union of the vertex 
groups of $F(X)$.
The groupoid it presents is defined as the category that is presented by
$cat\lan \bar{\Gamma} | \bar{R} \ran$ where $\bar{\Gamma}$ and $R_0$ are as above 
and $\bar{R}:=R_0 \cup \{(r,\id_{src(r)} : r \in R\}$.
(To verify that this is a groupoid note that any element has the form 
$[{a_1}^{\ep_1}\cdots{a_n}^{\ep_n}]_{\bar{R}}$ where $a_1,\ldots,a_n \in \Gamma, \,
\ep_1,\ldots,\ep_n \in \{+,-\}$ and so has inverse 
$[{a_n}^{-\ep_n}..{a_1}^{-\ep_1}]_{\bar{R}}$ where $-(+):=-,-(-):=+$.)\\

Some motivation for considering groupoid presentations is given by the fact that
a presentation $grp\lan X | R \ran$ of a group $G$ 
lifts to a presentation $gpd\lan \widetilde{X} | \widetilde{R} \ran$ of the 
\emph{covering groupoid} of the Cayley graph $\widetilde{X}$ of the group $G$
 \cite{Higgins64}.
In detail: let $\theta:F(X) \to G$ be the quotient map, and let
$\ob\widetilde{X}=\{g:g \in G\}$, 
$\arr\widetilde{X}=\{[g,x]:g \in G, x \in X\}$
where $src([g,x]):=g, \, tgt([g,x]):=g\theta(x)$,
and $\widetilde{R}=G \times R$.
(This is referred to in detail in Chapter 5).
  A monoid (or group) can be regarded as a category (or groupoid) with one 
object.  Let $mon\lan X|R \ran$ present a monoid $M$. Then the presentation
$cat\lan \Gamma X|R \ran$, where $\Gamma X$ is the one object graph and 
$\arr\Gamma X:=X$, is a category presentation for the monoid $M$.\\

  A category $\bC$ is \textbf{finitely presented} if it has a presentation 
$cat\lan \Gamma | R \ran$ where $\ob\Gamma, \arr\Gamma$ and $R$ are finite sets.
  The \textbf{word problem for a category presentation} $cat\lan \Gamma|R\ran $ 
is as follows:
\begin{center}
\begin{tabular}{lll}
INPUT: & $u, v \in \arr(P\Gamma)$ & (two arrows in the free category),\\
QUESTION: & $u =_R v$? & (do they represent the same element in the category 
presented?)\\
\end{tabular}
\end{center}
Terminology: The trivial category, with category presentation 
$cat\lan \bullet|\ran ~$ has only one object $\bullet$ 
and one arrow -- the identity $\id_{\bullet}$. 
The null functor maps a category to the trivial category, by mapping all 
the objects to $\bullet$ and the arrows to $\id_{\bullet}$. 
The hom-set of all arrows between two particular objects $A$ and $B$ of a 
category $\bP$ will be denoted $\bP(A,B)$.

\section{Abstract Reduction Relations}

We recall the definitions of reduction relations on abstract sets and some
of their properties. 
This is a brief exposition of the introductory material in \cite{BoOt}, 
the results stated are proved there. 
These results will be generalised to $\bP$-sets, 
where $\bP$ is a category, in Section 2.4\\ 

 Let $T$ be a set. A \textbf{reduction relation} $\to$ on a set 
$T$ is a subset of $T \times T$. 
We write $l \to r$ when $(l,r)$ is an element (rule) of $\to$.
The pair $(T,\to)$ will be called a \textbf{reduction system}.
\textbf{Reduction} is the name given to the procedure of applying rules to a 
given term to obtain another term  i.e. we \ ``reduce $t_1$ to $t_2$ in one  
step'' if $(t_1,t_2)$ is an element of the reduction relation. 
An element $t_1$ of $T$ is said to be \textbf{reducible} if there is another  
element $t_2$ of $T$ such that $t_1 \to t_2$, otherwise it is  
\textbf{irreducible}. 
The reflexive, transitive closure of a reduction relation $\to$ is denoted 
$\stackrel{*}{\to}$ 
i.e. if $t_1\to t_2\to \cdots \to t_n$ then we write 
$t_1 \stackrel{*}{\to} t_n$. 

The reflexive, symmetric, transitive closure of $\to$ is denoted 
$\stackrel{*}{\lra}$
This is the smallest equivalence relation on $T$ that contains $\to$.
The equivalence class of an element $t$ of $T$ under $\stackrel{*}{\lra}$ will be denoted $[t]$.\\

  The \textbf{word problem} for a reduction system $(T,\to)$ is:
\begin{center}
\begin{tabular}{lll}
INPUT: & $t_1,t_2 \in T$ &(two elements of $T$).\\
QUESTION: & $t_1 {\stackrel{*}{\lra}}_R t_2$ &
            (are they equivalent under ${\stackrel{*}{\lra}}_R$)?\\
\end{tabular}
\end{center} 

  Let $\to$ be a reduction relation on a set $T$.
  A \textbf{normal form} for an element $t\in T$ is an irreducible element 
  $t_N\in T$ such that $t \stackrel{*}{\lra} t_N$.
  A \textbf{set of unique normal forms}  
is a subset of $T$ which contains exactly one normal form for each equivalence
class of $T$ with respect to $\stackrel{*}{\lra}$.
  A \textbf{unique normal form function} is a 
function $N:T \to T$ whose image is a set of unique normal forms.
  One approach to solving the word problem is to attempt to choose a set of 
unique normal forms as representatives of the classes of the equivalence 
relation. Given any pair of elements, if their normal forms can be computed,
it can be seen that the elements are equivalent if and only 
if their normal forms are equal.\\

The definitions above indicate that if the irreducible elements are to be 
unique normal forms we require exactly one irreducible in each equivalence 
class. Further, if reduction is to be the unique normal form function then 
we should be able to obtain the normal
form of any element by a finite sequence of reductions.
 We consider conditions that guarantee these properties.
It is essential that equivalent elements reduce to the same irreducible.
  A reduction system $(T,\to)$ is \textbf{confluent}, if for all terms 
$t,u_1,u_2 \in T$ such that 
$t \stackrel{*}{\to} u_1 $ and $t \stackrel{*}{\to} u_2$ 
there exists an element $v \in T$ such that 
$u_1 \stackrel{*}{\to} v$ and $u_2 \stackrel{*}{\to} v$.
  The following picture illustrates the confluence condition.
  
$$
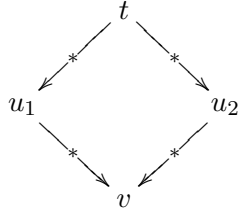
\xymatrix{ & & t \ar[dl]|{*} \ar[dr]|{*} & &\\
           & u_1 \ar[dr]|{*} &&  u_2  \ar[dl]|{*} & \\   & & v & & \\}
$$
 
The following facts may be found in \cite{BoOt}.
\begin{fact}
If a reduction system $(T,\to)$ is confluent then for each $t\in T$, $[t]$ has
at most one normal form.
\end{fact}

We require that the  irreducibles be obtainable by a finite sequence of 
reductions.
A reduction system $(T,\to)$ is \textbf{Noetherian} (or terminating) 
if there is no infinite sequence $t_1, t_2, \ldots \in T$ such that for all 
$i \in \bN, \, t_i\to t_{i+1}$.
A reduction system $(T,\to)$ is \textbf{locally confluent} if for all elements  
$t,u_1,u_2 \in T$ such that $t \to u_1 $ and $t \to u_2$ there exists  
a term $v \in T$ such that 
$u_1 \stackrel{*}{\to} v$ and 
$u_2 \stackrel{*}{\to} v$. 

\begin{fact} 
A Noetherian reduction system is confluent if it is locally confluent.
\end{fact}

\begin{fact}
If a reduction system $(T,\to)$ is Noetherian then for every
$t\in T$, $[t]$ has a normal form (not necessarily unique).
\end{fact}

A reduction system $(T,\to)$ is \textbf{complete} (or convergent) if it is
confluent and $\to$ is Noetherian. 

\begin{fact}
Let $(T,\to)$ be a reduction system. If it is complete then for every 
$t\in T$, $[t]$ has a unique normal form.
\end{fact}

Some motivation for considering complete reduction systems is that they enable the solution of the word problem through a
\textbf{normal form algorithm}. The normal forms are the irreducible elements 
(completeness ensures that there is exactly one irreducible in each equivalence 
class). The normal form function is repeated reduction (the Noetherian property
ensures that the irreducible is reached in finitely many reductions).
So: given two terms, we reduce them to irreducibles, 
the words are equivalent only if the irreducibles are equal.

\begin{fact} If a reduction system $(T,\to)$ is complete and $T$ is  finite, then the word problem for $(T,\to)$ is decidable.
\end{fact}

It is not in general possible to determine whether a finite reduction 
system is Noetherian, confluent or complete. However, if a finite system is 
known to be Noetherian, we can determine whether or not it is complete.
Non-confluence occurs when different rules apply to the same term, 
giving different reduced terms. 
A \textbf{critical pair} is a pair $(u_1,u_2)$ where there exists a term $t\in 
T$
such that $t\to u_1$ and $t\to u_2$.
A critical pair $(u_1,u_2)$ is said to \textbf{resolve} if there exists a term 
$v\in T$ such that 
$u_1 \stackrel{*}{\to} v$ and 
$u_2 \stackrel{*}{\to} v$. 

\begin{fact}
Let $(T,\to)$ be a reduction system. Let $N:T \to T$ be the normal form 
function where $N(s)$ is the irreducible form of $s$ with respect to $\to$. 
If for all $t\to s_1,t\to s_2$, 
$N(s_1)=N(s_2)$ then $(T,\to)$ is complete.
\end{fact}

A Noetherian system may sometimes be made confluent by adding in extra rules 
(the unresolvable critical pairs). This procedure will be discussed in the next
chapter in the particular setting with which we are concerned.

\chapter{Using Rewriting to Compute Kan Extensions of Actions} 
 
This chapter defines rewriting procedures for 
terms $x|w$ where $x$ is an element of a set and $w$ is a word. 
Two kinds of rewriting are involved here. 
The first is the familiar $x|ulv \to x|urv$. 
The second is given by an action of certain words on elements, so
allowing
rewriting 
$x|F(a) v \to x \cdot a|v$. 
Further, the elements $x$ and $x \cdot a$ are allowed to belong to
different
sets. 
The natural setting for this rewriting is a  ``presentation'' $kan \lan
\Gamma | \Delta | RelB | X | F \ran$ where $\Gamma,\Delta $ are
(directed) graphs and $X: \Gamma \to \sets$ and $F: \Gamma \to P \Delta$
are graph morphisms to the category of sets, and the free category on
$\Delta$ respectively, and $RelB$ is a set of relations on $P\Delta$.
The main result defines rewriting procedures on the $\bP$-set 
\begin{equation}  
T := \bigsqcup_{B\in\ob\Delta}\bigsqcup_{A\in\ob\Gamma} XA \times
\bP(FA, B)
\end{equation}
in order to
attempt the computation of Kan extensions of actions of categories given
by presentations (see section 5).

So the power of rewriting theory may now be brought to bear on a much
wider range of combinatorial enumeration problems.
Traditionally rewriting is used for solving the word problem for
monoids. It may now also be used in the specification of  
\begin{enumerate}[i)] 
\item equivalence classes and equivariant equivalence classes,  
\item arrows of a category or groupoid,  
\item action of a group on the cosets given by a subgroup,
\item right congruence classes given by a relation on a monoid, 
\item orbits of an  action of a group or monoid. 
\item conjugacy classes of a group,
\item coequalisers, pushouts and colimits of sets,
\item induced permutation representations of a group or monoid. 
\end{enumerate} 
and many others.

\section{Kan Extensions of Actions}

The concept of the Kan extension of an action will be central to this
chapter. It will therefore be defined here with some familiar examples
to 
motivate the construction listed afterwards. 
There are two types of Kan extension (the details are in
Chapter 10 of \cite{Mac}) known as right and left. Which type is right
and 
which left varies according to authors' chosen conventions. In this text
only 
one type is used (left according to \cite{CaWa1}, right according to
other 
authors) and to save conflict it will be referred to simply as ``the Kan 
extension'' - it is the colimit one, so there is an argument for calling
it a 
co-Kan, and the other one simply Kan, but we shall not presume to do
that 
here.\\

  Let $\bA$ be a category. 
  A \textbf{category action} $X$ of $\bA$ is a functor $X:\bA \to
\sets$. 
This means that for every object $A$ there is a set $XA$ and 
the arrows of $\bA$ act on the elements of the sets associated to their
sources
to return elements of the sets associated to their targets. 
So if $a_1$ is an arrow in $\bA(A_1,A_2)$  then $XA_1$ and $XA_2$ are
sets and $Xa_1 : XA_1 \to XA_2$ is a function where $Xa_1(x)$ is denoted
$x \cdot a_1$.
Furthermore, if $a_2\in\bA(A_2,A_3)$ is another arrow then 
$(x\cdot a_1)\cdot a_2=x.(a_1a_2)$ so the action preserves the
composition. 
This is equivalent to the fact that $Xa_2(Xa_1(x))=X(a_1a_2)(x)$ i.e.
$X$ is a functor. Also $F(\id_A)=\id_{FA}$ so $x\cdot \id=x$ when
defined.\\

Given the category $\bA$ and the action defined by $X$,
let $\bB$ be a second category and let $F:\bA \to \bB$ be a functor.
Then an \textbf{extension of the action $X$ along $F$} is a pair
$(K,\ep)$
where $K:\bB \to \sets$ is a functor 
and $\ep:X \to F \circ K$ is a natural transformation. 
This means that $K$ is a category action of $\bB$ and $\ep$ makes sure
that the
action defined is an extension with respect to $F$ of the action already
defined
on $A$. So $\ep$ is a collection
of functions, one for each object of $\bA$, such that $\ep_{src(a)}(Xa)$
and
$K(F(a))$ have the same action on elements of $K(F(src(a))$.\\

The \textbf{Kan extension of the action $X$ along $F$} is an extension of 
the action $(K,\ep)$ with the universal property that for any other
extension
of the action $(K',\ep')$ there exists a unique natural transformation
$\alpha:K \to K'$ such that $\ep'=\ep \circ \alpha$.
Here $K$ may thought of as the universal extension of the action of $\bA$ to
an action of $\bB$.
\vspace{-0.5cm}
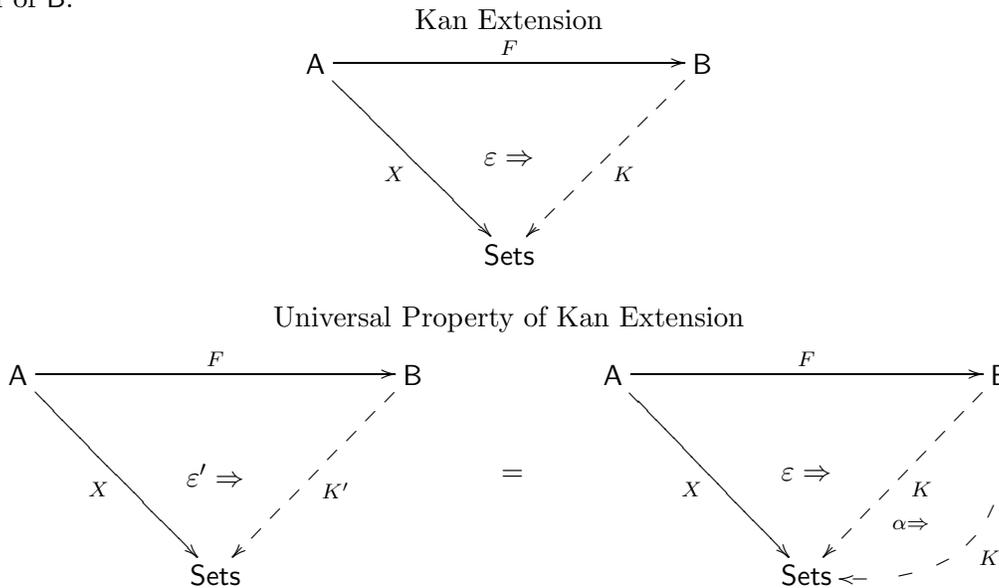
\begin{center}Kan Extension\end{center}
\vspace{-0.3cm}
$$
\xymatrix{ {\bA} \ar[rrrr]^F \ar[ddrr]_X &&&& {\bB} \ar@{-->}[ddll]^K \\ 
           && \ep \Rightarrow && \\
           && {\sets} && \\ }
$$

\begin{center}Universal Property of Kan Extension\end{center}
\vspace{-0.3cm}
$$
\xymatrix{ {\bA} \ar[rrrr]^F \ar[ddrr]_X &&&& {\bB} \ar@{-->}[ddll]^{K'}
           &&  {\bA} \ar[rrrr]^F \ar[ddrr]_X &&&& {\bB} \ar@{-->}[ddll]^K  
             \ar@/^3pc/@{-->}[ddll]^{K'} \\
           && \ep' \Rightarrow && 
           & = &  && \ep \Rightarrow &  \ar@{}[d]|{\alpha\Rightarrow}&\\
           && {\sets} &&  &&  && {\sets} && \\ }
$$
\section{Examples}
Some familiar problems will now be expressed in terms of Kan extensions. 
This is not a claim that these problems can always be computed, it
merely 
demonstrates that they are all special cases of the general problem of 
computing a Kan extension. MacLane wrote that ``the notion of Kan
extensions 
subsumes all the other 
fundamental concepts of category theory'' in section 10.7 of \cite{Mac} 
(entitled ``All Concepts are Kan Extensions''). This list helps to
illustrate 
his statement. Throughout these examples we use the same notation as the 
definition, so the pair $(K,\ep)$ is the Kan extension of the action $X$
of 
$\bA$ along the functor $F$ to $\bB$. By a monoid (or group)
``considered as 
a category'' we mean the one object category with arrows corresponding
to the 
monoid elements and composition defined by composition in the monoid.\\

\textbf{1) Groups and Monoids}\\
Let $\bB$ be a monoid regarded as a category. Let $\bA$ be the trivial
category,
acting trivially on a one point set $X\bullet$, and let $F:\bA \to \bB$
be the 
inclusion map.
Then the set $K\bullet$ is bijective with the set of elements of the
monoid and
the right action of the arrows of $\bB$ is right multiplication by the
monoid
elements. The natural transformation maps the unique element of
$X\bullet$ to 
the element of $K\bullet$ representing the monoid identity.\\

\textbf{2) Groupoids and Categories}\\
Let $\bB$ be a category. Let $\bA$ be the (discrete) category of objects
of 
$\bB$ with identity arrows only. Let $X$ define the trivial action of
$\bA$
on a collection of one point sets $\sqcup_A XA$ (one for each object $A 
\in \ob \bA$),
and let $F:\bA\to\bB$ be the inclusion map.
Then the set $KB$ for $B\in\bB$ is isomorphic to the set of arrows of
$\bB$ with 
target $B$ and the right action of the arrows of $\bB$ is defined by
right 
composition. The natural transformation maps the unique element of a set
$XA$ 
to the representative identity arrow for the object $FA$ for every
$A\in\bA$.\\

\textbf{3) Cosets, and Congruences on Monoids}\\
Let $\bB$ be a group considered as a category, and let $\bA$ be a
subgroup of $\bB$, 
with inclusion $F$. Let $X$ map the object of $\bA$ to a one point set.  
The set $K\bullet$ represents the (right) cosets of $\bA$ in $\bB$, with
the
right action of any group element $b$ of $\arr \bB$ taking
the representative of the coset $Hg$ to the representative of the coset
$Hgb$.
The left cosets can be similarly represented, defining the right action
$K$ by
a left action on the cosets. The natural transformation picks out the 
representative for the subgroup $H$. 
\\
Alternatively, let $\bB$ be a monoid considered as a category and $\bA$
be 
generated by arrows which map under $F$ to a set of generators for a
right 
congruence. Then the set $K\bullet$ represents the congruence classes,
the 
action of any monoid element $b$ of $\arr \bB$ taking
the representative (in $K\bullet$) of the class $[m]$ to the
representative 
of the class $[mb]$. The natural transformation picks out the
representative
for the class $[id]$.
(As above, left congruence classes may also be expressed in terms of a
Kan 
extension.)\\

\textbf{4) Orbits of Group Actions}\\
Let $\bA$ be a group thought of as a category and let $X$ define the
action of the group on a set $X\bullet$.
Let $\bB$ be the trivial category and let $F$ be the null functor.
Then the set $K\bullet$ is a set of representatives of the distinct
orbits of the action
and the action of $\bB$ on $K\bullet$ is trivial. The natural
transformation $\ep$ maps any 
element of the set $X\bullet$ to its orbit representative in $\bB$.\\

\textbf{5) Colimits in Sets}\\
Let $\bA$ be any category and let $\bB$ be the trivial category, with $F$
being
the null functor and $X$ being a functor to sets. Then the Kan extension 
corresponds to the colimit of (the diagram) $X:\bA\to\sets$; $K\bullet$
is the 
colimit object, and $\ep$ defines the colimit functions from each set
$XA$ to
$K\bullet$.
Examples of this are when
$\bA$ has two objects $A_1$ and $A_2$, 
and two non-identity arrows $a_1,a_2 :A_1 \to A_2$, \ 
(\emph{coequaliser} of the functions $Xa_1$ and 
$Xa_2$ in $\sets$);
$\bA$ has three objects $A_1$, $A_2$ and $A_3$
and two arrows $a_1:A_1 \to A_2$ and $a_2:A_1 \to A_3$ 
(\emph{pushout} of the functions $Xa_1$ and $Xa_2$ in $\sets$).\\

\textbf{6) Induced Permutation Representations}\\
Let $\bA$ and $\bB$ be groups thought of as categories, $F$ being a 
group morphism and $X$ being a right action of the group $\bA$ on 
the set $X\bullet$. The Kan extension of the action along $F$ is known
as the action of 
$\bB$ \emph{induced} from that of $\bA$ by $F$ (sometimes written
$F_*(X)$). 
There are simple methods of constructing the set $K \bullet$ when $\bA$
and $\bB$ are groups, but this is more difficult for monoids.\\

This last example is very close to the full definition of a Kan
extension.
A Kan extension \emph{is} the action of the category $\bB$ induced from
the
action of $\bA$ by $F$ together with $\ep$ which shows how to get from
the $\bA$-action to the $\bB$-action.
The point of the other examples is to show that Kan extensions can be
used as a 
method of representing a variety of situations.

\section{Presentations of Kan Extensions of Actions}
The problem that has been introduced is that of ``computing a Kan
extension''. In order to keep the analogy with computation and rewriting
for 
presentations of monoids we
propose the following definition of a presentation of a Kan extension. 
This formalises ideas used in \cite{CaWa2}.  

First, we define `Kan extension data'.

\begin{Def}
A \textbf{Kan extension data} $(X',F')$ consists of small categories
$\bA$, 
$\bB$ and functors $X':\bA\to\sets$ and $F':\bA\to\bB$.
\end{Def}
\begin{Def}
  A \textbf{Kan extension presentation} is a quintuple
$\mathcal{P}:=kan\lan \Gamma|\Delta|RelB|X|F \ran$ where
\begin{enumerate}[i)]
\item 
$\Gamma$ and $\Delta$ are graphs,
\item
$cat\lan \Delta | RelB \ran$ is a category presentation,
\item 
$X: \Gamma \to U \sets$ is a graph morphism,
\item
$F: \Gamma \to U P\Delta$ is a graph morphism.
\end{enumerate}

$\mathcal{P}$ \textbf{presents the Kan extension
data} $(X',F')$ where $X':\bA\to\sets$ and $F':\bA\to\bB$ 
if 
\begin{enumerate}[i)]
\item
$\Gamma$ is a generating graph for $\bA$ and $X:\Gamma\to\sets$ is the 
restriction of $X':\bA\to\sets$,
\item
$cat\lan \Delta|RelB\ran$ is a category presentation of $\bB$,
\item
$F:\Gamma\to P\Delta$ induces $F':\bA\to\bB$.
\end{enumerate}
We also say $\mathcal{P}$ \textbf{presents} the Kan extension $(K,\ep)$
of the 
Kan extension data $(X',F')$.
The presentation is \textbf{finite} if $\Gamma$, $\Delta$ and $RelB$ are finite.
\end{Def}

\begin{rem}\emph{
The fact that $X,\,F$ induce $X',\,F'$ implies extra conditions on
$X,\,F$ in 
relation to $\bA$ and $\bB$. In practice we need only the values of
$X',\,F'$ on 
$\Gamma$. This is analogous to the fact that for coset enumeration of a
subgroup 
$H$ of $G$ where $G$ has presentation $grp\lan \Delta|R\ran$ we need
only that 
$H$ is generated by certain words in the set $\Delta$.
}\end{rem}

\section{$\bP$-sets}

In this section we extend some of the usual concepts and terminology 
of rewriting in order to apply them to the new situation.
 
\begin{Def}
\label{Pset}
For a category $\bP$, a
\textbf{$\bP$-set} is a set $T$ together with a function 
$\tau:T \to \ob\bP$ and a partial action $\cdot$ of the arrows of $\bP$
on $T$. 
The action $t\cdot p$ is defined  for $t\in T$, $p\in\arr\bP$ when
$\tau(t)=src(p)$ and satisfies
\begin{align*} 
& i) \,  \tau(t\cdot p)=tgt(p),\\
\intertext{Further, for all $t\in T$, 
$p,q \in \arr\bP$ such that $(t\cdot p)\cdot q$ is
defined the following properties hold}
& ii) \,  t \cdot \id_{\tau(t)} = t,\\
& iii) \,(t\cdot p)\cdot q=t \cdot (pq).\\
\end{align*}
\end{Def}

\begin{Def}
\label{redrel}
A \textbf{reduction relation} on a $\bP$-set $T$ is a relation $\to$ on
$T$ such that for all $t_1,t_2 \in T$, $t_1 \to t_2$ implies
$\tau(t_1)=\tau(t_2)$.
\end{Def}

\begin{Def}
A reduction relation $\to$ on the $P$-set $T$ is \textbf{admissible} 
if for all $t_1, t_2 \in T$, $t_1 \to t_2$ implies 
$t_1 \cdot q \to t_2 \cdot q$ for all $q \in \arr\bP$ such that
$src(q)=\tau(t_1)$.
\end{Def}

For the rest of this chapter we assume that 
$\mathcal{P}=kan \lan \Gamma | \Delta | RelB |X |F \ran$ is a
presentation of a Kan extension.
The following definitions will be used throughout. 
Let $\bP$ denote the free category $P\Delta$.
Then define  
\begin{equation}  
T := \bigsqcup_{B\in\ob\Delta}\bigsqcup_{A\in\ob\Gamma} XA \times
\bP(FA, B)
\end{equation}
It is convenient to write an element $(x,p)$ of $XA \times \bP(FA,B)$ 
as $x|p$, a kind of ``tagged word'' -- with $x$ being the tag and $p$ the
word.
The function $\tau:T\to\ob\bP$ is defined by 
$$\tau(x|p):=tgt(p) \text{ for } x|p \in T.$$
The action of $\bP$ on $T$ is given by right multiplication
$$x|p \cdot q := x|pq \text{ for } x|p \in T, \ q \in \arr\bP \text{
when } src(q)=\tau(x|p).$$
It is routine to verify that
$\tau(x|p\cdot q)=tgt(q)$ and $(x|p\cdot q)\cdot r=(x|p)\cdot(qr)$,
whenever these terms are defined, hence proving the following lemma.

\begin{lem}
$T$ is a $\bP$-set.
\end{lem}

Now we define some `rewriting procedures' which require two types of
rule.\\

The first type is the `$\ep$-rules' $R_\ep\subseteq T\times T$. 
They are to ensure that the action is an extension of the action of
$\bA$ 
-- this is the requirement for $\ep:X \to KF$ to be a natural
transformation. 
For each arrow $a:A_1\to A_2$ in $\Gamma$ we get a set of $\ep$-rules. 
In this set there is one rule for each element $x$ of $XA_1$. Formally
\begin{equation}
R_\ep:=\{ (x|Fa,x\cdot a|id_{FA_2}) | x\in XA_1, a \in \Gamma(A_1, A_2),
A_1, A_2 \in \ob \Gamma \}.
\end{equation}
The other type is the `$K$-rules' 
$R_K \subseteq \arr \bP \times \arr \bP.$
They are to ensure that the action  preserves the structure of $\bB$ 
-- this is the requirement for $K$ to be a functor/category action.
These are simply the relations $(l,r)$ of $\bB$, formally:
\begin{equation}
R_K:=RelB.
\end{equation}
Now define $R_{init}:=(R_\ep, R_K)$. 
This we call the \textbf{initial rewrite system} that results from the
presentation.
A \textbf{rewrite system} for a Kan presentation $\mathcal{P}$ is a pair
$R$ of sets $R_T$, $R_P$ where $R_T \subseteq T \times T$ and $R_P
\subseteq \arr\bP \times \arr\bP$ such that for all $(s,u)\in R_T$, 
$\tau(s)=\tau(u)$ and for all $(l,r) \in R_P$, $src(l)=src(r)$ and 
$tgt(l)=tgt(r)$.

\begin{Def}
The \textbf{reduction relation generated by} a rewrite system $R=(R_T, R_P)$ on the $\bP$-set
$T$ is defined as 
$t_1 \to_R t_2$ if and only if one of the following is true:
\begin{enumerate}[i)]
\item 
There exist $(s,u)\in R_T, q \in \arr \bP$
such that $t_1= s \cdot q$ and $t_2:= u \cdot q$.
\item 
There exist $(l,r) \in R_P$, $s \in T$, $q \in \arr \bP$ such that $t_1
=s \cdot l q$ and $t_2= s\cdot rq$.
\end{enumerate}
Then we say $t_1$ \textbf{reduces} to $t_2$ by the rule $(s,u)$ or by $(l,r)$
respectively.
\end{Def}

Note that $\to_R$ is an admissible reduction relation on $T$ -- the proof
of this is part of the next lemma. The relation $\stackrel{*}{\to}_R$ is
the reflexive, transitive closure of $\to_R$, and $\stackrel{*}{\lra}_R$
is the reflexive, symmetric, transitive closure of $\to_R$.

\begin{rem}\emph{
Essentially, the rules of $R_P$ are two-sided and apply to any substring
to 
the right of the separator $|$. This distinguishes them from the one-
sided 
rules of $R_T$. The one-sided rules are not simply `tagged rewrite
rules' (tags being the part to the left of $|$) because the tags are
being rewritten.}
\end{rem}

\begin{lem} 
Let $R$ be a rewrite system on a $\bP$-set $T$. Then
$\stackrel{*}{\lra}_R$ is an admissible equivalence relation on the
$\bP$-set $T$.
\end{lem}
\begin{proof} 
By definition ${\stackrel{*}{\lra}}_R$ is symmetric, reflexive and
transitive. 
Now let $t_1,t_2 \in T$ be such that $t_1 \to_R t_2$ and let $v \in \arr \bP$.
be such that $src(v)= \tau(t_1)$.
Then there are two possibilities. 
  For the first case suppose 
(i) there exist $(s,u) \in R_T, q \in \arr \bP$
such that $t_1 = s \cdot q$ and $t_2 = u \cdot q$. 
Then it follows that 
$t_1 \cdot v = s \cdot q v$ and 
$t_2 \cdot v = u \cdot q v$,
(by $\bP$-set properties).
For the second case suppose (ii) there exist $s \in T$,
$(l_1,r_1) \in R_P$, $q \in \arr \bP$ such that 
$t_1 = s \cdot l q$ and 
$t_2 = s \cdot r q$.
Then it follows that
$t_1 \cdot v = s \cdot l q v$ and 
$t_2 \cdot v = s \cdot r q v$.
In either case $t_1 \cdot v \to_R t_2 \cdot v$ 
by the definition of $\to_R$.
Therefore $\to_R$ is admissible, and hence 
${\stackrel{*}{\lra}}_R$
is admissible.
\end{proof}
      
Notation: the equivalence class of $t \in T$ under $\stackrel{*}{\lra}_R$ will
be denoted $[t]$.\\

A Kan extension $(K, \ep)$ is given by a set $KB$ for
each $B \in \ob \Delta$ and a function $Kb: KB_1 \to KB_2$ for each
$b:B_1 \to B_2 \in \bB$, (defining the functor $K$) together with a
function $\ep_A: XA \to KFA$ for each $A \in \ob \bA$ (the natural
transformation). This information can be given in four parts: the set
$\sqcup KB$, a function $\bar{\tau}: \sqcup KB \to \ob \bB$, a partial
function 
(action) $\sqcup KB \times \arr\bP \to \sqcup KB$ and a function $\ep:
\sqcup XA \to \sqcup KB$. Here $\sqcup KB$ and $\sqcup XA$ (by a small
abuse of notation) are the
disjoint unions of the sets $KB$, $XA$ over $\ob \bB$, $\ob \bA$
respectively; $\bar{\tau}(z)=B$ for $z \in KB$ and if $src(p)=B$
for $p \in \arr\bP$ then $z \cdot p$ is defined.

\begin{thm}
Let $\mathcal{P}=kan \lan \Gamma | \Delta | RelB | X F \ran$ be a Kan
extension presentation, and let $\bP$, $T$, $R_{init}=(R_\ep, R_K)$ be
defined as above.
Then the Kan extension $(K,\ep)$ presented by $\mathcal{P}$ is given by
the following data:
\begin{enumerate}[i)]
\item
the set $\sqcup KB = T/ \stackrel{*}{\lra}_R$,
\item
the function $\bar{\tau}: \sqcup KB \to \ob \bB$ induced by $\tau:T \to
\ob \bP$,
\item
the action of $\bB$ on $ \sqcup KB$ induced by the action of $\bP$ on
$T$,
\item
the natural transformation $\ep$ determined by $x \mapsto [x| \id_{FA}]$
for $x \in XA$, $A \in \ob \bA$.
\end{enumerate}
\end{thm}

\begin{proof}
The initial rules $R$ on $T$ generate a reduction relation $\to$ on $T$.
Let $\stackrel{*}{\lra}$ denote the reflexive, symmetric, transitive
closure of $\to$.
 
\mbox{ }\textbf{Claim} 
$\stackrel{*}{\lra}$ preserves the function $\tau$.

\mbox{ }\textbf{Proof} 
Let $[x|p]$ denote the class of elements equivalent under
$\stackrel{*}{\lra}$ to $x|p\in T$.
We prove that $\lra$, the symmetric closure of $\to$ preserves $\tau$.
Let $t_1,t_2\in T$ so that  $t_1 \lra t_2$.
>From the definition of $\to$ there are two possible situations.
For the first case suppose that there exist
$(s_1,s_2) \in R_\ep$ such that $t_1 = s_1 \cdot p$ and $t_2 = s_2 \cdot
p$ for some $p \in \arr \bP$. 
Clearly $\tau(t_1)=\tau(t_2)$.
For the other case suppose that there exist
$(l,r)\in R_K$ such that $t_1=s \cdot(lp)$ and $t_2 = s\cdot(rp)$ for
some
$s \in T$, $p \in \arr \bP$. Again, it is clear that
$\tau(t_1)=\tau(t_2)$.
Hence $\bar{\tau}:T/\stackrel{*}{\lra}_R \; \to \ob\bP$ is well-defined
by $\bar{\tau}[t]=\tau(t)$. \hfill $\Box$\\

\mbox{ }\textbf{Claim} 
$T/\stackrel{*}{\lra}$ is a $\bB$-set.

\mbox{ }\textbf{Proof}
First we prove that $\bB$ acts on the equivalence classes of $T$ with
respect to 
$\stackrel{*}{\lra}$. 
An arrow of $\bB$ is an equivalence class $[p]$ of arrows of $\bP$ with
respect
to $RelB$.
It is required to prove that $[t]\cdot p:=[t\cdot p]$ is a well defined
action
of $\bP$ on $T/\stackrel{*}{\lra}$ such that $[t]\cdot p=[t]\cdot q$ for
all 
$p=_{RelB}q$.
Let $t\in T, p\in\arr\bP$ be such that $\tau[t]=src[p]$ i.e.
$\tau(t)=src(p)$.
Then $t\cdot p$ is defined. 
Suppose $s \stackrel{*}{\lra} t$. Then $[s\cdot p]=[t\cdot p]$ since 
$s\cdot p \stackrel{*}{\lra} t\cdot p$, whenever $s\cdot p, t\cdot p$
are 
defined.
Suppose $p =_{RelB} q$. Then $[t\cdot p]=[t\cdot q]$ since
$t\cdot p {\stackrel{*}{\lra}}_{R_K} t\cdot q$, whenever $t\cdot p,
t\cdot q$ 
are defined 
and ${\stackrel{*}{\lra}}_{RelB}$ is contained in $\stackrel{*}{\lra}$.
Therefore $\bP$ acts on $T/\stackrel{*}{\lra}$ and this action preserves
the 
relations of 
$\bB$ and so defines an action of $\bB$ on $T/\stackrel{*}{\lra}$.
Furthermore $\bar{ \tau }( [ t ] \cdot p ) = \bar{\tau}[t \cdot p ] =
tgt( p )$ and if $q \in \bP$ such that $src(q)=tgt(p)$ then 
$( [ t ] \cdot p ) \cdot q = [ ( t \cdot p ) \cdot q ] = [ t \cdot ( pq
) ] = [ t ] \cdot pq$. \hfill $\Box$\\

  The Kan extension may now be defined. For $B\in\ob\bB$ define 
\begin{equation}
\label{Kobjects}
KB:=\{ [x|p] : \bar{\tau}[x|p] = B \}.
\end{equation}
For $b:B_1\to B_2$ in $\bB$ define 
\begin{equation}
\label{Karrows}
Kb:KB_1\to KB_2 : [t] \mapsto [t\cdot p] \text{ for }[t]\in KB_1\text{
where } p \in [b].
\end{equation}
It is now routine to verify, since $p_1 =_{RelB} p_2$ implies $t \cdot p_1 
\stackrel{*}{\lra}_R t \cdot p_2$, for all $t$ where $t cdot p_1$ is defined,
that this definition of the action is a functor 
$K:\bB \to \sets$. 
Then define
\begin{equation}
\label{nattran}
\ep:X\to KF:x\mapsto[x|\id_{FA}]\text{ for }x\in XA,A\in\ob\bA.
\end{equation}
It is straightforward to verify that this is a natural transformation since 
$x|\id_{FA_1} \cdot Fa \stackrel{*}{\lra}_R x \cdot a| \id_{FA_2}$ 
for all $x \in XA_1$, $a:A_1 \to A_2 \in \ob \bA$. 

Therefore $(K,\ep)$ is an extension of the action $X$ of $\bA$.
The proof of the universal property of the extension is as follows.  
Let $K':\bB\to\sets$ be a functor 
and $\ep':X \to K'F$ be a natural transformation. 
Then there is a unique natural transformation $\alpha:K \to K'$, defined
by
$$\alpha_B[x|p] = K'(f) (\ep_A'(x))  \text{ for } [x|p]\in KB,$$
which clearly satisfies $\ep \circ \alpha = \ep'$.
\end{proof}
 
\begin{rem}\emph{
If the Kan extension presentation is finite then $R$ is finite. The
number of initial rules is by definition
$(\Sigma_{a\in\arr\Gamma}|Xsrc(a)|)+|Rel\bB|$. }
\end{rem}

\section{Rewriting Procedures for Kan Extensions}

In the next section we will explain the completion process for the
initial rewrite system. It is convenient for this procedure to have a
notation for the implementation of the data structure for a
\emph{finite} presentation $\mathcal{P}$ of a Kan extension. This we do
here.

\subsection{Input Data}

\begin{enumerate}[]
\item 
$\mathtt{ObA}$ \quad This is a list of integers $[1,2,\ldots]$, 
where each entry $i$ 
corresponds uniquely to an object $A_i$ of $\Gamma$.
\item 
$\mathtt{ArrA}$ \quad This is a list of pairs of integers 
$[[i_1,j_1],[i_2,j_2],\ldots]$, one for each
arrow $a_k : A_{i_k} \to A_{j_k}$ of $\arr\Gamma$. The first element of
each pair is the 
source of the arrow it represents, and the other entry is the target.
\item 
$\mathtt{ObB}$ \quad Similarly to $\ob\Gamma$, this is a list of
integers 
representing the objects of $\Delta$.
\item
$\mathtt{ArrB}$ \quad This is a list of triples 
$[[b_1,i_1,j_1],[b_2,i_2,j_2],\ldots]$, one triple for each arrow $b_k :
B_{i_k} \to B_{j_k}$ of 
$\mathtt{\arr\Delta}$. The first entry of each triple is a label for the
arrow (in $\mathsf{GAP}$ this is called a generator), 
and the other entries are integers representing the source and target
respectively. Note that the arrows of $\Gamma$ did not have labels. The
arrows of $\Delta$ will form parts of the terms of $T$ whilst those of
$\Gamma$ do not, so this is why we have labels here and not before. 
\item
$\mathtt{RelB}$ \quad This is a finite list of pairs of paths. Each path
is 
represented by a finite list $[b_1,b_2,\ldots,b_n]$ of labels of
composable arrows of $\mathtt{\arr\Delta}$.
In $\mathsf{GAP}$ it is convenient to consider these lists as
words $b_1 \cdots b_n$ in the generators that are labels for the arrows
of $\Delta$.
\item
$\mathtt{FObA}$ \quad This is a list of $|\ob\Gamma|$ integers.
The $k$th entry represents the object of $\Delta$ which is the image of
the 
object $A_k$ under $F$.
\item
$\mathtt{FArrA}$ \quad This is a list of paths where the entry at the
$k$th position is the path of $\bP$
which is the image of the arrow $a_k$ of $\Gamma$ under $F$. The length of the list is $|\arr
\Gamma|$.
\item 
$\mathtt{XObA}$ \quad This is a list of lists of distinct (GAP)
generators. 
There is one list of elements for each object in $\Gamma$. The list at
position 
$k$ represents the set which is the image of $A_k$ under $X$.
\item
$\mathtt{XArrA}$ \quad This is a list of lists of generators.
There is 
one list for each arrow $a$ of $\Gamma$. It represents the image under
the 
action
$Xa$ of the set $X(src(a))$. Suppose $a_k : A_{i_k} \to A_{j_k}$ is the
arrow at
entry $k$ 
in $\arr\Gamma$, and $[x_1,x_2,\ldots,x_m]$ is the $i$th entry in
$X \ob \Gamma$ 
(the image set $X(A_i)$). Then the $k$th entry of $X \arr \Gamma$ is the
list
$[x_1\cdot a,x_2\cdot a,\ldots,x_m\cdot a]$ where $x_i\in X(A_j)$.
\end{enumerate} 

Note: All the above lists are finite since the Kan extension is finitely
presented.

\subsection{Initial Rules Procedure}

The programmed function $\mathtt{InitialRules}$ extracts from the above
data the initial rewrite system $R_{init} := (R_\ep, R_K)$.
\begin{verbatim}
INPUT:     (ObA,ArrA,ObB,ArrB,RelB,FObA,FArrA,XObA,XArrA);
PROCEDURE: ans:=RelB;
           i:=1;
           while(i>Length(ArrA)) do
               a:=ArrA[i];                  ## arrow
               A:=a[1];                     ## source
               XA:=XObA[Position(ObA,A)];   ## set
               for j in [1..Length(XA)] do
                   x:=XA[j];                ## element
                   xa:=XArrA[i][j];         ## element after action
                   Fa:=FArrA[i][j];         ## image of arrow
                   rule:=[[x,Fa],[xa]];     ## epsilon-rule
                   Add(ans,rule);
               od;
           i:=i+1;
           od;
OUTPUT:    R:=ans;                          ## initial rewrite system
\end{verbatim}

We continue with the notation introduced so far, and apply the standard
terminology of reduction relations to the reduction relation $\to_R$ on
$T$.

\subsection{Lists}

In our $\mathsf{GAP}$ implementation terms of $T$ are represented by words in
generators, the generators may be thought of as labels, and the words as
lists. The first entry in the list must be a label for an element of
$XA$ for some $A \in \ob \Gamma$. The following entries will be labels
for composable arrows of $\Delta$, with the source of the first being
$FA$. Formally:\\

Let $L$ be the set of lists $l=\mathtt{[x,b1,\ldots,bn]}$, $n \geq 1$, 
such that $p = b_1 \cdots b_n$ is a reduced path (i.e. with no identity 
arrows) of $\bP$ and $x|p \in T$ or $l=\mathtt[x]$ and $x|\id_{\tau(x)} 
\in T$. We will refer to $\mathtt{List}(t)$ as the unique list associated 
with the element $t \in T$.
We will make use of the computer notation to extract particular 
elements of the list. So $t[1]$ means the first element $x$ when
$t=x|b_1 \cdots b_n$
and
$t[2..5]$ is the sublist which is $[b_1,\ldots,b_4]$ in the example,
which is an arrow in $\bP$.
Also, $\mathtt{Length}(t)$ means the number of elements in the list $t$. 
A sublist of the list for a tagged string $t \in T$ will be referred to as
a \emph{part} of $t$.

\subsection{Orderings}

To work with a rewrite system $R$ on $T$ we will require certain
concepts of order on $T$. We show how to use an ordering $>_X$ on
$\sqcup XA$
together
with an ordering $>_P$ on $\arr\bP$, these having certain
properties, to   
construct an ordering $>_T$ on $T$ with the properties needed for the
rewriting procedures.

\begin{Def}
A binary operation $>$ on the set is called a \textbf{strict partial
ordering} if it is irreflexive, antisymmetric and transitive.
\end{Def}

\begin{Def}
Let $>_X$ be a strict partial ordering on the set $\sqcup XA$. It is
called a 
\textbf{total ordering} if for all $x,y\in\sqcup XA$ 
either $x>_Xy$ or $y>_Xx$ or else $x=y$.
\end{Def}

\begin{Def}
Let $>_P$ be a strict partial ordering on $\arr \bP$. It is called
a \textbf{total path ordering} if for all $p, q \in \arr\bP$ such
that 
$src(p)=src(q)$ and $tgt(p)=tgt(q)$ either $p>_Pq$ or $q>_Pp$ or else
$p=q$.
\end{Def}

\begin{Def}
The ordering $>_P$ is \textbf{admissible on $\arr\bP$} 
if $p>_Pq \ \Rightarrow \ upv>_Puqv$ for all $u,v\in\arr\bP$ such that 
$upv,uqv\in\arr\bP$.
\end{Def}

\begin{Def}
An ordering $>$ is \textbf{well-founded} on a set of elements if there
is no
infinite sequence $x_1>x_2>\cdots$.
An ordering $>$ is a \textbf{well-ordering} on a structure if it is 
well-founded and a total ordering with respect to that structure.
\end{Def}

\begin{lem}
Let $>_X$ be a well-ordering on the finite set $\sqcup XA$ and let $>_P$
be an 
admissible well-ordering on $\bP$.
For $t_1,t_2\in T$ define $t_1 >_T t_2$ if
$$ 
 \Leftarrow t_1[2..Length(t_1)] >_P t_2[2..Length(t_2)]
\text{ or } t_1[2..Length(t_1)] = t_2[2..Length(t_2)] \text{ and }
 t_1[1] >_X t_2[1].
$$
Then $>_T$ is an admissible well-ordering on the $\bP$-set $T$. 
\end{lem}

\begin{proof}
It is straightforward to verify that irreflexivity, antisymmetry and 
transitivity of $>_X$ and $>_P$ imply those properties for $>_T$.
The ordering $>_T$ is admissible on $T$ because it is made compatible
with
the right action (defined by composition between arrows on $\bP$) by the
admissibility of $\>_P$ on $\arr\bP$. The ordering is linear, since if 
$t_1,t_2\in T$ such that neither $t_1>_Tt_2$ nor $t_2>_Tt_1$, it follows
by the linearity of $>_X$ and linearity of $>_P$ on $\arr\bP$ that
$t_1=t_2$.
That $>_T$ is well-founded is easily verified using the fact that any
infinite
sequence in terms of $>_T$ implies an infinite sequence in either $>_X$
or $>_P$
and $>_X$ and $>_P$ are both well-founded, so there are no such
sequences.
\end{proof}

The last result shows that there is some scope for choosing different
orderings on 
$T$. The actual choice is even wider than this but it is not relevant to
discuss this here. We are not concerned here with considering ranges of
possible 
orderings, but work with the one that is most straightforward to use.
The ordering implemented is a variation on the above. It corresponds to
the
length-lexicographical ordering and is defined in the following way.

\begin{Def}[Implemented Ordering]
Let $>_X$ be any linear order on (the finite set) $\sqcup XA$. 
Let $>_\Gamma$ be a linear ordering on (the finite set) $\arr\Delta$.
This induces an admissible ordering $>_P$ on $\arr\bP$ where $p>_Pq$ if
and only if $Length(p)>Length(q)$ or $Length(p)=Length(q)$ and there
exists 
$k>0$ such that $p[i]>_\Gamma q[i]$ for all $i<k$ and $p[k]=q[k]$.
The ordering $>_T$ is then defined as follows:
$t_1>_Tt_2$ if $Length(t_1)>Length(t_2)$ or
if $Length(t_1)=Length(t_2)$ and $t_1[1]>_Xt_2[1]$,
or if $Length(t_1)=Length(t_2)$ and there exists $k \in
[1..Length(t_1)]$ such 
that $t_1[i]=t_2[i]$ for all $i<k$ and $t_1[k]>_\Gamma t_2[k]$.
\end{Def}

\begin{prop} 
The definitions above give an admissible, length-non-increasing 
well-order $>_T$ 
on the $\bP$-set $T$.
\end{prop}

\begin{proof}
It is immediate from the definition that $>_T$ is length-non-increasing.
It is straightforward to verify that $>_T$ is irreflexive, antisymmetric
and 
transitive. It can also be seen that $>_T$ is linear (suppose neither 
$t_1 >_T t_2$ nor $t_2 >_T t_1$ then $t_1=t_2$, by the definition, and
linearity of
$>_X$, $>_\Gamma$). It is clear from the definition that $>_T$ is
admissible on 
the $\bP$-set $T$ (if $t_1 >_T t_2$ then $t_1.p >_T t_2.p$). 
To prove that $>_T$ is well-founded on $T$, suppose that 
$t_1 >_T t_2 >_T t_3 >_T \cdots$ is an infinite sequence. Then for each $i>0$
either 
$Length(t_i)>Length(t_{i+1})$ or if $Length(t_i)=Length(t_{i+1})$ and 
$t_i[1]>_Xt_{i+1}[1]$, or if $Length(t_i)=Length(t_{i+1})$ and there
exists $k \in [1..Length(t_i)]$ such that $t_i[j]=t_{i+1}[j]$ for all 
$j<k$ and $t_i[k] >_\Gamma t_{i+1}[k]$.
This implies that there is an infinite sequence of type
$n_1 > n_2 > n_3 > \cdots$ of
positive integers from some finite $n_1$, or of type 
$x_1 >_X x_2 >_X x_3 > \cdots$ of elements of $\sqcup XA$ or else of type 
$p_1 >_\Gamma p_2 >_\Gamma p_3 >_\Gamma \cdots$ of arrows of $\Delta$,
none
of which 
is possible as $>$, $>_X$, and $>_\Gamma$ are well-founded on $\bN$,
$\sqcup XA$ 
and $\arr\Delta$ respectively. Hence $>_T$ is well-founded.
\end{proof}

\begin{prop}
\label{ordarrow}
Let $>_T$ be the order defined above. Then $p_1 >_P p_2 \Rightarrow s\cdot 
p_1 >_T s \cdot p_2$.
\end{prop}
\begin{proof}
This follows immediately from the definition of $>_T$.
\end{proof}

\begin{rem}
\emph{The proposition can also be proved for the earlier definition of $>_T$
induced from $>_X$ and $>_P$.}
\end{rem}

\subsection{Reduction}

Now that we have defined an admissible well-ordering on $T$ it is
possible to discuss when a reduction relation generated by a rewrite
system is compatible with this ordering. 

\begin{lem} 
Let $R$ be a rewrite system on $T$. Orientate the rules of $R$ so that
for all $(l,r)$ in $R$, if $l,r\in\arr\bP$ then $l >_P r$ and if $l,r\in
T$ then $l >_T r$. 
Then the reduction relation $\to_R$ generated by $R$ is compatible with
$>_T$.
\end{lem}
\begin{proof}
Let $t_1,t_2\in T$ such that $t_1\to_Rt_2$.
There are two cases to be considered \ref{redrel}.
For the first case let $t_1=s_1\cdot p$, $t_2=s_2\cdot p$ for some
$s_1,s_2\in 
T$, $p\in\arr\bP$ such that $(s_1,s_2)\in R$. Then $s_1>_Ts_2$. It
follows that
$t_1>_Tt_2$ since $>_T$ is admissible on $T$.
For the second case let $t_1=s \cdot p_1q$, $t_2=s \cdot p_2q$ for some
$s \in T$, $p_1,p_2,q\in\arr\bP$ such that $(p_1,p_2)\in T$. Then
$p_1>_Pp_2$ 
and so by Proposition \ref{ordarrow} $s\cdot p_1>_Ts\cdot p_2$. Hence
$t_1>_Tt_2$ by
admissibility of $>_T$ on $T$. Therefore, in either case $t_1>_Tt_2$ so
$\to_R$ is compatible with $>_T$.
\end{proof}

\begin{rem}\emph{
A reduction is the replacement of a part of a tagged string $x|p \in T$ 
according to a rule of $R$. Rules from $R_T$ replace the tag $x|$ and part 
of the string $p$ whilst rules from $R_P$ replace substrings of $p$.
The reduction relation $\to_R$ is the successive replacement of parts of a 
tagged string.}
\end{rem}

It is a standard result that if a reduction relation is compatible with
an admissible well-ordering, then it is Noetherian. The next pseudo
program 
shows the function $\mathtt{Reduce}$ which returns from a term $t \in T$
and a rewrite system $R \subseteq T \times T \sqcup \arr \bP \times \arr
\bP$ a term $t_n \in [t]$ which is irreducible with respect to $\to_R$.

\begin{verbatim}
INPUT:(t,R);
PROCEDURE: new:=t; old:=[];
           while not(new=old) do
               old:=new;
               for rule in R do
                   lhs:=rule[1]; rhs:=rule[2];
                   if lhs is a sublist of new 
                       replace lhs in new by rhs
                   fi;
               od;
           od;  
OUTPUT: tn                # irreducible term in T #                 
\end{verbatim}

\subsection{Critical Pairs}

We can now discuss what properties of $R$ will make $\to_R$ a complete
(i.e. Noetherian and confluent) reduction relation. By standard abuse of
notation the rewrite system $R$ will be called \textbf{complete} when 
$\to_R$ is complete. In this case $\stackrel{*}{\lra}_R$ admits a normal 
form function. 

\begin{lem}[Newman's Lemma]
A Noetherian reduction relation on a set is confluent if it is 
locally confluent \cite{TAT}.
\end{lem}

Hence, if $R$ is compatible with an admissible well-ordering on
$T$ and  $\to_R$ is locally confluent then $\to_R$ is complete.
By orientating the pairs of $R$ with respect to the chosen ordering
$>_T$ on 
$T$, $R$ is made to be Noetherian. The remaining problem is testing for
local confluence of $\to_R$ and changing $R$ in order to obtain an
equivalent confluent reduction relation.\\

We will now explain the notion of critical pair for a rewrite system for
$T$, extending the traditional notion to out situation. In particular
the overlaps involve either just $R_T$, or just $R_P$ or an interaction
between $R_T$ and $R_P$.\\
 
A term $crit\in T$ is called \textbf{critical} if it may be reduced by
two or 
more different rules 
i.e. $crit \to_R crit1$, $crit \to_R crit2$ and $crit1 \not= crit2$. 
The pair $(crit1,crit2)$ resulting from two single-step 
reductions of the same term is called a \textbf{critical pair}.
A critical pair for a reduction relation $\to_R$ is said to
\textbf{resolve} if 
there 
exists a term $res$ such that both $crit1$ and $crit2$ reduce to a
common term $res$ i.e. $crit1 \stackrel{*}{\to}_R res$, 
$crit2 \stackrel{*}{\to}_R res$.

We now define overlaps of rules for our type of rewrite system, and show
how each kind results in a critical pair of the reduction relation.
Let $R=(R_T,R_P)$ be a rewrite system,
where $R_T \subseteq T \times T$ and $R_P \subseteq \arr\bP \times \arr \bP$.

\begin{Def}
Let $(rule1,rule2)$ be a pair of rules of $R$ such that
$rule1$ 
and $rule2$ may both be applied to the same term $crit$ in such a way
that there 
is a part of the term $crit$ that is affected by both the rules.
When this occurs the rules are said to \textbf{overlap}. There are five 
types of overlap for this kind of rewrite system.
\begin{align*}
\intertext{Suppose $rule1,rule2 \in R_T$.
Put $rule1:=(s_1,u_1)$, $rule2:=(s_2,u_2)$. Then there is one type of
overlap:}
i)\; 
& s_1 = s_2 \cdot q \text{ for some } q \in \arr \bP,
\text{ with resulting critical pair } 
  (u_1, u_2 \cdot q).\\
\intertext{Suppose $rule1,rule2 \in R_P$.
Put $rule1:=(l_1, r_1)$, $rule2:=(l_2,r_2)$. 
Then there are two possible types of overlap: }
ii)\;
& l_1 = p l_2 q \text{ for some } p, q \in \arr \bP,
\text{ with resulting critical pair }
  (r_1, p r_2 q).\\
iii)\;
& l_1 q = p l_2 \text{ for some } p, q \in \arr\bP,
\text{ with resulting critical pair }
  (r_1 q, p r_2).\\ 
\intertext{Suppose $rule1 \in R_T$, $rule2 \in R_P$.
Put $rule1:=(s_1,u_1)$, $rule2:=(l_1,r_1)$.
Then there are two possible types of overlap:}
iv)\;
& s_1 \cdot q = s \cdot l_1 \text{ for some } s \in T, \ q \in \arr
\bP,
\text{ with resulting critical pair } (u_1 \cdot q, s \cdot
r_1).\\
v)\;
& s_1 = s \cdot(l_1 q) \text{ for some } s \in T, \  q \in \arr \bP,
\text{ with resulting critical pair } (u_1, s \cdot r_1 q).
\end{align*}
\end{Def}

One pair of rules may overlap in more than one way, giving more than one
critical pair. For example the rules $(x|a^2ba, y|ba)$ and $(a^2,b)$
overlap 
with critical term $x|a^2ba$ and critical pair $(y|ba,x|b^2a)$ and also
with critical term $x|a^2ba^2$ and critical pair $(y|ba^2,x|a^2b^2)$.

\begin{lem}
Let $R$ be a finite rewrite system on the $\bP$-set $T$.
If $(t_1,t_2)$ is a critical pair then either the pair resolves
immediately
or there is an overlap between two rules $(rule1,rule2)$ such that 
if the critical pair $(crit1,crit2)$ resulting from that overlap
resolves
then $(t_1,t_2)$ resolves.
\end{lem}
\begin{proof}
Let $(t_1,t_2)$ be a critical pair. Then there exists a critical term
$t$
and two rules $rule1$, $rule2$
such that $t$ reduces to $t_1$ with respect to $rule1$ and to $t_2$ with 
respect to $rule2$.
There are seven cases that must be considered.\\

Suppose $rule1:=(s_1,u_1), rule2:=(s_2,u_2) \in R_T$.
Then the rules must overlap on $t$ as shown: 
$$
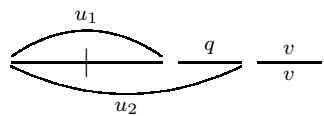
\xymatrix{\ar@{-}[rr] 
          \ar@{-}@/^1pc/[rr]^{u_1} 
          \ar@{-}@/_1pc/[rrr]_{u_2} 
        & | & \ar@{-}[r]^q 
        & \ar@{-}[r]^v \ar@{ }[r]_v &\\}
$$ 
and there exist $q, v \in\arr\bP$ such that
$t= s_1 \cdot q v = s_2 \cdot v$ and then $t_1 = u_1 \cdot q v$ and
$t_2 = u_2\cdot v$.
The critical pair resulting from this overlap (i) is 
$(u_1 \cdot q, u_2)$ and if this resolves to a common term $r$ then
$(t_1,t_2)$ 
resolves to $r \cdot v$.\\

Suppose $rule1:=(l_1,r_1)$, $rule2:=(l_2,r_2) \in R_P$.
Then there are three possible ways in which the rules may apply to $t$.
In the first case the rules do not overlap:
$$
\xymatrix{\ar@{-}[rr]^*+{s} \ar@{ }[rr]_*+{s} 
        & | 
        & \ar@{-}@/^1pc/[r]^{r_1} \ar@{-}[r]_{l_1} 
        & \ar@{-}[r]^p \ar@{ }[r]_p
        & \ar@{-}@/_1pc/[r]_{r_2} \ar@{-}[r]^{l_2}  
        & \ar@{-}[r]^q \ar@{ }[r]_q
        &\\}
$$ 
and there exist $s \in T$, $p, q \in \arr \bP$ such that
$t = s \cdot l_1 p\, l_2 q$ and then 
$t_1 = s \cdot r_1 p\, l_2 q$ and 
$t_2 = s \cdot l_1 p r_2 q$.
The pair $(t_1,t_2)$ immediately resolves to $u \cdot r_1pr_2q$ by
applying $rule2$ to $t_1$ and $rule1$ to $t_2$. 

In the second case one rule is contained within the other:
$$
\xymatrix{\ar@{-}[rr]^*+{s} \ar@{ }[rr]_*+{s}
        & | & \ar@{-}@/^1pc/[rrr]^{r_1} \ar@{-}[r]_p 
        & \ar@{-}@/_1pc/[r]_{l_2} \ar@{-}[r]
        & \ar@{-}[r]_q 
        & \ar@{-}[r]^{v} \ar@{ }[r]_v 
        &\\}
$$ 
and there exist $s \in T$, $p,q,v \in \arr \bP$ such that
$t = s \cdot l_1 v = s \cdot p\, l_2 q v$ and then 
$t_1 = s \cdot r_1 v$ and 
$t_2 = s \cdot p r_2 q v$.
The critical pair resulting from the overlap of the rules (ii) is 
$(r_1, p r_2 q)$ and if this resolves to a common term $r$ then
$(t_1,t_2)$ resolves to $s \cdot r v$. 

In the third case one part of the term is changed by both rules:
$$
\xymatrix{\ar@{-}[rr]^*+{s}  \ar@{ }[rr]_*+{s}
        & | 
        & \ar@{-}@/^1pc/[rr]^{r_1} \ar@{-}[r]_p
        & \ar@{-}@/_1pc/[rr]_{r_2} \ar@{-}[r]
        & \ar@{-}[r]^q
        & \ar@{-}[r]^{v} \ar@{ }[r]_v
        &\\}
$$ 
and there exist $s \in T$, $p, q, v \in \arr \bP$ such that
$t = s \cdot l_1 q v = s \cdot p l_2 v$ and then 
$t_1 = s \cdot r_1 q v$ and
$t_2 = s \cdot p r_2 v$.
The critical pair resulting from the overlap of the rules (iii) is 
$(r_1 q,p r_2)$ and if this resolves to a common term $r$ then
$(t_1,t_2)$ 
resolves to $s\cdot r v$.\\

Suppose finally that $rule1:=(s_1,u_1) \in R_T$ and $rule2:=(l_1,r_1)
\in R_P$.
Then there are (again) three possible ways in which the rules may apply
to $t$.
In the first case the rules do not overlap:
$$
\xymatrix{\ar@{-}[rr]_*+{s_1}
          \ar@{-}@/^1pc/[rr]^{u_1}
        & | 
        & \ar@{-}[r]^p \ar@{ }[r]_p
        & \ar@{-}@/_1pc/[r]_{r_1} \ar@{-}[r]^{l_1} 
        & \ar@{-}[r]^q \ar@{ }[r]_q
        &\\}
$$ 
and there exist $p,q\in\arr\bP$ such that
$t = s_1 \cdot p l_1 q$ and then 
$t_1 = u_1 \cdot p l_1 q$ and 
$t_2 = s_1 \cdot p r_1 q$.
The pair $(t_1,t_2)$ immediately resolves to $u_1 \cdot p r_1 q$ by
applying $rule2$ to $t_1$ and $rule1$ to $t_2$. 

In the second case one rule is contained within the other:
$$
\xymatrix{\ar@{-}[rr]_*+{s}
          \ar@{-}@/^1pc/[rrrr]^{u_1}
        & | 
        & \ar@{-}@/_1pc/[r]_{r_1} \ar@{-}[r]
        & \ar@{-}[r]_q
        & \ar@{-}[r]^v \ar@{ }[r]_v &\\}
$$ 
and there exist $s \in T$, $q, v \in \arr \bP$ such that
$t = s_1 v = s \cdot l_1 q v$ and then 
$t_1 = u_1 v$ and 
$t_2 = s r_1 q v$.
The critical pair resulting from the overlap of the rules (iv) is 
$(u_1, s \cdot r_1 q)$ and if this resolves to a common term $r$ 
then $(t_1,t_2)$ resolves to $r \cdot v$. 

In the third case there is one part of the term changed by both rules:
$$
\xymatrix{\ar@{-}[rr]_*+{s} 
          \ar@{-}@/^1pc/[rrr]^{u_1} 
        & | 
        & \ar@{-}@/_1pc/[rr]_{r_1} \ar@{-}[r]
        & \ar@{-}[r]^q 
        & \ar@{-}[r]^v \ar@{ }[r]_v 
        &\\}
$$ 
and there exist $s \in T$, $q, v \in \arr \bP$ such that
$t = s_1 \cdot q v = s \cdot l_1 v$ and then 
$t_1 = u_1 \cdot q v$ and 
$t_2 = s \cdot r_1 v$.
The critical pair resulting from the overlap of the rules (v) is 
$(u_1 \cdot q, s \cdot r_1)$ and if this resolves to a common term $r$ 
then $(t_1,t_2)$  resolves to $r \cdot v$.\\ 

Thus we have considered all possible ways in which a term may be reduced
by
two different rules, and shown that resolution of the critical pair
(when not
immediate) depends upon the resolution of the critical pair resulting 
from a particular overlap of the rules.
\end{proof}

\begin{cor} If all the overlaps between rules of a rewrite system $R$
on $T$ resolve then all the critical pairs for the reduction relation
$\to_R$
resolve, and so $\to_R$ is confluent.
\end{cor}
\begin{proof} Immediate from the Lemma.
\end{proof}

\begin{lem}
All overlaps of a pair of rules of $R$ can be found by 
looking for two types of overlap between the lists representing the left
hand sides of rules.
\end{lem}
\begin{proof}
Let $rule1=(l_1,r_1)$ and $rule2=(l_2,r_2)$ be a pair of rules. Recall
that
$\mathtt{List}(t)$ 
is the representation of a term $t \in T$ as a list.
The first type of list overlap occurs when $\mathtt{List}(l_2)$ is a
sublist of $\mathtt{List}(l_1)$ (or vice-versa). This happens in cases
(i), (ii) and (v).
The second type of list overlap occurs when the end of
$\mathtt{List}(l_1)$
matches
the beginning of $\mathtt{List}(l_2)$ (or vice-versa). This happens in
cases (iii) and (iv). 
\end{proof}

The program for finding overlaps and the resulting critical pairs is
called $\mathtt{CriticalPairs}$.
The outline of part of it is reproduced here:
Let $rule1:=(l1,r1)$ and $rule2:=(l2,r2)$ be a pair of rules. 
The program compares $rule1$ with $rule2$ to look for overlaps.
This part of the program shows how to determine whether $l1$ contains
$l2$ or the beginning of $l1$ overlaps with the end of $l2$.
To find other critical pairs the program can compare $rule2$ with
$rule1$.
\begin{verbatim}
    l1 := List(l1); len1 := Length(l1);
    l2 := List(l2); len2 := Length(l2);

    # Search for type 1 pairs  (l2 is contained in l1).
    if len1 >= len2 then
        for i in [1..len1-len2] do
        if l1{[i..i+len2-1]} = l2 then
                if i=1 then   u := IdWord;
                else          u := Product( Sublist(l1,1,i-1) );
                if i+len2-1 = len1 then    v := IdWord;
                else          v := Product( Sublist(l1,i+len2,len1) );
                [ u*r2*v, r1 ]   ## critical pair found

    # Search for type 2 pairs: (right of l1 overlaps the left of l2)
    for i in [1..len1] do
        while not( i>len1 or i>len2 ) do
            if ( l1{[len1-i+1..len1]} = l2{[1..i]} ) then
                if i = len1 then  u := IdWord;
                else              u := Product( Sublist(l1,1,len1-i) ); 
                if i = len2 then  v := IdWord;
                else              v := Product( Sublist(l2,i+1,len2) );
                [ r1*v, u*r2 ]   ## critical pair found
\end{verbatim}

It has now been proved that all the critical pairs of a finite rewrite
system $R$ on $T$ can be listed.
To test whether a critical pair resolves, each side of it is reduced
using the function $\mathtt{Reduce}$. If $\mathtt{Reduce}$ returns the
same term for each side then the pair resolves.

\subsection{Completion Procedure}

We have shown how to (i) find overlaps between rules of $R$ and 
(ii) test whether
the resulting critical pairs resolve. Further we have shown that if all
critical pairs for $R$ resolve then
$\to_R$ is confluent. We now show that critical pairs
which do not resolve may be added to $R$ without affecting the
equivalence $R$ defines on $T$.

\begin{lem} 
Any critical pair $(crit1,crit2)$ of a rewrite system $R$ 
may be added to the rewrite system without 
changing the equivalence relation ${\stackrel{*}{\lra}}_R$.
\end{lem}
\begin{proof} 
This result is proved by considering any critical pair $(t_1,t_2)$. By  
definition this pair is the result of two different single-step
reductions being
applied to a critical term $t$. Therefore $t\to_Rt_1$ and $t\to_Rt_2$.
It is
immediate that $t_1{\stackrel{*}{\lra}}_Rt{\stackrel{*}{\lra}}_Rt_2$,
and so 
adding $(t_1,t_2)$ to $R$
does not add anything to the equivalence relation
${\stackrel{*}{\lra}}$.
\end{proof}

We have now set up and proved everything necessary for a variant of the
Knuth-Bendix procedure, which will add rules to a rewrite system $R$
resulting from a presentation of a Kan extension, to attempt to find an
equivalent complete rewrite system. The benefit of such a system is that
$\to_R$ then acts as a normal form function for $\stackrel{*}{\lra}_R$
on $T$.

\begin{thm}
\label{proc}
Let $\mathcal{P}= \lan \Gamma | \Delta | RelB |X |F \ran$ be a finite
presentation of a Kan extension $(K, \ep)$. 
Let $P:=P \Delta$, 
$$T:= \bigsqcup_{B \in \ob \Delta} \bigsqcup_{A \in \ob \Gamma}
        XA \times \bP(FA, B),$$
and let $R=(R_\ep,R_P)$ be the
initial rewrite system for $\mathcal{P}$ on $T$. Let $>_T$ be an
admissible well-ordering on $T$.
Then there exists a procedure which, if it terminates,  will return a
rewrite system
$R^C$ which is complete with respect to $>_T$ such that
the admissible equivalence relations $\stackrel{*}{\lra}_{R^C}$ and 
$\stackrel{*}{\lra}_R$ coincide.
\end{thm}

\begin{proof}
The procedure finds all critical pairs resulting from overlaps of rules
of $R$. It attempts to resolve them. When they do not resolve it adds
them to the system as new rules. Critical pairs of the new system are
then examined.
When all the critical pairs of a system resolve, then the procedure
terminates, the final rewrite system $R^C$ obtained is complete.
This procedure has been verified in the preceding results of this
section.
\end{proof}

\begin{verbatim}
INPUT: (R,>T);
PROCEDURE:  NEW:=R; OLD:=[];
            while not OLD=NEW do
                CRIT:=CriticalPairs(R)
                for crit in CRIT do
                    crit[1]:=Reduce(crit[1],R);
                    crit[2]:=Reduce(crit[2],R);
                    if crit[1]=crit[2] then Remove(CRIT,crit);
                    if crit[1]<crit[2] then crit:=(crit[2],crit[1]);          
                od;
                Add(NEW,CRIT);
            od;
OUTPUT: NEW;  ## complete rewrite system.
\end{verbatim}

The whole procedure, which takes as input the presentation of a Kan
extension 
and yields as output a complete rewrite system with respect to the
ordering 
$>_T$, when this can be found, has been implemented in $\mathsf{GAP}$ in the file
$kan.g$.
We will now briefly discuss how to interpret a complete rewrite system
on $T$,
supposing that the program has returned one.

\section{Interpreting the Output}

\subsection{Finite Enumeration of the Kan Extension}

When every set $KB$ is finite we may catalogue the elements of all of
the sets
$\sqcup KB$ in stages. The 
first stage consists of all the elements $x|\id_{ FA}$ where $x \in XA$ 
for some $A \in \ob \Gamma$. These elements are considered to have
length
zero.
The next stage builds on the set of irreducible elements from the last
block 
to construct elements of the form $x|b$ where $b:FA\to B$ for some 
$B \in \ob \Delta$. This is effectively acting on the sets with the
generating 
arrows to define new (irreducible) elements of length one. The next 
builds on the irreducibles from the last block by acting with the
generators again. 
When all the elements of a block of elements of the same length
are reducible then the enumeration terminates 
(any longer term will contain one of these terms and therefore be
reducible). 
The set of irreducibles is a set of normal forms for $\sqcup KB$.
The subsets $KB$ of $\sqcup KB$ are determined by the function
$\bar{\tau}$,  i.e. if $x|b_1 \cdots b_n$ is a normal form in $\sqcup
KB$ and 
$\tau(x|b_1 \cdots b_n):=tgt(b_n)=B_n$ then $x|b_1 \cdots b_n$ is a
normal form in $KB_n$. 
Of course if one of the sets $KB$ is infinite then this may 
prevent the enumeration of other finite sets $KB_i$. The same problem 
would obviously prevent a Todd-Coxeter completion. This cataloguing
method only 
applies to finite Kan extensions. It has been implemented in the
function $kan$, which currently has an enumeration limit of 1000 on 
$\sqcup KB$ set in the program. If this limit is exceeded, the program 
returns the completed rewrite system -- provided the completion procedure
terminates.

\subsection{Regular Expression for the Kan Extension}

Let $R$ be a finite complete rewrite system on $T$ for the Kan extension 
$(K,\ep)$. Then the theory of languages and regular expressions may be 
applied. The set of irreducibles in $T$ is found after the construction
of an automaton from the rewrite system and the derivation of a language
from this automaton. Details of this method may be found in Chapter
Four.

\subsection{Iterated Kan Extensions}
One of the pleasant features of this procedure
is that the input and output are of similar form. 
The consequence of this is that if the extended action $K$
has been defined on $\Delta$ then given a second functor $G':\bB\to\bC$
and a 
presentation $cat\lan \Lambda|RelC \ran$ for $\bC$ it
is straightforward to consider a presentation for the Kan extension data 
$(K',G')$. This new extension is in fact the Kan extension with data 
$(X',F'\circ G')$ 

\begin{lem}
Let $kan\lan \Gamma|\Delta|RelB|X|F\ran$ be a presentation for a Kan
extension $(K,\ep)$.
Then let \newline $cat\lan\Lambda|RelC\ran$ present a category $\bC$ and let
$G':\bB\to\bC$. Then the Kan extension
presented by \newline
$kan\lan \Gamma|\Lambda|RelC|X|F\circ G|\ran$ is equal to the Kan
extension 
presented by
$kan\lan \Delta|\Lambda|RelC|K|G\ran$.
\end{lem}

\begin{proof}
Let $kan\lan \Gamma|\Delta|RelB|X|F\ran$ present the Kan extension data 
$(X',F')$ for the Kan extension $(K,\ep)$. 
Let $\bC$ be a category finitely presented by $cat\lan\Lambda|RelC\ran$
and let 
$G':\bB\to\bC$. Then $kan\lan \Delta|\Lambda|RelC|K|G\ran$ presents the
Kan 
extension data $(K',G')$  for the Kan extension $(L,\eta)$.

We require to prove that $(L,\ep\circ \eta)$ is the Kan 
extension presented by $kan\lan \Gamma|\Lambda|RelC|X|F\circ G\ran$
having data 
$(X',F'\circ G')$. 
It is clear that 
$(L,\epsilon\circ\eta)$ defines an extension of the action $X$ along
$F\circ G$ 
because $L$ defines an action of $\bC$ and $\ep\circ\eta:X\to F\circ
G\circ L$ 
is a natural transformation. \\
For the universal property, let $(M,\nu)$ be another extension of the
action $X$ 
along $F\circ G$. Then consider the pair $(G\circ M,\nu)$, it is an
extension of 
$X$ along $F$. 
Therefore there exists a unique natural transformation $\alpha:X \to
F\circ 
G\circ M$ such that $\ep\circ\alpha=\nu$ by universality of $(K,\ep)$.
Now consider the pair $(M,\alpha)$, it is an extension of $K$ along $G$.
Therefore there exists a unique natural transformation $\beta:L\to M$
such that 
$\eta\circ\beta=\alpha$ by universality of $(L,\eta)$.
Therefore $\beta$ is the unique natural transformation such that 
$\ep\circ\eta\circ\beta=\nu$, which proves the universality of the
extension 
$(L,\ep\circ\eta)$.
\end{proof}

\section{Example of the Rewriting Procedure for Kan Extensions} 
\label{kaneg}

Let $\bA$ and $\bB$ be the categories generated by the graphs below,
where 
$\bB$ has the relation $b_1b_2b_3=b_4$. 
{\large
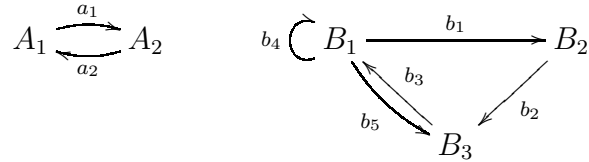
$$ 
\xymatrix{A_1 \ar@/^/[r]^{a_1}  
        & A_2 \ar@/^/[l]^{a_2}  
       && B_1  \ar@(dl,ul)^{b_4}  
             \ar[rr]^{b_1}  
             \ar@/_/[dr]_{b_5}  
       && B_2 \ar[dl]^{b_2} \\ 
     &&&& B_3 \ar[ul]_{b_3} \\} 
$$} 
Let $X:\bA \to \sets$ be defined by $ XA_1 = \{ x_1, x_2, x_3 \},  
                                   \ XA_2 = \{ y_1, y_2 \}$   with\\ 
$Xa_1:XA_1 \to XA_2: x_1 \mapsto y_1, x_2 \mapsto y_2, x_3 \mapsto
y_1$,\\ 
$Xa_2:XA_1 \to XA_2: y_1 \mapsto x_1, y_2 \mapsto x_2,$\\ 
and let $F:\bA \to \bB$ be defined by  
$FA_1=B_1, \ FA_2=B_2, \ Fa_1=b_1$ and $Fa_2 = b_3 b_2$. 
The input to the computer program takes the following form.
First we set up the variables:
\begin{verbatim}
gap> F := FreeGroup("b1","b2","b3","b4","b5","x1","x2","x3","y1","y2");;
gap> b1 := F.1;; b2 := F.2;; b3 := F.3;; b4 := F.4;; b5 := F.5;;
gap> x1 := F.6;; x2 := F.7;; x3 := F.8;; y1 := F.9;; y2 := F.10;;
\end{verbatim}
Then we input the data:
\begin{verbatim}
gap> ObA := [1,2];;
gap> ArrA := [ [1,1], [2,2] ];;
gap> ObB := [1,2,3];;
gap> ArrB := [ [b1,1,2], [b2,2,3], [b3,3,1], [b4,1,1], [b5,1,3] ];;
gap> RelB := [ [b1*b2*b3,b4] ];;
gap> FObA := [1,2];;
gap> FArrA := [b1,b2*b3];;
gap> XObA := [ [x1,x2,x3], [y1,y2] ];;
gap> XArrA := [ [y1,y2,y1],[x1,x2] ];;
\end{verbatim}
To combine all this data in one record do:
\begin{verbatim}
gap> KAN := rec( ObA:=ObA, ArrA:=ArrA,  ObB:=ObB, ArrB:=ArrB, RelB:=RelB, 
                 FObA:=FObA, FArrA:=FArrA, XObA:=XObA, XArrA:=XArrA );;
\end{verbatim}
To calculate the initial rules do
\begin{verbatim}
gap> IR := InitialRules( KAN );
\end{verbatim}
The output will be
\begin{verbatim}
i= 1, XA= [ x1, x2, x3 ], Ax= x1, rule= [ x1*b1, y1 ]
i= 1, XA= [ x1, x2, x3 ], Ax= x2, rule= [ x2*b1, y2 ]
i= 1, XA= [ x1, x2, x3 ], Ax= x3, rule= [ x3*b1, y1 ]
i= 2, XA= [ y1, y2 ], Ax= y1, rule= [ y1*b2*b3, x1 ]
i= 2, XA= [ y1, y2 ], Ax= y2, rule= [ y2*b2*b3, x2 ]
[ [ b1*b2*b3, b4 ], [ x1*b1, y1 ], [ x2*b1, y2 ], [ x3*b1, y1 ], 
  [ y1*b2*b3, x1 ], [ y2*b2*b3, x2 ] ]
\end{verbatim}
This means that there are five initial $\ep$-rules from: 
$( \ x_1|Fa_1,  x_1.a_1|\id_{FA_2} \ ), \ 
 ( \ x_2|Fa_1,  x_2.a_1|\id_{FA_2} \ ),$\\ 
$( \ x_3|Fa_1,  x_3.a_1|\id_{FA_2} \ ), \ 
 ( \ y_1|Fa_2,  y_1.a_1|\id_{FA_1} \ ), \ 
 ( \ y_2|Fa_2,  y_2.|a_11_{FA_1} \ ),$ \
i.e. $ \ x_1|b_1 \to y_1|\id_{B_2}, \ x_2|b_1 \to y_2|\id_{B_2}, \  
x_3|b_1 \to y_1|\id_{B_2}, \ y_1|b_2b_3  \to  x_1|\id_{B_1}, \
y_2|b_2b_3
\to
x_2|\id_{B_1} \ $ 
and one initial $K$-rule: $b_1b_2b_3 \to  b_4$. 
To attempt to complete the Kan extension presentation do:
\begin{verbatim}
gap> KB( IR );
\end{verbatim}
The output is:
\begin{verbatim}
[ [ x1*b1, y1 ], [ x1*b4, x1 ], [ x2*b1, y2 ], [ x2*b4, x2 ], [ x3*b1, y1 ], 
  [ x3*b4, x1 ], [ b1*b2*b3, b4 ], [ y1*b2*b3, x1 ], [ y2*b2*b3, x2 ] ]
\end{verbatim}
In other words to complete the system we have to add the rules 
$$x_1|b_4  \to  x_1, \quad x_2|b_4  \to  x_2, \text{ and } x_3|b_4  \to
x_1.$$ 
The result of attempting to compute the sets by doing:
\begin{verbatim}
gap> Kan(KAN);
\end{verbatim}
is a long list and then:
\begin{verbatim}
enumeration limit exceeded: complete rewrite system is:
[ [ x1*b1, y1 ], [ x1*b4, x1 ], [ x2*b1, y2 ], [ x2*b4, x2 ], [ x3*b1, y1 ], 
  [ x3*b4, x1 ], [ b1*b2*b3, b4 ], [ y1*b2*b3, x1 ], [ y2*b2*b3, x2 ] ]
\end{verbatim}
This means that the sets $KB$ for $B$ in $\bB$ are too large (the limit
set in 
the program is  1000). In fact this example is infinite. The complete
rewrite 
system is output instead of the sets. We can in fact use this to obtain
regular 
expressions for the sets. In this case the regular expressions are:
\begin{center}
\begin{tabular}{lcl}
$KB_1$ & $:=$ & $(x_1+x_2+x_3)|(b_5(b_3{b_4}^*b_5)^*b_3{b_4}^*+\id_{B_1}).$\\
$KB_2$ & $:=$ & $(x_1+x_2+x_3)|b_5(b_3{b_4}^*b_5)^*b_3{b_4}^*(b_1) +
(y_1+y_2)|\id_{B_2}.$\\
$KB_3$ & $:=$ & $(x_1+x_2+x_3)|b_5(b_3{b_4}^*b_5)^*(b_3{b_4}^*b_1b_2
+\id_{B_3}) + (y_1 + y_2)|b_2.$\\
\end{tabular}
\end{center}
The actions of the arrows are defined by concatenation followed by
reduction.
For example $x_1|b_5b_3b_4b_4b_5$ is an element of $KB_3$, so $b_3$ acts
on it 
to give $x_1|b_5b_3b_4b_4b_5b_3$ which is irreducible, and an element of
$KB_1$.\\

Details of how, in general, to obtain regular expressions will be given
in Chapter Four.

\section{Special Cases of the Kan Rewriting Procedure} 

\subsection{Groups and Monoids}

ORIGINAL PROBLEM:
Given a monoid presentation $mon\lan \Sigma|Rel\ran$, find a set of
normal forms 
for the monoid presented.\\
KAN INPUT DATA:
Let $\Gamma$ be the graph with one object and no arrows. Let $X\bullet$
be a one 
point set. Let $\bB$ be generated by the graph $\Delta$ with one object
and 
arrows labelled by $\Sigma$, it has relations $Rel\bB$ given by the
monoid 
relations. The functor $F$ maps the object of $\Gamma$ to the object of 
$\Delta$.\\
KAN EXTENSION: 
The Kan extension presented by $kan\lan \Gamma|\Delta|RelB|X|F\ran$ is
such that
$K\bullet$ is 
a set of normal forms for the elements of the
monoid, the arrows of $\bB$ (elements of $PX$) act on the right of $\bB$ 
by right multiplication. 
The natural transformation $\ep$ makes sure that the identity of $\bB$
acts trivially and helps to define the normal form function.
The normal form function is 
$w \mapsto \ep_\bullet(1)\cdot(w):=Kw(\ep_\bullet(1))$.\\

In this case the method of completion
is the standard Knuth-Bendix procedure used
for many years for working with monoid presentations of groups and
monoids.
This type of calculation is well documented.

\subsection{Groupoids and Categories}

ORIGINAL PROBLEM: 
To specify a set of normal forms for the elements of a 
groupoid or category given by a finite category presentation 
$cat\lan \Lambda|Rel \ran$.\\
KAN INPUT DATA: 
Let $\Gamma$ be the discrete graph with no arrows and object set equal
to 
$\ob\Lambda$. Let $XA$ be a distinct one object set for each
$A\in\ob\Gamma$. 
Let $\bB$ be the category generated by $\Delta:=\Lambda$ with relations 
$Rel\bB:=Rel$. Let $F$ be defined by the identity map on the objects.\\
KAN EXTENSION: 
Then the Kan extension presented by $kan\lan \Gamma | \Delta | RelB | X
| F \ran$
is such that $KB$ is a 
set of normal forms for the arrows of the category with target $B$, the
arrows 
of $\bB$ (elements of $P\Gamma$) act on the right of $\bB$ 
by right multiplication. 
The natural transformation $\ep$ makes sure that the identities of $\bB$
act
trivially and helps to define the normal form function.
The normal form function is 
$w \mapsto \ep_A\cdot(w):=Kw(\ep_A)$.\\

\begin{example}
\emph{Consider the group $S_3$ presented by $\lan x,y|x^3,y^2,xyxy\ran.$
The elements are \\
$\{ \id, x, y, x^2, xy, yx \}$.
The covering groupoid is generated by the Cayley graph.
The 12 generating arrows of the groupoid are $G \times X$: 
$$
\{ [\id,x],[x,x],[y,x],\ldots,[yx,x],[\id,y],[x,y],\ldots,[yx,y] \}.
$$ 
To make calculations clearer, we relabel them 
$\{ a_1, a_2, a_3, \ldots, a_6, b_1, b_2, \ldots,  b_6 \}$.}\\

\emph{The groupoid has 18 relators $G \times R$
-- the boundaries of irreducible cycles of the graph. 
The cycles may be written $[\id,x^3]$ and the corresponding boundary is
$[\id,x][x,x][x^2,x] ~ $ i.e. $a_1a_2a_4$.
For the category presentation of the group we could add in the inverses
$\{ A_1, A_2, \ldots, A_6, B_1, B_2, \ldots, B_6 \}$ with the relators
$A_1a_1$
and $a_1A_1$ etc and end up with a category presentation with 24
generators and the 42 relations. In this case however the groupoid is
finite and so there is no need to do this. For example there would be no 
need for $A_2$ because $(a_2)^{-1}=a_4a_1$.}

{\Large
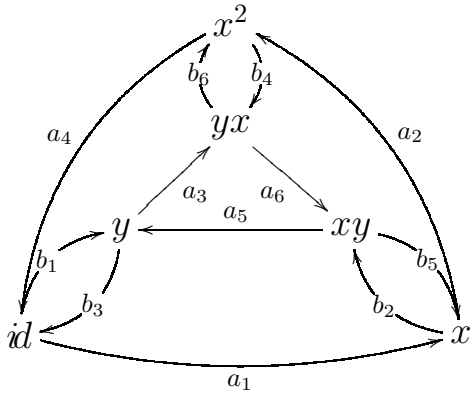
$$\xymatrix{&& x^2 \ar@/_1.7pc/[dddll]_{a_4} \ar@/^1pc/[d]|{b_4} && \\
            && yx \ar@/^1pc/[u]|{b_6} \ar[dr]_{a_6} && \\
            & y \ar@/^1pc/[dl]|{b_3} \ar[ur]_{a_3} &
            & xy \ar@/^1pc/[dr]|{b_5} \ar[ll]_{a_5} & \\
            \id  \ar@/^1pc/[ur]|{b_1} \ar@/_1pc/[rrrr]_{a_1} &&&&
            x \ar@/_1.7pc/[uuull]_{a_2} \ar@/^1pc/[ul]|{b_2} \\}
$$}

\emph{Now suppose the left hand sides of two rules overlap 
(for example $(a_1a_2a_4, \id)$ and $(a_4b_1a_3b_6, \id)$) 
in one of the two possible ways previously described 
then we have a critical pair
$(b_1a_3b_6,a_1a_2)$ ).
The following is $\mathsf{GAP}$ output of the completion of the rewrite system
for the
covering groupoid of our example:}

\begin{verbatim}
gap> Rel;                         ## Input rewriting system:
[ [ a1*a2*a4, IdWord ], [ a2*a4*a1, IdWord ], [ a4*a1*a2, IdWord ],
  [ a3*a6*a5, IdWord ], [ a6*a5*a3, IdWord ], [ a5*a3*a6, IdWord ],
  [ b1*b3, IdWord ], [ b3*b1, IdWord ], [ b2*b5, IdWord ], 
  [ b5*b2, IdWord ], [ b4*b6, IdWord ], [ b6*b4, IdWord ], 
  [ a1*b2*a5*b3, IdWord ], [ a2*b4*a6*b5, IdWord ], 
  [ a3*b6*a4*b1, IdWord ], [ a4*b1*a3*b6, IdWord ],
  [ a5*b3*a1*b2, IdWord ], [ a6*b5*a2*b4, IdWord ] ]
gap> KB( Rel );                   ## Completed rewriting system:
[ [ b1*b3, IdWord ], [ b2*b5, IdWord ], [ b3*b1, IdWord ], 
  [ b4*b6, IdWord ], [ b5*b2, IdWord ], [ b6*b4, IdWord ], 
  [ a1*a2*a4, IdWord ], [ a1*a2*b4, b1*a3 ], [ a1*b2*a5, b1 ], 
  [ a2*a4*a1, IdWord ], [ a2*a4*b1, b2*a5 ], [ a2*b4*a6, b2 ], 
  [ a3*a6*a5, IdWord ], [ a3*a6*b5, b3*a1 ], [ a3*b6*a4, b3 ], 
  [ a4*a1*a2, IdWord ], [ a4*a1*b2, b4*a6 ], [ a4*b1*a3, b4 ], 
  [ a5*a3*a6, IdWord ], [ a5*a3*b6, b5*a2 ], [ a5*b3*a1, b5 ], 
  [ a6*a5*a3, IdWord ], [ a6*a5*b3, b6*a4 ], [ a6*b5*a2, b6 ], 
  [ b1*a3*a6, a1*b2 ],  [ b1*a3*b6, a1*a2 ], [ b2*a5*a3, a2*b4 ], 
  [ b2*a5*b3, a2*a4 ],  [ b3*a1*a2, a3*b6 ], [ b3*a1*b2, a3*a6 ], 
  [ b4*a6*a5, a4*b1 ],  [ b4*a6*b5, a4*a1 ], [ b5*a2*a4, a5*b3 ], 
  [ b5*a2*b4, a5*a3 ],  [ b6*a4*a1, a6*b5 ], [ b6*a4*b1, a6*a5 ] ]
\end{verbatim}

\emph{
It is possible from this to enumerate elements of the category.
One method is to start with all the shortest arrows
($a_1,a_2,\ldots,b_6$) and see 
which ones reduce and build inductively on the irreducible ones:\\
Firstly we have the six identity arrows
$\id_{\id}, \ \id_x, \ \id_y, \ \id_{x^2}, \ \id_{xy}, \ \id_{yx}$.\\
Then the generators $a_1, \ a_2, \ a_3, \ a_4, \ a_5, \ a_6, \ b_1, \
b_2, \
b_3, \ b_4, \ b_5, \ b_6$ are all irreducible.\\
Now consider paths of length 2:\\
$a_1a_2, \ a_1b_2, \ a_2a_4, \ a_2b_4, \
 a_3a_6, \ a_3b_6, \ a_4a_1, \ a_4b_1, \
 a_5a_3, \ a_5b_3, \ a_6a_5, \ a_6b_5, \
 b_1a_3, \ b_1b_3 \to \id_{\id},$\\
$b_2a_5, \ b_2b_5 \to \id_x, \ b_3a_1, \ b_3b_1 \to \id_y, \
 b_4a_6, \ b_4b_6 \to \id_{x^2}, \
 b_5a_2, \ b_5b_2 \to \id_{xy}, \ b_6a_4, \ b_6b_4 \to \id_{yx}$.\\
Building on the irreducible paths we get the paths of length 3:
$a_1a_2a_4 \to \id_{\id}, \ a_1a_2b_4 \to b_1a_3,$\\
$a_1b_2a_5 \to b_1, \ a_1b_2b_5 \to a_1, \ a_2a_4a_1 \to \id_x,\ldots$\\
All of them are reducible, and so we can't build any longer paths; the 
covering groupoid has 30 morphisms and 6 identity arrows and is the tree 
groupoid with six objects.
}\end{example}

\begin{example}\emph{
This is a basic example to show how it is possible to specify the arrows
in an infinite small category with a finite complete presentation.
Let $\bC$ be the category generated by the following graph $\Gamma$
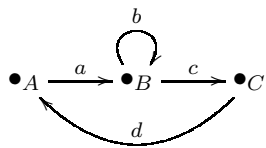
$$
\xymatrix{ \bullet_A \ar[r]^a & \bullet_B \ar@(ul,ur)^b \ar[r]^c 
         & \bullet_C \ar@/^2pc/[ll]_d}
$$
with the relations $b^2c=c,\ ab^2=a$. This rewriting system is complete,
and so we can
determine whether two arrows in the free category $P\Gamma$ are
equivalent in 
$\bC$. An automaton can be drawn (see chapter 3), and from this we can
specify
the language which is the set of normal forms. It is in fact 
$$
a(cd(acd)*ab+bcd(acd)*ab) + b^{\dagger} + cd(acd)^*ab + d(acd)^*ab 
$$
(and the three identity arrows) where $(acd)^*$ is
used to denote the set of elements of $\{ acd \}^*$ (similarly
$b^\dagger$), 
so $d(acd)*$, for example, denotes the set $\{d, dacd, dacdacd,
dacdacdacd,
\ldots \}$, 
$+$ denotes the union and $-$ the difference of sets.
This is the standard notation of languages and regular expressions. 
}\end{example}

\subsection{Coset systems and Congruences}

ORIGINAL PROBLEM:
Given a finitely presented group $G$ and a finitely generated subgroup
$H$
find a set of normal forms for the coset representatives
of $G$ with respect to $H$.\\
KAN INPUT DATA:  
Let $\Gamma$ be the one object graph $\Gamma$ with arrows labelled by 
the subgroup generators. Let $X\bullet$ be a one point set on which the
arrows 
of $\Gamma$ act trivially. Let $\bB$ be the category generated by the
one object 
graph $\Delta$ with arrows labelled by the group generators, with the
relations
$Rel\bB$ of $\bB$ being the group relations. Let $F$ be defined on
$\Gamma$ by 
inclusion of the subgroup elements to the group.\\
KAN EXTENSION:
The Kan extension presented by $kan\lan \Gamma|\Delta|RelB|X|F\ran$ is
such that 
the
set 
$K\bullet$ is a set of representatives for the cosets, $Kb$ defines the
action 
of the group on the cosets $Hg\mapsto Hgb$ and $\ep_\bullet$ maps the
single 
element of $X\bullet$ to the representative for $H$ in $K\bullet$.
Therefore it follows that the Kan extension defined is computable if and
only if 
the coset system is computable. \\ 

In the monoidal case $F$ is the inclusion of the submonoid $\bA$ of the 
monoid $\bB$, and the action is trivial as before. The Kan extension of
this 
action gives the quotient of $\bB$ by the right congruence generated by
$\bA$, 
namely the equivalence relation generated by $ab\sim 
b$ for all $a \in \bA, b \in \bB$, with the induced right action of
$\bB$.\\

It is appropriate to give a calculated example here. The example is
infinite so 
standard Todd-Coxeter methods will not terminate, but the Kan extension
/ 
rewriting procedures enable the complete specification of the coset
system.

\begin{example}\emph{ 
Let $\bB$ represent the infinite group presented by  
$$grp \lan a, b, c \ | \ a^2b=ba, a^2c=ca, c^3b=abc, caca=b \ran$$ 
and let $\bA$ represent the subgroup generated by $\{ c^2 \}$. \\
We obtain one initial $\ep$-rule (because $\bA$ has one generating
arrow) 
i.e. $ H|c^2 \to H|\id. $\\ 
We also have four initial $K$-rules corresponding to the relations of
$\bB$:}
$$a^2b \to ba, \ a^2c \to ca, \ c^3b \to abc, \ caca \to b.$$
\emph{
Note: On completion of this rewriting system for the group, we find 24
rules 
and for all $n \in \bN$ both $a^n$ and $c^n$ are irreducibles with
respect to 
this  system (one way to prove that the group is infinite).}\\

\emph{
The five rules are combined and an infinite complete system for the Kan  
extension of the action is easily found (using Knuth-Bendix with the 
length-lex order).  
The following is the $\mathsf{GAP}$ output of the set of 32 rules:}

\begin{verbatim} 
[ [ H*b, H*a ], [ H*a^2, H*a ], [ H*a*b, H*a ], [ H*c*a, H*a*c ], 
  [ H*c*b, H*a*c ], [ H*c^2, H ], [ a^2*b, b*a ], [ a^2*c, c*a ], 
  [ a*b^2, b^2 ], [ a*b*c, c*b ], [ a*c*b, c*b ], [ b*a^2, b*a ], 
  [ b*a*b, b^2 ], [ b*a*c, c*b ], [ b^2*a, b^2 ], [ b*c*a, c*b ], 
  [ b*c*b, b^2*c ], [ c*a*b, c*b ], [ c*b*a, c*b ], [ c*b^2, b^2*c ], 
  [ c*b*c, b^2 ], [ c^2*b, b^2 ], [ H*a*c*a, H*a*c ], [ H*a*c^2, H*a ], 
  [ b^4, b^2 ], [ b^3*c, c*b ], [ b^2*c^2, b^3 ], [ b*c^2*a, b^2 ], 
  [ c*a*c*a, b ], [ c^2*a^2, b*a ], [ c^3*a, c*b ], [ c*a*c^2*a, c*b ] ] 
\end{verbatim} 

\emph{
Note that the rules without $H$ i.e. the two-sided rules, constitute a  
complete rewriting system for the group. 
The set $KB$ (recall that there is only one object $B$ of $\bB$) 
is infinite. It is the set of (right) cosets of the subgroup in the
group. 
Examples of these cosets include:  
$$H, Ha, Hc, Ha^2, Hac, Ha^3, Ha^4, Ha^5,\ldots$$ 
A regular expression for the coset representatives is: 
$$a^*+c+ac.$$
Alternatively consider the subgroup generated by $b$. Add the rule 
$Hb \to H$ and the complete system below is obtained:}

\begin{verbatim}
[ [ H*a, H ], [ H*b, H ], [ H*c*a, H*c ], [ H*c*b, H*c ], [ H*c^2, H ], 
  [ a^2*b, b*a ], [ a^2*c, c*a ], [ a*b^2, b^2 ], [ a*b*c, c*b ], 
  [ a*c*b, c*b ], [ b*a^2, b*a ], [ b*a*b, b^2 ], [ b*a*c, c*b ], 
  [ b^2*a, b^2 ], [ b*c*a, c*b ], [ b*c*b, b^2*c ], [ c*a*b, c*b ], 
  [ c*b*a, c*b ], [ c*b^2, b^2*c ], [ c*b*c, b^2 ], [ c^2*b, b^2 ], 
  [ b^4, b^2 ], [ b^3*c, c*b ], [ b^2*c^2, b^3 ], [ b*c^2*a, b^2 ], 
  [ c*a*c*a, b ], [ c^2*a^2, b*a ], [ c^3*a, c*b ], [ c*a*c^2*a, c*b ] ] 
\end{verbatim} 
\emph{
Again, the two-sided rules are the rewriting system for the group. 
This time the subgroup has index 2, and the coset representatives are 
$\id$ and $c$.
}\end{example}

\subsection{Equivalence Relations and Equivariant Equivalence Relations} 
ORIGINAL PROBLEM:
Given a set $\Omega$ and a relation $Rel$ on $\Omega$. 
Find a set of representatives for the equivalence classes of the set
$\Omega$ 
under the equivalence relation generated by $Rel$.\\
KAN INPUT DATA:
Let $\Gamma$ be the graph with object set $\Omega$ and
generating arrows $a:A_1 \to A_2$ if $(A_1,A_2) \in Rel$.
Let $XA:=\{A\}$ for all $A\in\Omega$. The arrows of $\Gamma$ act
according to 
the relation, so $src(a)\cdot a=tgt(a)$.
Let $\Delta$ be the graph with one object and no arrows so that $\bB$ is
the 
trivial category with no relations. Let $F$ be the null functor.\\
KAN EXTENSION:
The Kan extension presented by $kan\lan \Gamma|\Delta|RelB|X|F\ran$
is such that $K\bullet:= \Omega/{\stackrel{*}{\lra}}_{Rel}$ is a set of 
representatives for
the equivalence classes of the set $\Omega$ under the equivalence 
relation generated by $Rel$.\\

Alternatively let $\Omega$ be a set with a group or monoid $M$ acting on
it.
Let $Rel$ be a relation on $\Omega$.
Define $\Gamma$ to have object set $\Omega$ and
generating arrows $a:A_1 \to A_2$ if $(A_1,A_2) \in Rel$ or if $A_1\cdot
m=A_2$
Again, $XA:=\{A\}$ for $A\in\ob\Gamma$ and the arrows act as in the case
above. 
Let $\Delta$ be the one object graph with arrows labelled by generators
of $M$ 
and for $\bB$ let $Rel\bB$ be the set of monoid relations. Let $F$ be
the null 
functor.
The Kan extension gives the action of $M$ on the 
quotient of $X$ by the $M$-equivariant equivalence relation generated by
$Rel$. 
This example illustrates the advantage of working in categories, since
this is a  
coproduct of categories which is a fairly simple construction.

\subsection{Orbits of Actions}

ORIGINAL PROBLEM:
Given a group $G$ which acts on a set $\Omega$,  find a set $KB$ of 
representatives for the orbits of the action of $\bA$ on $\Omega$.\\
KAN INPUT DATA:
Let $\Gamma$ be the one object graph with arrows labelled by the
generators of 
the group. Let $X\bullet:=\Omega$. Let $\Delta$ be the one object, zero
arrow 
graph generating the trivial category $\bB$ with $Rel\bB$ empty. Let $F$
be the 
null functor.\\
KAN EXTENSION:
The Kan extension presented by $kan\lan \Gamma|\Delta|RelB|X|F\ran$ is
such that
$K\bullet$ is a set of representatives for the orbits of the 
action of the group on $\Omega$.

We present a short example to demonstrate the procedure in this case.

\begin{example}\emph{
Let $\bA$ be the symmetric group on three letters with presentation\\
$mon\lan a,b| a^3,b^2,abab\ran$ and let $X$ be the set $\{v,w,x,y,z\}$.
Let
$\bA$ act on $X$ by giving $a$ the effect of the permutation 
$(v \ w \ x)$ and $b$ the effect of $(v \ w)(y \ z)$.}\\ 

\emph{
In this calculation we have a number of $\ep$-rules and no $K$-rules. 
The $\ep$-rules just list the action, namely (trivial actions omitted):
$$
v \to w, ~ w \to x, ~ x \to v, ~ v \to w, ~ w \to v, ~ y \to z, ~ z \to y. 
$$
The system of rules is complete and reduces to 
$\{ w \to v, ~ x \to v,~  z \to y\}$.
Enumeration is simple: $v, \ w \to v, \ x \to v, \ y, \ z \to y$, so
there are
two orbits of $\Omega$ represented by $v$ and $y$.
}\end{example}

This is a small example. With large 
examples the idea of having a minimal element (normal form) in each
orbit
to act as an anchor or point of comparison makes a lot of sense. 
This situation serves as another illustration of rewriting in the
framework of
a Kan extension, showing not only that rewriting gives a result, but
that it is
the procedure one uses naturally to do the calculation.\\

One variation of this is if $\Omega$ is the set of elements of the group
and the 
action is conjugation: $x^a := a^{-1}xa$. Then the orbits are the  
\emph{conjugacy classes} of the group.

\begin{example}\emph{ Consider
the quarternion group, presented by $\lan a,b \ | \ a^4, b^4,
abab^{-1}, a^2b^2\ran ~$ and  
$\Omega=\{ \id, ~ a, ~ b, ~ a^2, ~ ab, ~ ba, ~ a^3, ~ a^2b \}$
-- enumerating the elements of the group using the method described 
in Example 3. 
Construct the Kan extension as above, where
the actions of $a$ and $b$ are by conjugation on elements of $\bA$.\\
There are 16 $\ep$-rules which reduce to 
$\{ a^3 \to a, ~ a^2b \to b, ~ ba \to ab \}$. 
The conjugacy classes are enumerated by applying these
rules to the elements of $\bA$. The irreducibles are 
$\{\id, ~ a, ~ b, ~ a^2, ~ ab\}$, and these are representatives of the five conjugacy
classes. 
}\end{example}

\subsection{Colimits of Diagrams of $\sets$}

ORIGINAL PROBLEM:
Given a presentation of a category action $act\lan \Gamma|X\ran$ find
the 
colimit of the diagram in $\sets$ on which the category action is
defined.\\
KAN INPUT DATA:
Let $\Gamma$ and $X$ be those given by the action presentation. Let
$\Delta$ be 
the graph with one object and no arrows that generates the trivial
category 
$\bB$ with $Rel\bB$ empty. Let $F$ be the null functor.\\
KAN EXTENSION:
The Kan extension presented by $kan\lan \Gamma|\Delta|RelB|X|F\ran$ is
such that
$K\bullet$ is 
the colimit object, and $\ep$ is the set of colimit functions of the
functor 
$X:\bA \to \sets$.\\
 
Particular examples of 
this are when $\bA$ has two objects $A_1$ and $A_2$, and two 
non-identity arrows $a_1$ and $a_2$ from $A_1$ to $A_2$, \ and $Xa_1$ 
and $Xa_2$ are functions from the set $XA_1$ to the set $XA_2$ 
(\emph{coequaliser} of $a_1$ and $a_2$ in $\sets$); $\bA$ has three 
objects $A_1$, $A_2$ and $A_3$ and two non-identity arrows $a_1:A_1 \to 
A_2$ and $a_2:A_1 \to A_3$. \ $XA_1$, $XA_2$ and $XA_2$ are sets, and 
$Xa_1$ and $Xa_2$ are functions between these sets (\emph{pushout} of 
$a_1$ and $a_2$ in $\sets$). The following example is included not as an
illustration of rewriting but to show another situation where
presentations of Kan extensions can be used to express a problem
naturally.

\begin{example}\emph{
Suppose we have two sets $\{ x_1, x_2, x_3 \}$ and $\{y_1, y_2,
y_3, y_4 \}$,
with two functions from the first to the second given by 
$( x_1 \mapsto y_1, ~ x_2 \mapsto y_2, ~ x_3 \mapsto y_3 )$ and
$( x_1 \mapsto y_1, ~ x_2 \mapsto y_1, ~ x_3 \mapsto y_3 )$.\\
Then we can calculate the coequaliser. 
We have a number of $\ep$-rules
$$ 
y_1|\id_\bullet \to x_1|\id_\bullet, ~ 
y_2|\id_\bullet \to x_2|\id_\bullet, ~
y_3|\id_\bullet \to x_3|\id_\bullet, ~
y_1|\id_\bullet \to x_1|\id_\bullet, ~
y_2|\id_\bullet \to x_1|\id_\bullet, ~
y_3|\id_\bullet \to x_3|\id_\bullet. ~
$$  
There is just one overlap, between 
$(y_2|\id_\bullet \to x_1|\id_\bullet)$ and 
$(y_2|\id_\bullet \to x_2|\id_\bullet)$: to 
resolve the critical pair we add the rule 
$(x_2|\id_\bullet \to x_1|\id_\bullet)$, 
and the system is complete: 
$$\{y_1|\id_\bullet \to x_1|\id_\bullet, ~
    y_2|\id_\bullet \to x_1|\id_\bullet, ~
    y_3|\id_\bullet \to x_3|\id_\bullet, ~
    x_2|\id_\bullet \to x_1|\id_\bullet
\}.$$ 
The elements of the set $K\bullet$ are easily enumerated:\\ 
$$
x_1|\id_\bullet, \ x_2|\id_\bullet \to x_1|\id_\bullet, \
x_3|\id_\bullet, \ y_1|\id_\bullet \to x_1|\id_\bullet, \ 
y_2|\id_\bullet \to x_1|\id_\bullet, \ 
y_3|\id_\bullet \to x_3|\id_\bullet, \ 
y_4|\id_\bullet.
$$ 
So the coequalising set is 
$$K\bullet=\{x_1|\id_\bullet, x_3|\id_\bullet, y_4|\id_\bullet \},$$ 
and the coequaliser function 
to it from $XA_2$ is given by $y_i\mapsto y_i|\id_\bullet$ for
$i=1,\ldots,4$ 
followed by reduction defined by $\to$ to an element of $K\bullet$. 
}\end{example}

\subsection{Induced Permutation Representations}

Let $\bA$ and $\bB$ be groups and let $F:\bA \to \bB $ be a morphism of
groups. 
Let $\bA$ act on the set $XA$. The Kan extension of this action along 
$F$ is known as the action of $\bB$ {\em induced} from that of $\bA$ 
by $F$, and is written $F_*(XA)$. It can be constructed simply as the 
set $X \times \bB$ factored by the equivalence relation generated by 
$(xa,b)\sim (x,F(a)b)$ for all $x \in XA, a \in \bA,b \in \bB$. The 
natural transformation $\ep$ is given by $ x \mapsto [x,1]$, where 
$[x,b]$ denotes the equivalence class of $(x,b)$ under the equivalence 
relation $\sim$. The morphism $F$ can be factored as an epimorphism 
followed by a monomorphism, and there are other descriptions of 
$F_*(XA)$ in these cases, as follows.\\  
 
Suppose first that $F$ is an epimorphism with kernel $N$. Then we can 
take as a representative of $F_*(XA)$ the orbit set $X/N$ with the 
induced action of $\bB$.\\  
 
Suppose next that $F$ is a monomorphism, which we suppose is an 
inclusion. Choose a set $T$ of representatives of the right cosets  
of $\bA$ in $\bB$, so that $1 \in T$. Then the induced representation 
can be taken to be $XA \times T$ with $\ep$ given by $x \mapsto (x,1)$ 
and the action given by $(x,t)^b= (xa,u)$ where $t,u \in T, b \in \bB, a 
\in \bA$ and $tb=au$.\\  
 
On the other hand, in practical cases, this factorisation of $F$ may not 
be a convenient way of determining the induced representation.  
In the case $\bA,\bB$ are monoids, so that $XA$ is a transformation 
representation of $\bA$ on the set $XA$, we have in general no 
convenient description of the induced transformation representation 
except by one form or another of the construction of the Kan 
extension.

\include{three}
\chapter{Reduction and Machines} 
 
In the first section automata are considered in the standard way, as acceptors, but applied to 
the Kan extensions of Chapter 2. We show how to construct automata which accept 
the unique normal forms of the elements of each set $KB$ for $B\in\ob\Delta$. 
Creating accepting automata for such structures is new, and we describe their 
construction from the complete rewriting systems as well as showing how to apply 
standard automata theory \cite{Hopcroft} to obtain a regular expression for the 
language which is the set of irreducible elements. 
Further, we extend the ideas to algebras. It appears that some work is being 
done in this line \cite{Patrik} (monomial acceptors) but it is still appropriate 
to include it here, to relate the concepts.\\

In the second section we move on to consider a more useful class of automata 
-- those 
with output. These machines not restricted to accepting or rejecting strings, 
but can reduce them into the unique irreducible representative forms. 
The best known example of this is the use of the Cayley Graph to work out 
multiplication of group elements. The use of the Cayley Graph 
as a reduction machine is the first thing to be described. 
Rewriting systems for Kan extensions can be translated into 
reduction machines for Kan extensions. These machines are defined as Moore 
machines. 
The next consideration is of reduction machines 
for algebras, which are constructed from the Gr\"obner bases. I believe this to 
be a new idea. The construction and operation of the ``Gr\"obner machines'' is 
explained, using a small Hecke Algebra as an example.\\ 
 
The final section introduces a third type of machine: a Petri net. 
There are many different classes of Petri nets, and we show how to consider the 
``Gr\"obner machine'' of the previous section as a Petri net. We also show 
how commutative Gr\"obner bases may be applied to successfully solve the 
standard problems posed for reversible Petri nets. 
This small section speculates on the relation between Petri nets and Gr\"obner 
bases and does not prove any results. It is hoped that it provides a starting 
point for further investigations into the relation between Petri nets and 
Gr\"obner bases.

\section{Normal Forms Acceptors} 
 
\subsection{Definitions and Notation} 
 
For a detailed introduction to automata theory refer to \cite{Cohen} or 
\cite{Hopcroft}. This section only outlines the essential ideas we use.\\

A (finite) \textbf{deterministic automaton} is a 5-tuple $\underline{A}=(S, 
\Sigma, s_0,  
\delta, Q)$ \ where $S$ is a finite set of \emph{states} (represented by circles), 
$s_0 \in S$ is the \emph{initial state} (marked with an arrow), 
$\Sigma$ is a finite \emph{alphabet},  
$\delta:S \times \Sigma \to S$ is the \emph{transition}, 
$Q\subseteq S$ is the set of \textbf{terminal states} (represented by double 
circles).
A deterministic automaton $\underline{A}$ is \textbf{complete} 
if $\delta$ is a function, and \textbf{incomplete} if it is only a 
partial function. If $\underline{A}$ is incomplete, then when 
$\delta(s,a)$ is undefined, the automaton is said to \textbf{crash}.\\ 
 
The \textbf{extended state transition} $\delta^*$ is the extension of $\delta$ 
to $\Sigma^*$. It is defined by 
$\delta^*(s,\id):=s$, $\delta^*(s,a):=\delta(s,a)$, $\delta^*(s,aw) := 
\delta^*(\delta(s,a),w)$ where $s \in S$, $a\in \Sigma$ and $w$ is a string
in $\Sigma^*$.
We are interested in the final state $\dd^*(s_0,w)$ of the machine after a 
string $w$ has been completely read. If the machine crashes 
or ends up at a non-terminal state then the string is said to have been 
\textbf{rejected}. If it ends up at a terminal state then we say the string 
is \textbf{accepted}. \\

A \textbf{language} over a given alphabet $\Sigma$ is a subset $L$ of
$\Sigma^*$. 
The set $L(\underline{A})$ of all acceptable strings is the 
 \textbf{language accepted by the automaton} $\underline{A}$. A language $L$ is a 
\textbf{recognisable} if it is accepted by some automaton $\underline{A}$. 
Two automata are \textbf{equivalent} if their languages are equal. 
The \textbf{complement} of a complete, deterministic automaton is found by 
making non-terminal states terminal and vice versa. If the language accepted by 
an automaton is $L$, then the language accepted by its complement is 
$\Sigma^*-L$.

\begin{lem}[\cite{Cohen}] 
Let $\underline{A}=(S, \Sigma, s_0, \delta, Q)$ be an incomplete deterministic 
automaton. Then there exists a complete deterministic automaton 
$\underline{A}^{CP}$ such that $L(\underline{A})=L(\underline{A}^{CP})$. 
\end{lem} 
\begin{out} 
Define $\underline{A}^{CP}=(S \sqcup d, \Sigma, s_0, \delta_1, Q)$ 
where the transition $\delta_1:S \times \Sigma \to S$ is defined by 
$\delta_1(s,a):=\delta(s,a)$ if $\delta(s,a)$ is defined, otherwise 
$\delta_1(s,a):=d$, and $\delta_1(d,a):=d$. 
\end{out}

Diagrammatically this means that automata may be completed by adding one further 
non-terminal (dump) state $d$ and adding in all the missing arrows so that they 
point to this state. \\

A \textbf{non-deterministic automaton} is a 5-tuple 
$\underline{A}=(S, \Sigma, S_0, \delta, Q)$ \  
where $S$ is a finite set of states, 
$S_0 \subseteq S$ is a set of initial states, 
$\Sigma$ is a finite alphabet,  $Q\subseteq S$ is the set of 
terminal states and 
$\delta:S \times \Sigma \to \mathbb{P}(S)$ is the transition mapping 
where $\mathbb{P}(S)$ is the power set.

\begin{lem}[\cite{Cohen}] 
Let $\underline{A}=(S, \Sigma, S_0, \delta_1, Q)$ be a non-deterministic 
automaton. Then there exists a deterministic automaton $\underline{A}^d$ such 
that $L(\underline{A}^d)=L(\underline{A})$. 
\end{lem} 
\begin{out} 
Define $\underline{A}^d:=(S^d,\Sigma,{S_0}^d,\delta^d,Q^d)$ where 
$S^d:=\mathbb{P}(S)$ then ${S_0}^d = S_0 \in S^d$, $Q^d:=\{U\in \mathbb{P}(S)|U\cap 
Q\not= \emptyset\}$. Define $\delta^d(U,a):=\bigcup_{u\in U}\delta(u,a)$ for 
$a\in \Sigma$. It can be verified that $L(\underline{A}^d)=L(\underline{A})$. 
\end{out} 
 
In practice a non-deterministic automaton may be made deterministic by drawing 
a \emph{transition tree} and then converting the tree into an automaton;  
for details of this see  \cite{Cohen}.\\

Let $\Sigma$ be a set (alphabet). The following notation is standard when 
working with languages. 
The empty word will be denoted $\id$. If $x \in \Sigma^*$ then we will write $x$ for $\{x\}$. 
If $A,B \in \mathbb{P}\Sigma^*$ then 
$A+B:=A \cup B$, $A-B:=A \,/ \,B$. Therefore, for example $(x+y)^*+z=\{x,y\}^* \cup \{z\}$.\\

A \textbf{regular expression} over $\Sigma$ is a string of symbols formed by 
the rules 
\begin{enumerate}[i)] 
\item $a_1\cdots a_n$ is regular for $a_1,\ldots,a_n \in \Sigma$, 
\item $\emptyset$ is regular, 
\item $\id$ is regular, 
\item if $x$ and $y$ are regular then $xy$ is regular, 
\item if $x$ and $y$ are regular then $x+y$ is regular, 
\item if $x$ is regular then $x^*$ is regular. 
\end{enumerate} 
 
A \textbf{right linear language equation} over $\Sigma$ is an 
expression $X=AX+E$ where $A,X,E \subseteq \Sigma^*$.
 
\begin{thm}[Arden's Theorem \cite{Cohen}]\mbox{ } 
Let $A,X,E \subseteq \Sigma^*$ such that $X=AX+E$ where $A$ and $E$ are known and $X$ is unknown. 
Then 
\begin{enumerate}[i)] 
\item $A^*E$ is a solution, 
\item if $Y$ is any solution then $A^*E\in Y$, 
\item if $\id \not\in A$ then $A^*E$ is the unique solution. 
\end{enumerate} 
\end{thm} 
 
\begin{thm}[\cite{Cohen}] 
A system of right linear language equations:
\vspace{-0.2cm} 
\begin{center} 
\begin{tabular}{lllllllll} 
$X_0$     & $=$ & $A_{0,0}X_0$   & $+$ & $\cdots$ & $+$ & $A_{0,n-1}X_{n-1}$   & $+$ & $E_0$,\\ 
$X_1$     & $=$ & $A_{1,0}X_0$   & $+$ & $\cdots$ & $+$ & $A_{1,n-1}X_{n-1}$   & $+$ & $E_1$,\\ 
$\cdots$  &     & $\cdots\cdots$ &     & $\cdots$ &     & $\cdots\cdots\cdots$ &     & $\cdots$\\ 
$X_{n-1}$ & $=$ & $A_{n-1,0}X_0$ & $+$ & $\cdots$ & $+$ & $A_{n-1,n-1}X_{n-1}$ & $+$ & $E_{n-1}$. 
\end{tabular}
\end{center} 
where $A_{i,j},E_i \in \mathbb(\Sigma^*)$ and $\id \not\in A_{i,j}$ for 
$i,j=0,\ldots,n-1$, has a unique solution. 
\end{thm} 
\begin{out} 
Begin with the last equation. 
By assumption $\id \not\in A_{n-1,n-1}$. 
So by Arden's theorem 
$X_{n-1} =  A_{n-1,n-1}^* ( A_{n-1,0}X_0 + \cdots + 
A_{n-1,n-2}X_{n-2} + E_{n-1})$. 
Substitute this value for $X_{n-1}$ into the remaining $n-1$ equations and 
repeat the procedure. Eventually an equation in $X_0$ only will be obtained which
can be solved explicitly. The back-substitution will give explicit values of 
$X_1,\ldots,X_{n-1}$. 
\end{out} 
 
\begin{thm}[\cite{Cohen}] 
Let $\underline{A}$ be a (non)-deterministic automaton. 
Then $L(\underline{A})$ is regular. 
\end{thm} 
\begin{out} (for the deterministic case)\\ 
Let $\underline{A}:=(S,\Sigma,s_0,\delta,Q)$ where $S=\{s_0,\ldots,s_{n-1}\}$. 
Define $X_i := \{z\in\Sigma^*:\delta(s_i,z)\in Q\}$ for $i = 0, \ldots, n-1$. 
It is clear that $L(\underline{A})=X_0$. 
Define $E_i:=\emptyset$ if $s_i \not\in Q$ and $E_i:=\{\id\}$ if $s_i \in Q$ for 
$i=0,\ldots,n-1$. 
Define $A_{i,j}:=\{a\in\Sigma:\delta(s_i,a)=s_j \}$ for $i,j=0,\ldots,n-1$. 
Form the following system:
\vspace{-.2cm}
\begin{center} 
\begin{tabular}{lllllllll} 
$X_0$     & $=$ & $A_{0,0}X_0$   & $+$ & $\cdots$ & $+$ & $A_{0,n-1}X_{n-1}$   & $+$ & $E_0$,\\ 
$X_1$     & $=$ & $A_{1,0}X_0$   & $+$ & $\cdots$ & $+$ & $A_{1,n-1}X_{n-1}$   & $+$ & $E_1$,\\ 
$\cdots$  &     & $\cdots\cdots$ &     & $\cdots$ &     & $\cdots\cdots\cdots$ &     & $\cdots$\\ 
$X_{n-1}$ & $=$ & $A_{n-1,0}X_0$ & $+$ & $\cdots$ & $+$ & $A_{n-1,n-1}X_{n-1}$ & $+$ & $E_{n-1}$. 
\end{tabular}
\end{center} 
This system of $n$ right linear equations in $n$ unknowns satisfies the conditions 
of the previous theorem and therefore has a unique solution. Moreover, the 
solution can easily be converted into regular expressions. 
\end{out} 
 
So every non-deterministic automaton gives rise 
to a system of language equations from whose solutions a description of the 
language may be obtained. 
 
\begin{thm}[Kleene's Theorem \cite{Cohen}]
A language $L$ is regular if and only if it is recognisable. 
\end{thm}

 
\subsection{Acceptors for Kan Extensions} 
Throughout this section we will use the notation introduced in Chapter Two.
Recall that a presentation of a Kan extension $(K,\ep)$ 
is a 
quintuple $\mathcal{P}:=kan \lan \Gamma | \Delta | RelB | X | F \ran$  
where $\Gamma$ and $\Delta$ are graphs, $RelB$ is a set of relations on 
$\bP:=P\Delta$, while 
$X:\Gamma \to \sets$ and $F:\Gamma \to \bP$ are graph morphisms.
Elements of the set
$$
T:=\bigsqcup_{B \in \ob\Delta} \bigsqcup_{A \in \ob\Gamma} XA \times \bP(FA,B)
$$
are written $t=x|b_1 \cdots b_n$ with $x \in XA$, 
and $b_1,\ldots,b_n \in \arr\Delta$ are composable with $src(b_1)=FA$.
The function $\tau:T \to \ob\Delta$ is defined by 
$\tau(x|b_1 \cdots b_n):=tgt(b_n)$ and the action of $\bP$ on $T$, written
$t\cdot p$ for $t \in T$, $p\in \arr\bP$, is defined when $\tau(t)=src(p)$.\\

In Chapter Two we defined an initial rewriting system $R_{init}:=(R_\ep,R_K)$ on $T$, and 
gave a procedure for attempting to complete this system. 
We will be assuming that the procedure has terminated, returning a complete 
rewriting system $R=(R_T,R_P)$ on $T$. 
In this section automata will be used to find regular expressions for 
each of the sets $KB$ for $B \in \ob \Delta$. \\
 
Recall that $\sqcup XA$ is the union of the images under $X$ of all the objects 
of $\Gamma$ and $\sqcup KB$ is the union of the images under $K$ of all the 
objects of $\Delta$.
In general the automaton for the irreducible terms which are accepted as members 
of $\sqcup KB$ is the complement 
of the machine which accepts any string containing undefined compositions of 
arrows of $\bB$, any string not containing a single $x_i$ on the left-most end, 
and any string containing the left-hand side of a rule. 
This essentially uses a semigroup presentation of the Kan extension. 
 
\begin{lem} 
Let $\mathcal{P}$ present the Kan extension $(K, \ep)$. 
Then the set $\sqcup KB$ may be identified with the non-zero elements of the 
semigroup having the presentation with generating set 
$$U:= (\sqcup XA) \sqcup \arr\Delta \sqcup 0$$ 
and relations 
\begin{center} 
\begin{tabular}{lll} 
$0u=u0=0$     &  for all & $u \in U$,\\ 
$ux=0$        &  for all & $u \in U, \ x \in \sqcup XA$,\\  
$xb=0$        &  for all & $x \in XA, \ A \in \ob\Gamma, 
                         b \in \arr\Delta$ ~ such that ~ $src(b)\not=FA$,\\ 
$b_1b_2=0$    &  for all & $b_1, b_2 \in \arr\Delta$ ~ such that ~ 
$src(b_2)\not=tgt(b_1)$\\ 
$x(Fa)=(x\cdot a)$ & for all & $x \in XA, \ a \in \arr \bA$ ~ such that ~ $src(a)=A$,\\ 
$l=r$         &  for all & $(l,r) \in RelB$.\\ 
\end{tabular} 
\end{center} 
\end{lem} 
 
\begin{proof} 
The semigroup defined is the set of equivalence classes of  
$T$ with respect to the second two relations (i.e. the Kan extension rules 
$R_\ep$ and $R_K$) 
with a zero adjoined and multiplication of any two classes of $T$ defined to be 
zero. 
\end{proof} 
 
\begin{lem} 
Let $\mathcal{P}$ be a presentation of a Kan extension $(K,\ep)$. 
Then $T$ is a regular language over 
the alphabet $\Sigma:= (\sqcup XA) \sqcup \arr\Delta$. 
\end{lem} 
\begin{proof} 
To prove that $T$ is regular over $\Sigma$ we define an automaton with input 
alphabet $\Sigma$ which recognises $T \subseteq \Sigma^*$. 
Define $\underline{A}:=(S,\Sigma,s_0,\delta,Q)$ where 
$S:=\ob\Delta \sqcup s_0 \sqcup d$, $Q:=\ob\Delta$ and $\delta$ is defined as follows: 
\begin{align*}
\delta(s_0,u):=& 
\left\{ \begin{array}{ll} 
    FA & \quad \text{ for } u \in XA, A \in \ob \Gamma\\ 
    d  & \quad \text{ otherwise.}\\ 
        \end{array} \right. \\
\text{for } B \in \ob \Delta, \quad
\delta(B,u):=&
\left\{ \begin{array} {ll} 
    tgt(u) & \text{ for } u \in \arr\Delta, src(u)=B\\ 
    d & \text{ otherwise.}\\ 
                          \end{array} \right.\\ 
\delta(d,u):= & \quad d \qquad \text{ for all } u \in \Sigma. 
\end{align*} 
It is clear from the definitions that the extended state transition $\delta^*$ 
is such that $\delta^*(s_o,t) \in\ob\Delta$ if and only if $t \in T$. Hence 
$L(\underline{A})=T$. 
\end{proof} 
 
\begin{thm} 
Let $\mathcal{P}$ be a presentation of a Kan 
extension $(K,\ep)$. 
Let $R$ be a finite rewriting system on $T$. 
Then the set of elements $\mathtt{IRR}(\to_R) \subseteq T$ which are irreducible with 
respect to $\to_R$ is a regular language over the alphabet $\Sigma:=\sqcup XA \sqcup 
\arr\Delta$. 
\end{thm} 
\begin{proof} 
We define an incomplete non-deterministic automaton $\underline{A}$ 
 with input alphabet $\Sigma$, and language 
$\Sigma^*-\mathtt{IRR}(\to_R)$ i.e. that rejects only the irreducible elements of $T$ and 
accepts all reducible and undefined elements. This is sufficient proof for the 
theorem, since a language recognised by an incomplete non-deterministic 
automaton $\underline{A}$ is recognisable and 
therefore regular. The complement of $\Sigma^*-\mathtt{IRR}(R)$ is $\mathtt{IRR}(R)$ and therefore 
if $\Sigma^*-\mathtt{IRR}(R)$ is regular then $\mathtt{IRR}(R)$ is regular.\\
 
Begin by defining $L(R_T)$ and $L(R_P)$ to be the sets of left hand sides of 
rules of $R_T$ and $R_P$ respectively. Then define $\mathtt{PL}(R_T)$ and $\mathtt{PL}(R_P)$ to 
be the sets of all prefixes of elements of $L(R_T)$ and $L(R_P)$ and 
define $\mathtt{PPL}(R_T)$ and $\mathtt{PPL}(R_P)$ 
to be the sets of all proper prefixes of elements of $L(R_T)$ and $L(R_P)$.
The proper prefixes of a term $x|b_1 \cdots b_n$ are the terms 
$x|b_1,\ldots,x|b_{n-1}$. Note that each $x$ has its own state and we do not 
require that $x|\id$ is a prefix. Similarly the proper prefixes of a path
$b_1 \cdots b_n$ are the elements $b_1, \ldots b_1 \cdots b_{n-1}$. The 
difference between proper prefixes and prefixes is that $x|b_1 \cdots b_n$ is
considered to be a prefix of itself (but not a proper one), similarly for
$b_1 \cdots b_n$. 
Note $\mathtt{PPL}(R_T) \cup L(R_T) = \mathtt{PL}(R_T)$, similarly for $R_P$.\\ 
 
Define $\underline{A}:=(S,\Sigma,s_0,\delta,Q)$ where 
$S ~ := ~ s_0 \sqcup (\ob\Delta \cup (\sqcup XA) \cup \mathtt{PPL}(R_T) \cup \mathtt{PPL}(R_P)) \sqcup D$, 
$Q ~ := ~ s_0 \sqcup D$.
Let $x,b \in \Sigma$ so that $x \in \sqcup XA$ and $b \in \arr \Delta$.
Let $x_1 \in \sqcup XA$, $B\in\ob\Delta$, $u\in \mathtt{PPL}(R_P)$ and $p\in \mathtt{PPL}(R_P)$.
Define the transition 
$\delta: S \times \Sigma \to \mathbb{P}(S)$ by:
\begin{align*}
\delta(s_0,x) ~ := ~ &\left\{ \begin{array}{ll}
                  \{ x \}                & \text{ if } x \not\in L(R_T),\\
                  \{ D \}              & \text{ if } x \in L(R_T),
                  \end{array} \right.\\  
\delta(s_0,b) ~ := ~ & \{ D \},\\ 
\delta(y,x) ~ := ~ & \{ D \},\\
\delta(y,b) ~ := ~ &\left\{ \begin{array}{ll} 
                  \{x_1|b, tgt(b)\} & \text{ if } x_1|b \in \mathtt{PPL}(R_T),\\ 
                  \{tgt(b)\}            & \text{ if } \tau(y)=src(b), y|b \not\in \mathtt{PL}(R_T),\\ 
                  \{D\}                 & \text{ if } x_1|b \in L(R_T),\\
                  \{D\}                 & \text{ if } \tau(y) \not = src(b),\\
                  \end{array} \right.\\  
\delta(B,x) ~ := ~ & \{D\},\\
\delta(B,b) ~ := ~ & \left\{ \begin{array}{ll} 
                  \{b, tgt(b) \}    & \text{ if } src(b)=B, b \in \mathtt{PPL}(R_P),\\ 
                  \{tgt(b)\}            & \text{ if } src(b)=B, b \not\in \mathtt{PL}(R_P),\\ 
                  \{D\}                 & \text{ if } src(b)=B, b \in L(R_P),\\
                  \{D\}                 & \text{ if } src(b) \not= B,\\
                   \end{array} \right.\\
\delta(u,x) ~ := ~ & \{D\},\\
\delta(u,b) ~ := ~ & \left\{ \begin{array}{ll}
                  \{u \cdot b, tgt(b)\} & \text{ if } u \cdot b \in \mathtt{PPL}(R_T),\\
                  \{tgt(b) \}               & \text{ if } \tau(u)=src(b), u \cdot b \not\in \mathtt{PL}(R_T),\\
                  \{ D \}                     & \text{ if } u \cdot b \in L(R_T),\\
                  \{ D \}                    & \text{ if } \tau(u) \not = src(b),\\
                   \end{array} \right.\\
\delta(p,x) ~ := ~ & \{ D \},\\
\delta(p,b) ~ := ~ &\left\{ \begin{array}{ll} 
                  \{pb, tgt(b)\} & \text{ if } pb \in \mathtt{PPL}(R_P),\\
                  \{tgt(b) \}         & \text{ if } tgt(p)=src(b), pb \not\in \mathtt{PL}(R_P),\\
                  \{D\}              & \text{ if } pb \in L(R_P),\\
                  \{D\}              & \text{ if } tgt(p) \not = src(b),\\
                  \end{array} \right. \\
\delta(D,x) ~ := ~ & \{D\},\\
\delta(D,b) ~ := ~ & \{D\}.
\end{align*}
It follows from these definitions that the extended state transition function 
$\delta^*$ is such that $\delta^*(s_0,t) \cap Q \not= \emptyset$ if and only if $t$ is in 
$\Sigma^*-T$ or if some part of $t$ is the left-hand side of a rule of $R$ (i.e.
if $t$ is reducible).
Therefore $\Sigma^*-\mathtt{IRR}(R)$ is regular, hence $\mathtt{IRR}(R)$ is regular. 
\end{proof} 
 
\begin{cor} 
Let $R$ be a finite complete rewriting system for a Kan extension $(K,\ep)$. 
Then regular expressions for the sets $KB$ of the extended action $K$ can be 
calculated. 
\end{cor} 
\begin{out} 
This follows from the preceding results. 
The automaton $\underline{A}$ of the theorem can be constructed using the 
specifications in the proof. By the results quoted in the introduction to this 
chapter a complete deterministic automaton that 
recognises the same language can be defined. The complement of this has a 
language that can be identified with $\sqcup KB$. Language equations for this 
automaton can be written down and Arden's theorem may be applied to find a 
solution, which gives the language of the automaton as a regular expression. 
\end{out} 
 
The following example illustrates the calculations outlined above. 
 
\begin{example} 
\emph{We construct simple automata which accept the terms which represent 
elements of some set $KB$ for $B \in \ob\bB$ for the general example of a Kan 
extension \ref{kaneg}. 
Recall that the graphs were} 
$$ 
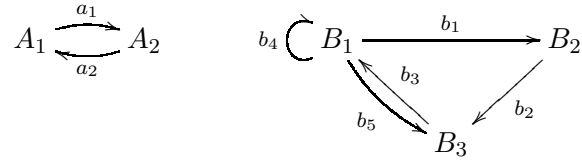
\xymatrix{A_1 \ar@/^/[r]^{a_1} 
        & A_2 \ar@/^/[l]^{a_2} 
       && B_1  \ar@(dl,ul)^{b_4} 
             \ar[rr]^{b_1} 
             \ar@/_/[dr]_{b_5} 
       && B_2 \ar[dl]^{b_2} \\ 
     &&&& B_3 \ar[ul]_{b_3} \\} 
$$ 
\emph{The relations are $RelB= \{ b_1b_2b_3=b_4 \}$, $X$ was defined by  
$XA_1 = \{ x_1, x_2, x_3 \}, XA_2 = \{ y_1, y_2 \}$   with 
$Xa_1:XA_1 \to XA_2: x_1 \mapsto y_1, x_2 \mapsto y_2, x_3 \mapsto 
y_1$, 
$Xa_2:XA_1 \to XA_2: y_1 \mapsto x_1, y_2 \mapsto x_2,$  \, and $F$ was defined by 
$FA_1=B_1$, $FA_2=B_2$, $Fa_1=b_1$ and $Fa_2 = b_2b_3$.}\\ 

{\em The completed rewriting system was:} 
\begin{center} 
\begin{tabular}{llll} 
$x_1|b_1 \to y_1|\id_{B_2}$,      &  $x_2|b_1 \to y_2|\id_{B_2}$, 
  &  $x_3|b_1 \to y_1|\id_{B_2}$, &  $y_1|b_2b_3 \to x_1|\id_{B_1}$,\\ 
$y_2|b_2b_3 \to x_2|\id_{B_1}$,      &  $x_1|b_4 \to x_1|\id_{B_1}$, 
  &  $x_2|b_4 \to x_2|\id_{B_1}$,    &  $x_3|b_4 \to x_1|\id_{B_1}$,\\ 
$b_1b_2b_3 \to b_4$.             &&&\\ 
\end{tabular} 
\end{center} 
\emph{The proper prefix sets are 
$\mathtt{PPL}(R_T):=\{y_1|b_2, y_2|b_2\}$ and 
$\mathtt{PPL}(R_P):=\{b_1,b_1b_2\}$.
The following table defines the incomplete non-deterministic automaton 
which rejects only the terms of $T$ that are irreducible with respect to the 
completed relation $\to$. The alphabet over which the automaton is defined is} 
$\Sigma:=\{x_1,x_2,x_3,y_1,y_2, b_1, b_2, b_3, b_4, b_5 \}$. 
\begin{center}
\begin{tabular}{l|llllllllll}
\emph{state/letter} &$ x_1 $&$ x_2 $&$ x_3 $&$ y_1 $&$ y_2 $&$ b_1 $&$ b_2 $&$ b_3 $&$ b_4 $&$ b_5$\\
\hline
$s_0 $&$ x_1 $&$ x_2 $&$ x_3 $&$ y_1 $&$ y_2 $&$ D $&$ D $&$ D $&$ D $&$ D$\\
$x_1 $&$ D $&$ D $&$ D $&$ D $&$ D $&$ D $&$ D $&$ D $&$ D $&$ B_3$\\
$x_2 $&$ D $&$ D $&$ D $&$ D $&$ D $&$ D $&$ D $&$ D $&$ D $&$ B_3$\\
$x_3 $&$ D $&$ D $&$ D $&$ D $&$ D $&$ D $&$ D $&$ D $&$ D $&$ B_3$\\
$y_1 $&$ D $&$ D $&$ D $&$ D $&$ D $&$ D $&$ y_1|b_2, B_3 $&$ D $&$ D $&$ D$\\
$y_2 $&$ D $&$ D $&$ D $&$ D $&$ D $&$ D $&$ y_2|b_2, B_3 $&$ D $&$ D $&$ D$\\
$y_1|b_2 $&$ D $&$ D $&$ D $&$ D $&$ D $&$ D $&$ D $&$ D $&$ D $&$ D$\\
$y_2|b_2 $&$ D $&$ D $&$ D $&$ D $&$ D $&$ D $&$ D $&$ D $&$ D $&$ D$\\
$B_1 $&$ D $&$ D $&$ D $&$ D $&$ D $&$ b_1, B_2 $&$ D $&$ D $&$ B_1 $&$ B_3$\\
$B_2 $&$ D $&$ D $&$ D $&$ D $&$ D $&$ D $&$ B_3 $&$ D $&$ D $&$ D$\\
$B_3 $&$ D $&$ D $&$ D $&$ D $&$ D $&$ D $&$ D $&$ B_1 $&$ D $&$ D$\\
$b_1    $&$ D $&$ D $&$ D $&$ D $&$ D $&$ D $&$ b_1b_2, B_3 $&$ D $&$ D $&$ D$\\
$b_1b_2 $&$ D $&$ D $&$ D $&$ D $&$ D $&$ D $&$ D $&$ D $&$ D $&$ D$\\
$D      $&$ D $&$ D $&$ D $&$ D $&$ D $&$ D $&$ D $&$ D $&$ D $&$ D$\\
\end{tabular}
\end{center}

\emph{By constructing the transition tree for this automaton, 
we will make it deterministic. The next picture is of the partial transition 
tree -- the arrows to the node marked $\{D\}$ are omitted.} 
$$ 
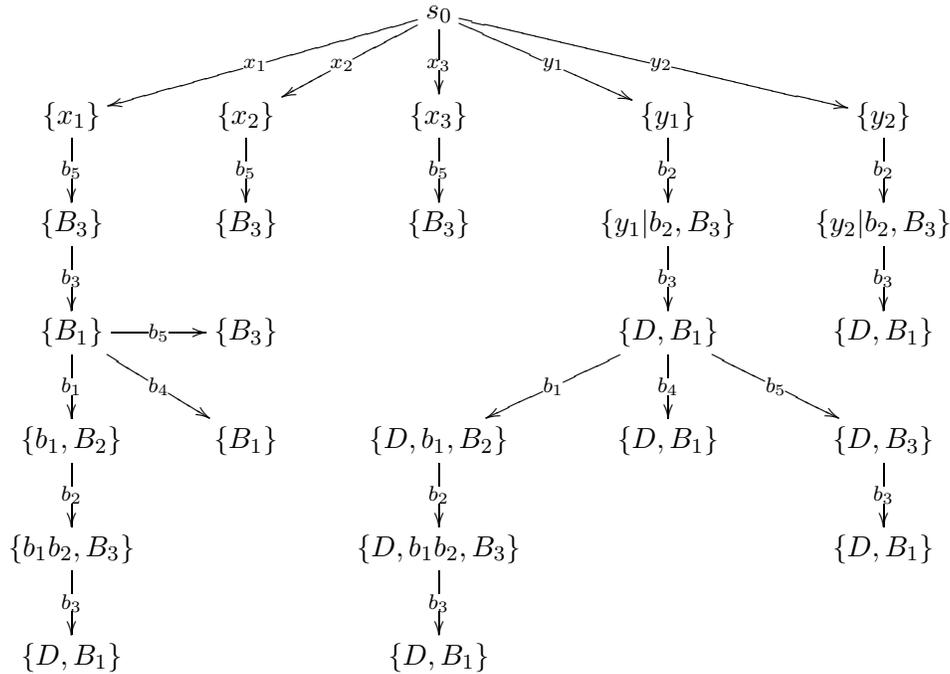
\xymatrix{&&s_0 \ar[dll]|{x_1} \ar[dl]|{x_2} \ar[d]|{x_3} \ar[dr]|{y_1} \ar[drr]|{y_2}&&\\
          \{x_1\} \ar[d]|{b_5} & \{x_2\} \ar[d]|{b_5} & \{x_3\} \ar[d]|{b_5} 
            & \{y_1\} \ar[d]|{b_2} & \{y_2\} \ar[d]|{b_2} \\
          \{B_3\} \ar[d]|{b_3} & \{B_3\} & \{B_3\}
            & \{y_1|b_2,B_3\} \ar[d]|{b_3} & \{y_2|b_2,B_3\} \ar[d]|{b_3}\\
          \{B_1\} \ar[d]|{b_1} \ar[dr]|{b_4} \ar[r]|{b_5} & \{B_3\} &
            & \{D,B_1\} \ar[dl]|{b_1} \ar[d]|{b_4} \ar[dr]|{b_5} & \{D,B_1\}\\
          \{b_1,B_2\} \ar[d]|{b_2} & \{B_1\} & \{D,b_1,B_2\} \ar[d]|{b_2}
            & \{D,B_1\} & \{D,B_3\} \ar[d]|{b_3}\\
          \{b_1b_2,B_3\} \ar[d]|{b_3} && \{D,b_1b_2,B_3\} \ar[d]|{b_3}
            & & \{D,B_1\}\\
          \{D,B_1\} && \{D,B_1\} &&\\}
$$ 
\emph{The tree is constructed with strict observation of the order 
on $\sqcup XA$ and $\arr\Delta$, all arrows are drawn from $\{s_0\}$ and then
arrows from each new state created, in turn. 
When a label e.g. $\{B_3\}$ occurs that branch 
of the tree is continued only if that state has not been defined previously.
Eventually the stage is reached where no new states are defined, all the 
branches have ended. The tree is then converted into an automaton by `gluing' 
all states of the same label.
The initial state is  $\{s_0\}$ and a state is 
terminal if its label contains a terminal state from the original automaton. 
The automaton can often be made smaller, for example, here 
all the terminal states may be glued together.
One possibility is drawn below:}\\
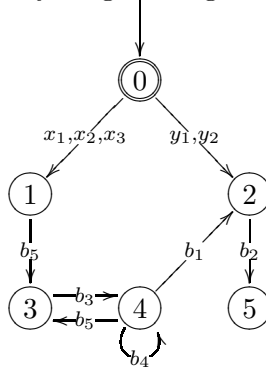
\vspace{-1cm} 
$$ 
\xymatrix{& \ar[d] & \\ 
         & *++[o][F=]{0} \ar[dl]|{x_1,x_2,x_3} \ar[dr]|{y_1,y_2} 
            & \\ 
            *++[o][F-]{1}  \ar[d]|{b_5}
         && *++[o][F-]{2}  \ar[d]|{b_2} \\
            *++[o][F-]{3}  \ar@<1ex>[r]|{b_3} &
            *++[o][F-]{4}  \ar@<1ex>[l]|{b_5} \ar@(dl,dr)|{b_4} \ar[ur]|{b_1}
            & *++[o][F-]{5} \\}
$$ 
\emph{Here the state $S_1$ is labelled $1$ and corresponds to the glueing together of 
$\{x_1\}$, $\{x_2\}$ and $\{x_3\}$ to form $\{x_1,x_2,x_3\}$ and the state $S_2$
is $\{y_1,y_2,b_1,B_2\}$.
States $S_3$ and $S_4$ represent $\{B_3\}$ and  $\{B_1\}$
respectively and state $S_5$ is $\{y_1|b_2,y_2|b_2,B_3,b_1b_2\}$.
The complement of this automaton accepts all irreducible elements of 
$\sqcup KB$. 
When $S_1$ and $S_4$ are terminal the language accepted is $KB_1$. 
When $S_2$ is terminal the language accepted is $KB_2$. 
When $S_3$ and $S_5$ are terminal the language accepted is $KB_3$. 
The language equations from the automaton for $KB_1$ are:}
\vspace{-0.2cm} 
\begin{align*} 
X_0 & =(x_1+x_2+x_3)X_1+(y_1+y_2)X_2,\\ 
X_1 & =b_5X_3+\id_{B_1},\\ 
X_2 & =b_2X_5,\\ 
X_3 & =b_3X_4,\\ 
X_4 & =b_1X_2+b_4X_4+b_5X_3+\id_{B_1},\\ 
X_5 & =\emptyset.
\intertext{\emph{Putting $X_2=\emptyset$ and eliminating $X_1$ and $X_3$ 
by substitution gives}}
X_0 &= (x_1+x_2+x_3)(b_5b_3X_4+\id_{B_1}),\\
X_4 &=(b_4+b_5b_3)X_4 + \id_{B_1}.
\intertext{\emph{Finally, applying Arden's Theorem to $X_4$
we obtain the regular expression}}
X_0 &=(x_1+x_2+x_3)|(b_5b_3(b_4+b_5b_3)^*+\id_{B_1}).
\intertext{\emph{The separator ``$|$'' may be added at this point. 
Similarly, we can obtain regular expressions for $KB_2$ and $KB_3$. 
For $KB_2$ we have}}
X_0 &=(x_1+x_2+x_3)|b_5b_3(b_4+b_5b_3)^*b_1+(y_1+y_2)|\id_{B_2}.
\intertext{\emph{For $KB_3$ the expression is}}
X_0 &=(x_1+x_2+x_3)|(b_5b_3(b_4+b_5b_3)^*(b_1b_2+b_5) + b_5) + (y_1+y_2)|b_2.
\end{align*} 
\end{example} 
 
\subsection{Accepting Automata for Algebras} 
 
We have discussed automata for rewriting systems which accept only irreducible 
words. The concept will now be generalised to Gr\"obner bases. 
The irreducibles of an algebra $K[S]/\lan P\ran $ in which we are interested are 
the 
irreducible monomials; reducibility of a polynomial is determined by 
reducibility of the monomials it contains. 
Therefore the automaton we draw is over the alphabet $X$, the generators of $S$ 
and the language it accepts is the set of irreducible monomials. 
The automaton below is for the infinite dimensional algebra 
$\mathbb{Q}[\{a,b\}^\dagger]$ factored by the ideal generated by the Gr\"obner 
basis $\{a^3-b+2, \ ba^2b-2b^2+4a \}$. 
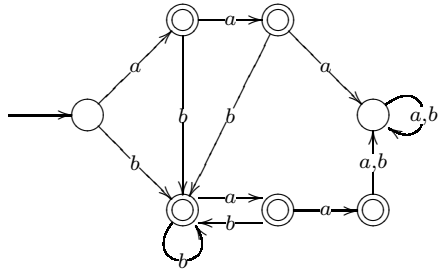
$$ 
\xymatrix{&& *++[o][F=]{} \ar[r]|a \ar[dd]|b 
           & *++[o][F=]{} \ar[dr]|a \ar[ddl]|b &\\ 
          \ar[r] & *++[o][F-]{} \ar[ur]|a \ar[dr]|b 
           &&& *++[o][F-]{} \ar@(ur,dr)|{a,b}\\ 
          && *++[o][F=]{} \ar@<1ex>[r]|a \ar@(dl,dr)|b 
           & *++[o][F=]{} \ar[r]|a \ar@<1ex>[l]|b 
           & *++[o][F=]{} \ar[u]|{a,b} \\} 
$$ 
The point of drawing acceptor automata is to find nice expressions for the sets 
of 
irreducibles. If an algebra is finite then the number of irreducible monomials 
it has is the \emph{dimension} of the algebra. In the infinite example above we 
can at least find a regular expression for the set of irreducible monomials. 
 
It is: 
$$(a^2b+ab+b)(ab+b)^*(a^2+a+\id)+(a^2+a)$$ 
Any element of the algebra is then uniquely expressible as a sum of 
$K$-multiples 
of these monomials. 
 
It is possible to adapt the automaton so that it accepts polynomials by allowing 
$+$ and $-$ to be elements of the input alphabet, with transitions (from each 
state) labelled by $+$ and $-$ going to the initial state, and by adding $k$ for 
$k \in K$ as a loop at the initial state. In this way it may be possible to 
define automatic algebras. One difficulty to such a definition is the fact that 
a multiplier/equality recogniser has to recognise that two polynomials are equal 
though the terms may be input in a different order ($b+a^2$ and $a^2+b$). 
There is not the option, as with the acceptor, of working only with monomials. 
The reason for this is that the normal form of a monomial $w$ multiplied by 
a generator $x$ (as if to define the multiplier automaton) may well not be 
a monomial. We mention these issues in passing, only here being concerned 
with the acceptors and with the reduction machines (next section).

\section{Reduction Machines} 
 
\subsection{Cayley Graphs}

The Cayley graph $\Gamma$ of a group $G$ with generating set $X$, and quotient 
morphism $\theta:F(X) \to G$  is the graph 
with vertex set $\ob\Gamma:=G$ and edge set $\arr\Gamma:=G\times X$ with 
$src[g,x]=g$, $tgt[g,x]=g \theta(x)$. 
The Cayley graph is a representation of the whole 
multiplication table for the group. In this section we indicate how to use the 
Cayley graph of a group to help with rewriting procedures. The results are not 
surprising, but formalise certain procedures which may sometimes be useful. 
 
\begin{prop} 
Let $G$ be the group given by the finite presentation $grp\lan X|Rel\ran$. Let 
$\Gamma$ be the Cayley graph of $G$. 
Let $\theta:F(X)\to G$ be the quotient map. 
Let $>$ be the length-lex order on $X^*$ induced by a linear order on $X$. 
Then $>$ specifies a tree in the Cayley graph and a vertex labelling $V\subseteq 
X^*$ where for all $w_1\in V$, $w_2\in F(X)$ such that $\theta(w_1)=\theta(w_2)$ 
it is the case that $w_2>w_1$ or $w_2=w_1$. 
\end{prop} 
\begin{proof} 
Since $G$ is finite the inverse of any generator can be represented by a 
positive power. So for any word $r \in F(X)$ there is a word $r^+$ obtained
by replacing each $x^{-1}$ with $x^{Order(x)-1}$, with $\theta(r)=\theta(r^+)$.
Therefore we consider the presentation $mon\lan X|R\ran$ where 
$R:=\{(r^+,\id):r\in Rel\}$ of $G$. 
Let $T:=\emptyset$, $V:=\emptyset$. 
Start at vertex $\id$ and add this label to $V$. 
Go through the elements of $X$ in order, adding the edge $[\id,x]$ to $T$ 
whenever it will not create a cycle in the graph. When an edge $[\id,x]$ is 
added to $T$ the target vertex label $x$ should be added to $V$. 
Clearly, if $x_i\in V$ and $\theta(x_i)=\theta(x_j)$ for some $x_j$ in $X$ then 
$x_j>x_i$ and $x_j\not\in V$ or else $x_j=x_i$.\\ 
 
Now repeat the following step until all the vertices of the graph are  
represented in $V$; that is until $\theta(V)=G$. Choose the vertex with least  
label $w$ of $V$ in the graph and go through the elements of $X$ in order adding  
edges $[w,x]$ to $T$ whenever $\theta(wx)\not\in \theta(V)$. This is the  
condition that to add that edge will not create a cycle. For each new edge 
$[w,x]$ added to $T$, add the vertex label $wx$ to $V$.\\ 
It is immediate from the inductive construction that the set of vertex labels 
$V$ is least in the sense that for any $w$ in $V$, $w$ is the least element of 
$F(X)$ with respect to $>$ with image $\theta(w)$. 
Furthermore, since $\Gamma$ is connected and edges are chosen so as not to 
create cycles, $T$ defines a spanning tree of $\Gamma$ with edges 
$[\theta(w),x]$. 
\end{proof} 
 
\begin{cor} 
The set of vertex labels $V$ is a set of unique normal forms for $G$ in $F(X)$  
and the tree $T$ defines a normal form function $N:F(X)\to V$. 
\end{cor} 
\begin{proof} 
It is immediate from the last result that $V$ is a set of unique normal forms 
for $R$ on $X^*$. The normal form function is defined by using the Cayley graph 
as a reduction machine operating on $F(X)$. 
Let $x_0^{\ep_0}x_1^{\ep_1}\cdots x_m^{\ep_m}$ be an input word where 
$\ep_i:=\pm 1$ and $x_i \in X$. 
Start at the vertex with label $\id$ and follow the path  
$[\id,x_0^{\ep_0}][\theta(x_0^{\ep_0}),x_1^{\ep_1}] 
\cdots [\theta(x_0^{\ep_0}\cdots x_{m-1}^{\ep_{m-1}}),x_m^{\ep_m}]$. 
The label of the target vertex  
$\theta(x_0^{\ep_0}x_1^{\ep_1}\cdots x_m^{\ep_m})$ is the least element  
$w \in F(X)$ such that $\theta(w)=\theta(x_0^{\ep_0}  x_1^{\ep_1} \cdots x_m^{\ep_m})$. 
This defines a normal form function $N$.
\end{proof} 
 
\begin{example} 
\emph{ Consider the Cayley graph for the dihedral group $D_8$ which is presented by  
\linebreak $grp\lan a,b|a^4, b^2, abab\ran$. 
The Cayley graph is depicted below, with the vertices labelled according to the 
ordering induced by $a<b$.} 

{\Large
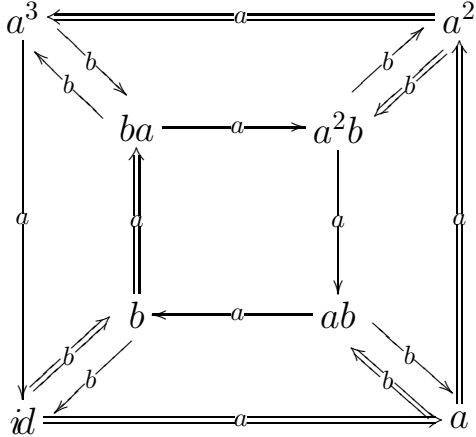
$$ 
\xymatrix{a^3 \ar[dddd]|a \ar@<1ex>[dr]|b 
            &&&& a^2 \ar@{=>}[llll]|a \ar@{=>}@<1ex>[dl]|b \\ 
          & ba \ar[rr]|a \ar@<1ex>[ul]|b  
            && a^2b \ar[dd]|a \ar@<1ex>[ur]|b & \\ 
          &&&& \\ 
          & b \ar@{=>}[uu]|a \ar@<1ex>[dl]|b 
            && ab \ar[ll]|a \ar@<1ex>[dr]|b & \\ 
          \id \ar@{=>}[rrrr]|a \ar@{=>}@<1ex>[ur]|b 
            &&&& a \ar@{=>}[uuuu]|a \ar@{=>}@<1ex>[ul]|b \\} 
$$}
\emph{Consider the word $aba^3b$. 
Beginning at $\id$ follow the path to $a$. 
Read $b$ and go to vertex $ab$. 
Read $a$ and so go to vertex $b$. 
When the final $b$ is read, it takes us to the vertex with label $a^2$, hence 
$N(aba^3b)=a^2$.}
\end{example}

\subsection{Reduction Machines for Kan Extensions} 
 
We now generalise the reduction machine idea to Kan extensions. 
Formally, standard output automata are defined in two ways, as \emph{Moore} 
machines or \emph{Mealy} machines (see \cite{Hopcroft}). 
The reduction machines here are Moore machines.\\ 
 
A \textbf{Moore machine} is a six-tuple 
$\underline{M}:=(S,\Sigma, s_0, \delta,\lambda, \Theta)$ 
where $S$ is the set of states with an initial state $s_0$, 
 $\Sigma$ is the input alphabet,  $\Theta$ is the output alphabet,  
 $\delta$ is the transition function from $S \times \Sigma \to S$ 
and $\lambda:S \to \Theta$ is a mapping which gives the output associated with 
each state. (All states are ``terminal''.)
As before $\dd^*$ denotes the extended state transition function.\\

We continue with the assumption that
$\mathcal{P}:=kan \lan \Gamma | \Delta | RelB | X | F \ran$ is the finite
presentation of the Kan extension $(K,\ep)$ and $R=(R_T,R_P)$ is a finite
complete rewriting system on the $\bP$-set $T$ given by $\mathcal{P}$.
We will only work with finite machines, so for the rest of this chapter 
the Kan extensions will be assumed to be finite i.e. $\sqcup KB$ is finite.

\begin{prop} 
Let $\mathcal{P}$ be a presentation of a finite Kan extension, with complete
rewriting system $R$.
Then there exists a Moore machine 
$\underline{M}=(S,\Sigma,s_0,\delta,\lambda,\Theta)$ such that 
$\lambda(\delta(w))$ is the irreducible form of $w$ with respect to $\to_R$ on  
$T$.  
\end{prop} 
\begin{proof} 
Define a Moore machine $M$ in the following way. 
Let $S := (T/{\stackrel{*}{\lra}_R}) \sqcup s_0 \sqcup d$, 
$\Sigma := XA\sqcup\arr\Delta$, and
$\Theta := T \sqcup 0$. Let $s_0$ be the initial state. 
Define $\delta:S\times T \to S$  
by $\delta(s_0,x) := [x|\id_{FA}]$ and 
   $\delta([t],x) = \delta(d,x) := d$ for all 
                     $x\in XA, A\in\ob\Gamma$ and $t\in T$; and 
   $\delta([t],b) := [t\cdot b]$ for all 
                     $t\in T,\ b\in\arr\Delta$ such that $\tau(t) = src(b)$ and 
   $\delta([t],b) = \delta(s,b) = \delta(d,b) := d$ otherwise. 
Then define $\lambda:S\to\Theta$ by $\lambda(s)=\lambda(d)=0$ and 
$\lambda([t]) := N(t)$. 
It is clear from these definitions that $\lambda(\delta(s,t))=N(t)$ for all 
$t\in T$.  
\end{proof}

\begin{example}
\emph{We conclude this subsection with an example of a 
reduction machine for a Kan extension. Let 
$\mathcal{P}$ be a Kan extension where 
$\Gamma$ and $\Delta$ are as follows:} 
$$
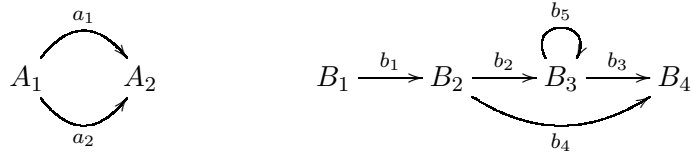
\xymatrix{A_1 \ar@/^1.5pc/[r]^{a_1} \ar@/_1.5pc/[r]_{a_2} & A_2
   && B_1 \ar[r]^{b_1} & B_2 \ar[r]^{b_2} \ar@/_1.5pc/[rr]_{b_4} &
        B_3 \ar@(ul,ur)^{b_5} \ar[r]^{b_3} & B_4 \\} 
$$ 
\emph{The relations of $\bB$ are } $RelB:=\{(b_2b_5b_3,b_4),(b_5^2,b_5)\}$.
\emph{The functors $F$ and $X$ are defined by:-} 
$FA_1:=B_1$, $FA_2:=B_4$, $Fa_1:=b_1b_2b_3$, $Fa_2:=b_1b_4$ {\em and }
$XA_1:=\{x_1,x_2,x_3\},\ XA_2:=\{y_1,y_2\}$,
$Xa_1:XA_1\to XA_2: x_1\mapsto y_1, x_2\mapsto y_1, x_3\mapsto y_2$, \, 
$Xa_2:XA_1\to XA_2: x_1\mapsto y_1, x_2\mapsto y_2, x_3\mapsto y_2$. 
\emph{The initial rewriting system is in fact complete. It is} 
\begin{center}
$\{ x_1|b_1b_2b_3 \to y_1|\id_{B4},  
\ x_2|b_1b_2b_3 \to y_1|\id_{B4},  
\ x_3|b_1b_2b_3 \to y_2|\id_{B4},  
\ x_1|b_1b_4 \to y_1|\id_{B4},$\\  
$\ x_2|b_1b_4 \to y_2|\id_{B4},  
\ x_3|b_1b_4 \to y_2|\id_{B4},  
\ b_2b_5b_3 \to b_4,
\ b_5^2 \to b_5 \}.$ 
\end{center}
\emph{Following the directions in the proof above we construct the Moore machine. 
There are 14 states $[t]\in S$ and also the initial state $s$ and the dump 
state $d$ which rejects any terms that are not defined in $T$.} 
\begin{center}
$\lambda(S) := \{d, x_1|\id_{B1}, x_2|\id_{B1}, x_3|\id_{B1}, y_1|\id_{B4},  
y_2|\id_{B4}, $\\
$x_1|b_1, x_2|b_1, x_3|b_1, x_1|b_1b_2, x_2|b_1b_2,  
x_3|b_1b_2, x_1|b_1b_2b_5, x_2|b_1b_2b_5, 
x_3|b_1b_2b_5 \}.$
\end{center}
\emph{The non-trivial part of the transition function is as follows:}
\begin{alignat*}{2} 
\delta(s,x_1)=[x_1|\id_{B1}] &   \quad
\delta(s,x_2)=[x_2|\id_{B1}] &&  \quad
\delta(s,x_3)=[x_3|\id_{B1}] \\  
\delta(s,y_1)=[y_1|\id_{B4}] &  \quad
\delta(s,y_2)=[y_2|\id_{B4}] &&  \quad
\delta([x_1|\id_{B1}],b_1)=[x_1|b_1] \\ 
\delta([x_2|\id_{B1}],b_1)=[x_2|b_1] &\quad
\delta([x_3|\id_{B1}],b_1)=[x_1|b_1] &&\quad
\delta([x_1|b_1],b_2) = [x_1|b_1b_2] \\  
\delta([x_1|b_1],b_4) = [y_1|\id_{B4}] & \quad
\delta([x_2|b_1],b_2) = [x_2|b_1b_2] &&  \quad
\delta([x_2|b_1],b_4) = [y_2|\id_{B4}] \\ 
\delta([x_3|b_1],b_2) = [x_3|b_1b_2] & \quad
\delta([x_3|b_1],b_4) = [y_2|\id_{B4}] && \quad
\delta([x_1|b_1b_2],b_3)=[y_1|b_1] \\
\delta([x_1|b_1b_2],b_5)=[y_1|b_2] & \quad
\delta([x_2|b_1b_2],b_3)=[y_2|b_1] && \quad
\delta([x_2|b_1b_2],b_5)=[y_2|b_2] \\ 
\delta([x_3|b_1b_2],b_3)=[y_1|b_1] & \quad
\delta([x_3|b_1b_2],b_5)=[y_2|b_2] &&
\end{alignat*} 
\emph{The machine can be represented by a diagram -- 
states have not been circled as the labels are too long, and the state $d$ which
rejects anything not defined is not drawn.}

\vspace{-0.8cm}
\large
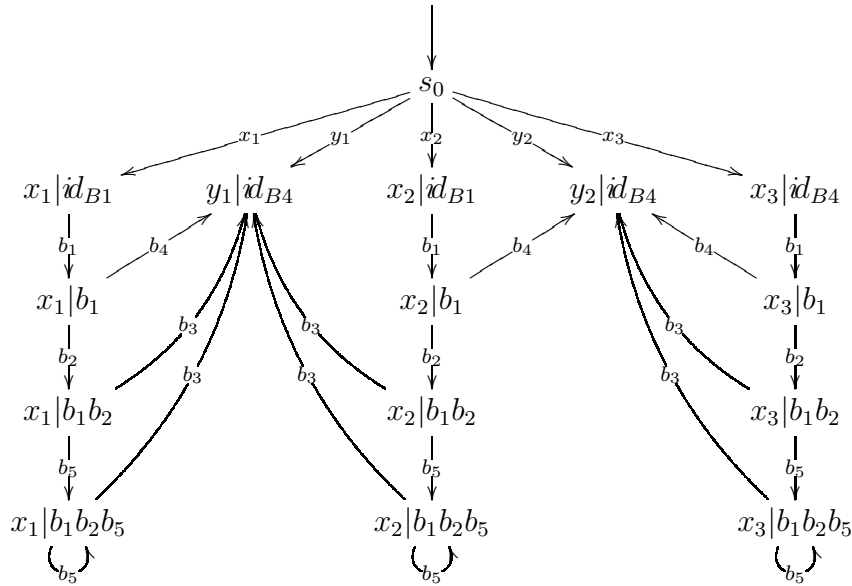
{ 
$$
\xymatrix{ && \ar[d] &&\\
&& {s_0} 
           \ar[dll]|{x_1}
           \ar[dl]|{y_1}
           \ar[d]|{x_2}
           \ar[dr]|{y_2}
           \ar[drr]|{x_3} &&\\
    {x_1|\id_{B1}} \ar[d]|{b_1} 
  & {y_1|\id_{B4}} 
  & {x_2|\id_{B1}} \ar[d]|{b_1}
  & {y_2|\id_{B4}} 
  & {x_3|\id_{B4}} \ar[d]|{b_1} \\
{x_1|b_1}   \ar[d]|{b_2} \ar[ur]|{b_4}
 && {x_2|b_1} \ar[d]|{b_2} \ar[ur]|{b_4}
 && {x_3|b_1} \ar[d]|{b_2} \ar[ul]|{b_4} \\
{x_1|b_1b_2} \ar[d]|{b_5} \ar@/_1.2pc/[uur]|{b_3}
 && {x_2|b_1b_2} \ar[d]|{b_5} \ar@/^1.2pc/[uul]|{b_3}
 && {x_3|b_1b_2} \ar[d]|{b_5} \ar@/^1.2pc/[uul]|{b_3} \\
  {x_1|b_1b_2b_5} \ar@/_1.2pc/[uuur]|{b_3} \ar@(dl,dr)|{b_5}
 && {x_2|b_1b_2b_5} \ar@/^1.2pc/[uuul]|{b_3} \ar@(dl,dr)|{b_5}
 && {x_3|b_1b_2b_5} \ar@/^1.2pc/[uuul]|{b_3} \ar@(dl,dr)|{b_5}
 \\}
$$}
\end{example}
 
This example serves to illustrate the principle of  
converting a complete rewriting system $R$ on $T$ for which there are a finite 
number of irreducibles 
into a machine which accepts terms of $T$ (which may be 
infinite) and gives as output their irreducible form i.e. 
representatives of elements of $\sqcup KB$.

\subsection{Reduction Machines for Algebras} 
 
We have shown how to use general rewriting systems to construct automata. 
In a similar way Gr\"obner bases may be used to construct reduction 
machines for finite dimensional algebras. 
The concepts of reduction machines for the previous structures were new but 
based on standard automata for semigroups. The Gr\"obner reduction 
machines for algebras are different from basic output automata. \\

Let $K$ be a field and let $X$ be a set.
Let $\to_R$ be a reduction relation on $K[X^\dagger]$.
We define a reduction machine $\underline{M}$ to be a marked graph whose
vertices $V$ are labelled by monomials of $X^*$ that are irreducible
with respect to $\to_R$. (The monoid identity $\id$ represents the algebra 
identity $1$.) 
Edges have the form $(c,x)$ with $c \in K$, $x \in X$ and from every vertex 
$m$ there will be at least one edge $(c,x)$ for each $x \in X$. 
The targets of these edges are the monomials of the reduced form
of $mx$ with respect to $\to_R$.\\

A state of the machine can be represented by a vector in 
$K[X^\dagger]^n$, where $n$ is the number of vertices. The value at each vertex
represents the unprocessed input.
When the Cayley graph machines were considered in this way, 
the state of a machine was essentially a function $V \to F(X)$. 
Thus it seems reasonable that the state of a Gr\"obner machine should be 
represented by a function $V \to K[X^\dagger]$. 
Essentially the state of a machine is the 
specification of a value $v \in K[X^\dagger]$ for each vertex $m$.\\

The machine acts by reading the first letter $x_1 \in X$ of a monomial 
$x_1 \cdots x_n$ of the value $v$ at a vertex $m$ and moves to a 
new state determined by all the edges leaving $m$ that are labelled $(c_i,x_1)$
and have target $m_i$.
The value at $m$ is decreased by $kx_1 \cdots x_m$ where 
$k$ is the coefficient of $x_1 \cdots x_n$ in $v$ and the value at each $m_i \in S$
is increased by $c_ix_2 \cdots x_n$.  
The vital difference between these machines and earlier ones is 
that monomials can reduce to polynomials, and so 
there may be more than one arrow with the same letter label coming from a  
vertex. This becomes clearer on examination of an example.

\begin{example} 
\emph{The third Hecke algebra is 
$\mathbb{Q}[\{e_1,e_2\}^*]/ \lan P \ran $ where} 
$$
P:=\{ e_1^2-e_1, e_2^2-e_2, e_2e_1e_2 - e_1e_2e_1 + 2/9\,e_2 - 2/9\,e_1 \}.
$$
\emph{In fact $P$ is a Gr\"obner basis for this algebra. 
The algebra has dimension 6, the irreducible monomials being $\id, e_1, e_2, 
e_1e_2, e_2e_1, e_1e_2e_1$. 
We draw a machine which acts to reduce polynomials in $\mathbb{Q}[\{e_1,e_2\}^*]$  
The edges have two labels; a  
generator $e_1$ or $e_2$ and a coefficient from $\mathbb{Q}$, 
(1 where unmarked).
For example $e_1e_2e_1e_2$ reduces to 
$e_1e_2e_1-\frac{2}{9}e_1e_2+\frac{2}{9}e_1$ 
so there are three arrows with letter label $e_2$ coming out of the vertex 
$e_1e_2e_1$.}\\ 

{\em The following diagram shows the ``Gr\"obner machine'' for the Hecke algebra defined 
above.} 
{\large
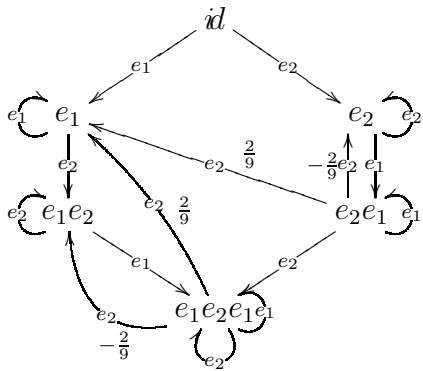
$$ 
\xymatrix{& \id \ar[dl]|{e_1} \ar[dr]|{e_2} &\\ 
          e_1 \ar@(dl,ul)|{e_1}  \ar[d]|{e_2}  
           && e_2 \ar@(dr,ur)|{e_2}  \ar@<1ex>[d]|{e_1} \\ 
          e_1e_2 \ar[dr]|{e_1}  \ar@(dl,ul)|{e_2}  
           && e_2e_1 \ar@(dr,ur)|{e_1}  \ar[ull]|{e_2}_(.4){\frac{2}{9}} 
                     \ar[dl]|{e_2} \ar@<1ex>[u]|{e_2}^{-\frac{2}{9}}\\ 
          & e_1e_2e_1 \ar@/^2pc/[ul]|{e_2}^(.4){-\frac{2}{9}} \ar@(dr,ur)|{e_1}  
                      \ar@/_/[uul]|{e_2}_(.4){\frac{2}{9}} \ar@(dr,dl)|{e_2} &  
\\} 
$$} 
 
\emph{The machine operates to reduce monomials, for example: $e_1e_2e_1e_2e_1$.
Start with the value $e_1e_2e_1e_2e_1$ at vertex $\id$. Read $e_1$ and the 
new state of the machine is given by the value $e_2e_1e_2e_1$ at $e_1$ and $0$ 
elsewhere. Read $e_2$ and the state is now given by the value $e_1e_2e_1$ at 
$e_1e_2$ and $0$ elsewhere. Read $e_1$ and the state of the machine is $e_2e_1$
at $e_1e_2e_1$ and $0$ elsewhere. Read $e_2$ and the new state is given by $e_1$
at $e_1e_2e_1$, $-2/9e_1$ at $e_1e_2$ and $2/9e_1$ at $e_1$ with $0$ elsewhere.
At vertex $e_1e_2e_1$ read $e_1$ and the new state of the machine is $1$ at 
$e_1e_2e_1$ and the values of the other vertices unchanged. At vertex
$e_1e_2$ read $-2/9e_1$ and the new state of the machine is given by 
$7/9$ at $e_1e_2e_1$ and $2/9e_1$ at $e_1$ and $0$ elsewhere. 
To finish, read $2/9e_1$ at $e_1$, and the final state of the machine is given 
by the values of $7/9$ at state $e_1e_2e_1$, $2/9$ at $e_1$ and $0$ elsewhere.
The output polynomial is therefore $7/9e_1e_2e_1 + 2/9e_1$, this is the
irreducible form of $e_1e_2e_1e_2e_1$.}
\end{example} 

The ``Gr\"obner Machines'' described are really no more than ``pictures'' of the 
Gr\"obner bases. We will formalise the ideas of reduction machines for 
algebras, for the general case, by using Petri nets. 
 
\section{Petri nets} 
 
This section introduces Petri nets and formalises the ``Gr\"obner machines'' 
devised in the previous section in terms of these well-defined structures. 

\subsection{Introduction to Petri nets} 
 
Petri nets are a graphical and mathematical modelling tool applicable to many 
systems. They may be used for specifying information processing systems that 
are concurrent, asynchronous, distributed, parallel, non-deterministic, 
and/or stochastic. 
Graphically, Petri nets are useful for illustrating and describing systems, and  
tokens can simulate the dynamic and concurrent activities. 
Mathematically, it is possible to set up models such as state equations and 
algebraic equations which govern the behaviour of systems. 
Petri nets are understood by practitioners and theoreticians and so provide a 
powerful link of communication between them. 
For example engineers can show mathematicians how to make practical and 
realistic models, and mathematicians may be able to produce theories to make 
the systems more methodical or efficient. 
A good introduction to the ideas of Petri nets is \cite{Murata}.\\ 
 
An integer-valued \emph{Petri net} is a kind of directed graph together with an initial state 
(called an \emph{initial marking} $M_0$). 
The underlying graph of a Petri net is a directed, weighted bipartite graph. 
The two kinds of vertices are \emph{places} (represented by circles) and 
\emph{transitions} (represented by rectangles). Edges go between places and 
transitions and are labelled with their weights. A \emph{marking} assigns a 
non-negative integer to each place. If a place $p$ is assigned $k$ in a marking 
then we say $p$ has $k$ \emph{tokens} (represented by black dots). 
In modelling, places represent conditions and transitions represent events. 
A transition has input and output places, which represent preconditions and 
postconditions (respectively) of the event.\\

A \textbf{Petri net} (without specific initial marking) is a 4-tuple 
$\underline{N}=(P,T,\cF,w)$ where:\\ 
$P=\{p_1,\ldots,p_m\}$ is a finite set -- the places,\\ 
$T=\{t_1,\ldots,t_n\}$ is a finite set -- the transitions,\\ 
$\cF \subseteq (P \times T) \cup (T \times P)$ is a set of edges -- the flow 
relation,\\ 
$w:\cF \to \mathbb{N}$ is a weight function,\\ 
and $P \cap T = \emptyset$, $P \cup T \not= \emptyset$.\\ 
 
The state of a Petri net is represented by a marking.
A \textbf{marking}  is a function $M:P \to \mathbb{N}+\{0\}$.
Let $\underline{N}$ be a Petri net where each place is given a 
distinct label $p_i$. To every marking $M$ we will associate a polynomial 
$pol(M):=\Sigma_P \; p M(p)$ that is the formal sum of terms where $M(p)$ is 
a non-negative integer and $p$ is a place label.\\

The behaviour of dynamic systems may be described in terms of system states and 
changes. A marking of a Petri net is changed according to the 
\textbf{firing rule}: 
\begin{enumerate}[i)] 
\item 
A transition $t$ is \textbf{enabled} if each input place $p$ of $t$ is marked 
with at least $w(p,t)$ tokens where $w(p,t)$ is the weight of the edge from $p$ to 
$t$. 
\item 
An enabled transition may or may not \textbf{fire} -- depending on whether or 
not 
the relevant event occurs. 
\item 
Firing of an enabled transition $t$ removes $w(p,t)$ tokens from each 
input place $p$ of $t$ and adds $w(t,q)$ tokens to each output place $q$ of $t$ 
where $w(t,q)$ is the weight of the edge from $t$ to $q$. 
\end{enumerate} 
\begin{example} 
\emph{The markings of the nets below are given by the polynomials $H_2+2O_2$ and
$2H_2+2O_2$ respectively. The transition $t$ is enabled in the second case 
and not in the first:} 
$$
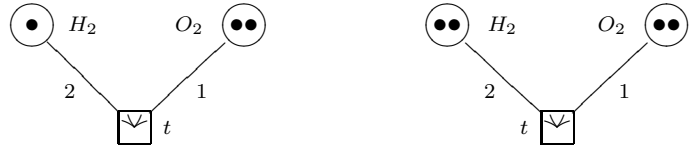
\xymatrix{*++[o][F-]{\bullet} \ar@{}[r]|{H_2} \ar[dr]_2 && 
             *++[o][F-]{\bullet \bullet} \ar@{}[l]|{O_2} \ar[dl]^1 && 
             *++[o][F-]{\bullet \bullet} \ar@{}[r]|{H_2} \ar[dr]_2 && 
             *++[o][F-]{\bullet \bullet} \ar[dl]^1 \ar@{}[l]|{O_2} \\ 
            & *++[F-]{} \ar@{}[r]|(.3){t} 
             &&&& *++[F-]{} \ar@{}[l]|(.3){t} & \\} 
$$ 
\end{example} 
 
Each transition $t$ has an associated polynomial 
$pol(t):= \Sigma_P \; p w(p,t) - \Sigma_P \; p w(t,p)$, that is the sum of the
weights of tokens that a firing of transition $t$ takes from each input place
minus the sum of weights of tokens that it adds to each output place.
A firing/occurrence sequence is denoted by  
$M_0 \stackrel{t_1}{\to} M_1 \stackrel{t_2}{\to} \cdots \stackrel{t_n}{\to} M_n$ where 
the $M_i$ are markings and the $t_i$ are transitions (events) transforming 
$M_{i-1}$ into $M_i$. 
For $i=1, \ldots, n$ it follows from the definitions that 
$pol(M_{i}) = pol(M_{i-1}) - pol(t_i)$.
Therefore the above firing sequence gives the information
$pol(M_n) = pol(M_0) - pol(t_1) - pol(t_2) - \cdots - pol(t_n)$.

\begin{example}
\emph{The formula $2H_2+O_2=2H_2O$ is represented by the transition in the 
diagrams below, the left diagram shows the initial marking and  
the right shows the marking after the transition has fired.} 
$$
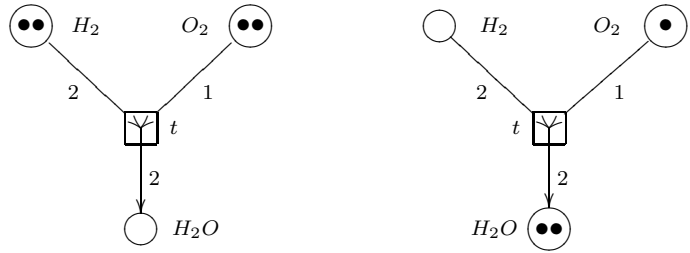
\xymatrix{*++[o][F-]{\bullet \bullet} \ar@{}[r]|{H_2} \ar[dr]_2 && 
             *++[o][F-]{\bullet \bullet} \ar@{}[l]|{O_2} \ar[dl]^1 && 
             *++[o][F-]{} \ar@{}[r]|{H_2} \ar[dr]_2 && 
             *++[o][F-]{\bullet} \ar[dl]^1 \ar@{}[l]|{O_2} \\ 
            & *++[F-]{} \ar@{}[r]|(.3){t} \ar[d]^2  
             &&&& *++[F-]{} \ar@{}[l]|(.3){t} \ar[d]^2 &\\ 
            & *++[o][F-]{} \ar@{}[r]|{H_2O} 
             &&&& *++[o][F-]{\bullet \bullet} \ar@{}[l]|{H_2O} &\\} 
$$ 
\emph{The polynomial for the transition is
$2H_2+O_2-2H_2O$ and the firing sequence would be denoted  
$2H_2+2O_2 \stackrel{t}{\to} O_2+2H_2O$.}
\end{example} 
 
One of the main problems in Petri net theory is \emph{reachability} 
(see \cite{Desel} for some examples).  
A marking $M$ is said to be \textbf{reachable} from a marking $M_0$ in a net  
$\underline{N}$, 
if there is a sequence of firings that transforms $M_0$ to $M$.

\begin{Def}
The \textbf{reachability problem} for a Petri net $\underline{N}$ is as follows:
\begin{center}
\begin{tabular}{ll}
INPUT:    & $M_1$, $M_2$, two markings of $\underline{M}$,\\
QUESTION: & is $M_2$ reachable from $M_1$?\\
\end{tabular}
\end{center}
\end{Def}

Often a Petri net comes with a specified \textbf{initial marking} $M_0$.
Then the reachability refers to reachability from $M_0$ and the
reachability problem refers to deciding whether a marking $M$ is
reachable from $M_0$. 
Note: For the type of Petri nets defined so far reachability is decidable 
\cite{Murata} (in exponential time and space). \\

A Petri net $\underline{N}$ is called
\textbf{reversible} if a marking $M_2$ is reachable from another marking $M_1$
implies that $M_1$ 
is reachable from $M_2$. A Petri net with initial marking may be called 
reversible if there is always a  
firing sequence of events that will transform the net from any reachable 
marking back to the initial marking. 

\begin{prop}
Let $\underline{N}$ be a reversible Petri net.
Define $F:=\{pol(t) : t \in T\}$ and let $\lan F \ran$ be the ideal generated 
by $F$ in $\mathbb{Z}[P]$.
Let $M$ and $M'$ be two markings of $\underline{N}$.
Then $M'$ is reachable from $M$ only if
$pol(M)-pol(M') \in \lan F \ran$.
\end{prop}
\begin{proof}
From the definitions above, if $M'$ is reachable from $M$ then there is a 
firing sequence
$M=M_0 \stackrel{t_1}{\to} M_1 \stackrel{t_2}{\to} \cdots 
\stackrel{t_n}{\to} M_n=M'$
so that $pol(M')=pol(M)-pol(t_1)- \cdots - pol(t_n)$. This implies that 
$pol(M)-pol(M') =pol(t_1) + \cdots + pol(t_n) \in \lan F \ran$.
\end{proof}

\begin{example}
\emph{Let $\underline{N}$ be the reversible Petri net given by the 
marked graph below:} 
\vspace{-1.2cm}
{\Large
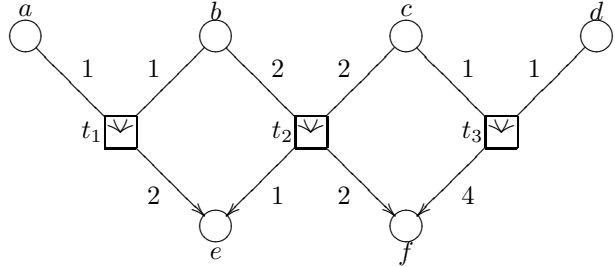
$$ 
\xymatrix{&&&&&&\\ 
            *++[o][F-]{} \ar@{}[u]|(.3){a} \ar[dr]^1 
            && *++[o][F-]{} \ar@{}[u]|(.3){b} \ar[dl]_1 \ar[dr]^2 
            && *++[o][F-]{} \ar@{}[u]|(.3){c} \ar[dl]_2 \ar[dr]^1 
            && *++[o][F-]{} \ar@{}[u]|(.3){d} \ar[dl]_1  \\ 
          & *++[F-]{}    \ar@{}[l]|(.3){t_1} \ar[dr]_2 
            && *++[F-]{} \ar@{}[l]|(.3){t_2} \ar[dl]^1 \ar[dr]_2 
            && *++[F-]{} \ar@{}[l]|(.3){t_3} \ar[dl]^4 &  \\ 
         && *++[o][F-]{} \ar@{}[d]|(.3){e} 
            && *++[o][F-]{} \ar@{}[d]|(.3){f} && \\ 
            &&&&&&\\ } 
$$} 
\vspace{-1cm}

\emph{The places are $P:=\{a,b,c,d,e,f\}$ and the polynomials defined by the 
transitions are 
$t_1 := a+b-2e$, 
$t_2 := 2b+2c-e-2f$ and 
$t_3 := c+d-4f$.
A Gr\"obner basis (using the order $f>e>d>c>b>a$) for the ideal generated in 
$\mathbb{Q}[P]$ is
$$F:= \{d-3c-3b+a, \ e-\frac{1}{2}a-\frac{1}{2}b, \ 
f+\frac{1}{4}a-\frac{3}{4}b-c\}.$$ 
For any marking $M$ the polynomial $pol(M)$ may be reduced, using the relation
$\to_F$ defined by the Gr\"obner basis, to an irreducible form 
$irr(M) \in \mathbb{Q}^{\geqslant 0}[\{a,b,c\}^*]$. Here are three examples.}
\begin{align*} 
pol(M_0) =& ~ 2a+2b+3c+d \to_F  2a+2b+3c-(-3c-3b+a) = a+5b+6c\\ 
pol(M_1) =& ~ 4e+2c+4f    \to_F 
4(\frac{1}{2}a + \frac{1}{2}b) + 2c + 4(-\frac{1}{4}a+\frac{3}{4}b+c) = 
  a+5b+6c \\
pol(M_2) =& ~ a+d+3e+5f \to_F 
a + (3c+3b-a) + 3(\frac{1}{2}a + \frac{1}{2}b) 
  + 5(-\frac{1}{4}a + \frac{3}{4}b +c) =\frac{1}{4}a+\frac{33}{4}b+8c
\end{align*} 
\emph{So $M_2$ is not reachable from $M_0$ because the corresponding polynomials 
do not reduce to the same form. It is here the case that $M_1$ is reachable from 
$M_0$ but this result does not necessarily follow from the
reduced polynomials for these markings being the same.}
\end{example}

\begin{rem}
\emph{We can draw a rational-valued Petri net that is equivalent to the original 
net $\underline{N}$ but 
whose transition polynomials are the Gr\"obner basis and whose markings are a
function $P \to \mathbb{Q}^{\geqslant 0}$. This is constructed 
by drawing a state for each letter and a transition for each polynomial. The 
arcs into a transition come from the letters with positive coefficient and are 
weighted with that coefficient. Similarly the arcs leaving a transition 
correspond to the negative terms in the polynomial.}
\end{rem}

\subsection{Gr\"obner Machines as Petri-Nets} 
 
The Gr\"obner machine for reducing polynomials which was described at the 
end of Section 4.2 can be expressed quite nicely as a Petri net. 
 
\begin{thm} 
Let $K$ be a field, let $X$ be a set and let $F \subseteq K[X^\dagger]$ be a Gr\"obner basis for the ideal $\lan F \ran$. 
Then there is a Petri net $\underline{N}$ which can be marked with a polynomial $f \in K[X^\dagger]$ so that any resulting sequence of firings can be extended to a finite sequence of firings that terminates with a unique non-live state. All states reachable from the initial marking may be identified with polynomials that are equivalent under $=_F$ to $f$. 
\end{thm} 
\begin{proof} 
We will define a type of Petri net and firing rule from the Gr\"obner basis. 
Let $\underline{N}:=(P,T,\cF,w)$. 
The set of places $P$ is the set of 
monomials $m$ of $K[X^\dagger]$ which are irreducible with respect to $\to_F$,
together with an `initial' place labelled $\id$. 
The set of transitions $T$ is identified with $P \times X$.\\

The flow relation $\cF$ is described as follows.  
The transition $(m,x)$ has a single input edge from $m$ with weight $x$.
If $mx \in P$ then $(m,x)$ has a single output edge to $mx$ with weight $1$.
If $mx \not\in P$ then $mx$ is the leading monomial of some 
$f=mx - \Sigma_{i=1}^{n}k_i m_i$ in $F$. 
In this case there is an output edge from $(m,x)$ to each non-leading term in 
$f$, the edge to $m_i$ having weight $k_i$.\\

The Petri net just defined differs from the standard type in that the weight function returns elements of $K$ or elements of $X$ rather than just integers. So $w: \cF \to K[X^\dagger]$.  
Similarly a marking is a function $M: P \to K[X^\dagger]$ 
and is identified with the polynomial  $pol(M):= \Sigma_P \; m M(m)$\\ 

Let $M_1$ be a marking, with $M_1(m) \in  K[X^\dagger]$ for each $m \in P$.
Let $(m,x)$ be an enabled transition, so that $M_1(m)$ contains a term $kxv$ 
for some $k \in K$, $v \in X^*$.  
If $mx$ is irreducible, then when $(m,x)$ fires, the term $kxv$ is removed from
$m$ while $mx$ gains a term $kv$, so the resulting marking $M_2$ is
such that $$pol(M_2):= \Sigma_P m M_2(m) =
\Sigma_P m M_1(m) - m(kxv)+mx(kv)= pol(M_1).$$
Alternatively, when $f=mx - \Sigma_{i=1}^{n} k_i m_i \in F$ and $(m,x)$ fires,
$M_2$ is such that 
$$pol(M_2)= pol(M_1) - m(kxv)+ \Sigma_{i=1}^{n} m_i (k k_i v) = pol(M_1)-kfv,$$
and so $pol(M_1) \to_F pol(M_2)$.\\

Thus a firing represents a single step reduction by $\to_F$. 
The relation is complete, since $F$ is a Gr\"obner basis, and therefore 
there exists a unique non-live marking (irreducible polynomial) which may 
be reached within a finite firing sequence (sequence of reductions). 
\end{proof} 
 
\begin{example}
\emph{The picture for the third Hecke Algebra Petri net 
(whose Grobner machine was Example 4.2.6) is as follows (with each transition label 
$(m,x)$ written $mx$):} 
\vspace{-1cm}

$$ 
\xymatrix{ &&&&&&\\
           &*++[F-]{} \ar@{}[l]|{e_1^2}                \ar@<1ex>[d]|{\,1}  
             & *++[F-]{} \ar@{}[u]|{e_1}            \ar[dl]|1 
             & *++[o][F-]{\ \ } \ar@{}[u]|{\text{\large$\id$}}  \ar[l]|{e_1}  \ar[r]|{e_2} 
             & *++[F-]{} \ar@{}[u]|{e_2}            \ar[dr]|1 
             & *++[F-]{} \ar@{}[r]|{e_2^2}            \ar@<1ex>[d]|1 &\\ 
           &*++[o][F-]{\ \ } \ar@{}[l]|{\text{\large$e_1$}}     \ar@<1ex>[u]|{e_1} \ar[d]|{e_2} 
             &&&& *++[o][F-]{\ \ } \ar@{}[r]|{\text{\large$e_2$}} \ar[d]|{e_1} \ar@<1ex>[u]|{e_2} &\\ 
           &*++[F-]{} \ar@{}[l]|{e_1e_2} \ar[d]|{\,1} 
             &&& *++[F-]{} \ar@{}[u]|{e_2e_1e_2} \ar[ulll]|{\frac{2}{9}}  
                                     \ar[ur]|{-\frac{2}{9}}      \ar[ddl]|1 
             & *++[F-]{}  \ar@{}[r]|{e_2e_1}           \ar[d]|{\,1} &\\ 
           &*++[o][F-]{\ \ } \ar@{}[l]|{\text{\large$e_1e_2$}}        \ar@<1ex>[d]|{e_2} \ar[dr]|{e_1}  
             && *++[F-]{}  \ar@{}[u]|{e_1e_2e_1e_2}          \ar[ll]|{\frac{-2}{9}}  
                                     \ar[uull]|{\frac{2}{9}}   \ar@<1ex>[d]|{\,1} 
            && *++[o][F-]{\ \ } \ar@{}[r]|{\text{\large$e_2e_1$}}   \ar[ul]|{e_2} \ar@<1ex>[d]|{e_1} &\\ 
           &*++[F-]{}  \ar@{}[l]|{e_1e_2^2}               \ar@<1ex>[u]|{\,1} 
             & *++[F-]{}  \ar@{}[d]|{e_1e_2e_1}           \ar[r]|1 
             & *++[o][F-]{\ \ } \ar@{}[d]|{\text{\large$e_1e_2e_1$}} \ar@<1ex>[u]|{e_2} \ar@<1ex>[r]|{e_1} 
             & *++[F-]{}  \ar@{}[d]|{e_1e_2e_1^2}           \ar@<1ex>[l]|1 
             & *++[F-]{}   \ar@{}[r]|{e_2e_1^2}          \ar@<1ex>[u]|{\,1} &\\
           &&&&&&\\} 
$$
\vspace{-0.6cm}

\emph{The states of the Petri net are labelled by the irreducible monomials. 
To reduce a polynomial $p$ take the initial marking $M_0$ to be such that
$M_0(\id)=p$ and $M_0(m)=0$ for all other $m \in P$.
A transition is enabled if the input states to it hold terms which are right 
multiples of the weight on their input arcs.  
Firing of a transition transforms the input and all output states 
simultaneously. For example, if in the 
situation illustrated here the state $s$ holds tokens to a value of $e_2v$ for some 
string $v$ then the transition $t$ is enabled (to the value of $v$). }

\vspace{-1.6cm}
{\Large
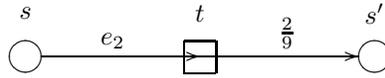
$$ 
\xymatrix{ &&&&\\
    *++[o][F-]{} \ar@{}[u]|{s} \ar[rr]^{e_2} &&  
    *++[F-]{} \ar@{}[u]|{t}    \ar[rr]^{\frac{2}{9}} &&  
    *++[o][F-]{} \ar@{}[u]|{s'}\\} 
$$}
\emph{If transition $t$ then fires, the output state $s'$ receives tokens to the 
value of $\frac{2}{9}v$, which is added to the token value it already holds. 
The marking remaining on the net when all enabled transitions have fired and the  
net is no-longer live (this happens due to the Noetherian property of the  
Gr\"obner basis), represents the irreducible form of the polynomial given by the 
initial marking. This polynomial is extracted from the Petri net by adding the  
token multiples of the states, i.e. if there are 9 tokens at state $e_1$ and  
$\frac{5}{3}$ tokens at state $e_1e_2$ then the polynomial is $9e_1+\frac{5}{3}e_1e_2$.} 
\end{example} 
 
\begin{rem}\emph{ 
The nature of Petri nets is to allow for concurrent  
operations, and this ties in well with the different ways in which a polynomial  
may be reduced by a set of other polynomials.  
A Petri net can be used to model reduction by a set of non-commutative  
polynomials.
It is only in those sets which are Gr\"obner bases, however, that the non-live 
state eventually reached is entirely determined by the initial marking.} 
\end{rem} 
 
\section{Remarks} 
 
The main theme of Chapter Four was the relation between rewrite  
systems / Gr\"obner bases and various types of machine.  \\
 
Automata can be useful for determining whether or not a structure is finite 
(has a finite number of elements). The automaton is drawn directly from the 
complete rewriting system, the equations for it (see \cite{Cohen}) can be 
solved (Arden's theorem) to obtain a regular expression for the language 
(i.e. the set of normal forms of the elements) which will be infinite if 
the free monoid (Kleene star) of some sub-expression occurs. 
Beyond acceptance or rejection of words, these automata have no output. 
It is more helpful to consider the type of machines (``Cayley machines'') 
which take any word as input and output its reduced form.  
We introduced such Cayley machines (or ``Gr\"obner machines'') for algebras. 
Input is a polynomial and the unique irreducible form of that algebra element 
is the output. 
These machines can be seen as types of automata with output or -- as illustrated  
for the polynomial ring case -- as Petri nets.\\ 
 
The main result of the second section was the definition of reduction machines 
for finite Kan extensions.\\ 
The final section of this chapter on machines introduced Petri nets. 
It is of interest to model Gr\"obner bases with Petri nets, because it would be 
extremely useful to find some equivalences between them, so that Petri nets 
could be analysed using Gr\"obner bases. With this aim in mind we showed how the 
``Gr\"obner machine'' for an algebra is a type of Petri net. 
An example of an application of commutative Gr\"obner bases to 
the reachability problem in reversible Petri nets is also given. There is much  
scope for further work in this area.

\chapter{Identities Among Relations}  
  
There is a large number  of papers on computing resolutions of  
groups, in the usual sense of homological algebra. Many of these  
computations are for particular classes of groups (e.g.  
$p$-groups, nilpotent groups) and some of these compute only  
resolutions mod $p$. In general, they do not compute modules of  
identities among relations because they are not specific to a  
presentation.  
  
This problem can be put more generally as that of extending a  
\emph{partial resolution} of a group. That is, we are given an  
exact sequence of free $\bZ G$-modules $C_n \to C_{n-1}\to \cdots  
\to C_1$, and we are asked to extend it by further stages. For the  
identities among relations for a presentation $\mathcal{P}=grp  
\lan X|R \ran$, the initial case is $n=2$ with the boundary given  
by  the Whitehead-Fox derivative 
$$
\partial_2 = 
  ( \partial r / \partial x ) : ( \bZ G )^R \to ( \bZ G )^X.
$$
The problem  
is to extend this by one or two more stages -- the boundaries of  
the free generators of $C_3$ then give generators for the module  
of identities. If also we find $C_4$ and the boundary to $C_3$,  
then we have a module presentation of the module of identities.\\  
  
This problem is usually expressed as `choose generators for the  
kernel of $\partial_2$'. However, it is not clear how this can be  
done algorithmically. The main result of Brown/Razak \cite{BrSa97}
relates this  
problem to the construction of a partial contracting homotopy for  
a partial free crossed resolution of the universal covering  
groupoid of the group $G$. This contracting homotopy is related to  
choices of what are often called 0- and 1-combings of the Cayley  
graph.\\  
  
The main results of this chapter show how to define an  ``extra
information rewriting system''  
or EIRS and how to use this to construct the homotopy $h_1$. The EIRS
records the steps that have been taken in  
rewriting. The `record' is a sequence of elements of the free  
crossed module of the presentation. This shows that the normal  
form function of a complete rewriting system for a group presentation  
determines (up to  
some choices) a set of free generators for the part $C_3$ of a  
resolution, together with the boundary to $C_2$. In fact the  
generators of $C_3$ are in one to one correspondence with the  
elements of $G \times R$, but the boundary depends on the choice  
of complete EIRS.  
This method of computing $h_1$ means that the computation of a set of  
generators for the module of identities among relations is completely  
algorithmic. This work was done with the help of Chris Wensley.  
The computer program $\mathtt{idrels.g}$ implements the procedure.\\  
  
The next problem is that of reducing the generating set of the  
$|R| \times |G|$ identities computed. When the group is small  
(e.g. $S_3$) this can be done by trial and error. In fact $S_3$ is  
a Coxeter group, and for these it has already been proven  
\cite{PrSt90,PrSt88} that the standard presentation yields  
a minimum of  4 generators for the module of identities. The  
methods of these papers do not, however, produce relations among  
these module generators.\\

The example of $S_3$ is used to demonstrate how reduced sets of  
generators at one level determine the identities at the next  
level, and the way in which the reducible elements are expressed  
in terms of the irreducibles allows the calculation of these new  
identities. The example is a good illustration because it is small  
enough to be done by hand, whilst illustrating that the crossed  
resolution for even a small group given by a familiar presentation  
may be quite complex.\\  
  
The final part of the chapter identifies why the problem of  
reducing the set of generators is difficult, and expresses it in  
terms of a Gr\"obner basis problem (the submodule problem).\\  
  
The crossed complex construction of \cite{BrSa97}, together with  
an enhanced rewriting procedure and noncommutative  
Gr\"obner basis theory over rings are brought together to indicate  
an algorithmic method for constructing a free crossed resolution  
of a group. This is an area that will require much further  
development.

\section{Background}  
  
There are strong geometrical and algebraic reasons for studying  
the \emph{module of identities among relations} \cite{BrHu,Pr90}.  
The following exposition gives some of the topological  
background.\\

We assume the usual notion of a presentation $\mathcal{P}  
:= grp\lan X|R \ran$ of a group $G$, where $X$ is a set generating $G$  
and $R \subseteq F(X)$ is called the set of relators. To allow for  
repeated relators we can also consider  presentations  of the form
$grp\lan X,\mathcal{R},w\ran$ where $w:\mathcal{R}\to F(X)$ is a  
function  such that $w(\mathcal{R})=R$. \\  
  
From $\mathcal{P}$ we form the \emph{cell-complex  
$K=K(\mathcal{P})$ of the presentation}. This is a 2-dimensional  
complex.  Its 1-skeleton $K^1$ is $\bigvee_{x \in X}S^1_x$, a  
wedge of directed circles - one for each generator $x \in X$: 
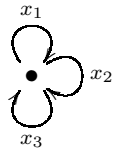
$$  
\xymatrix{\bullet \ar@(ul,ur)^{x_1} \ar@(ur,dr)^{x_2}  
\ar@(dr,dl)^{x_3}} $$ This  topological space has fundamental  
group $\pi_1(K^1,*)$ isomorphic to  the free group $F(X)$ on the  
set $X$. Now $K$ is formed as $$K=K^1 \cup_{\{f_r\}} \{e^2_r \},$$  
by attaching to $K^1$ a 2-cell by a map $f_r:S^1_r \to K^1$ chosen in  
the homotopy class  $w(r) \in F(X)=\pi_1(K^1)$ for each $r\in  
\mathcal{R}$. The homotopy type of $K$ is independent of the  
choice of $f_r$ in its homotopy class.

In the next section we shall define the free crossed module $(\dd_2:
C(w)  
\to F(X))$ on a function $w:\mathcal{R} \to F(X)$.   Whitehead  
\cite{Whitehead1,Whitehead2,Whitehead3} proved that  
$(\pi_2(K^2,K^1,*) \to \pi_1(K^1,*))$ is the free crossed module  
on $w:\mathcal{R} \to \pi_1(K^1,*) \ = \ F(X)$, and so is  
isomorphic to $(C(w) \to F(X))$. In particular $ker \dd_2 \cong  
\pi_2(K,*)$, the second homotopy group of the geometrical model of  
the presentation, and so this homotopy group is also called the  
module of identities among relations for the group presentation.

\begin{example} {\em  
The torus $T=S^1\times S^1$ has a cell structure $(S^1 \vee S^1)  
\cup_{f_r} \{ e^2_r \}$  and its fundamental group  is presented  
by $\mathcal{P} :=  grp \lan a,b \ | \ aba^{-1}b^{-1} \ran$. In this  
case $\pi_2(T)=0$, since $\pi_2 (S^1) = 0$, but it is not so obvious
that $\ker \dd_2 = 0$. }  
\end{example}  
  
More background to these topological ideas may be found in  
\cite{Brown}. There have been many papers written on  
$\pi_2(K^2,*)=ker(C(R) \to F(X))$ (some examples are 
\cite{BaPr92,RB80,BrHi78a,Whitehead1,Whitehead2,Whitehead3,Gru60,Gru70}). 
The methods often use a geometrical notion of ``pictures'' 
\cite{BoPr,Pr90,Pr92,Pr93a,Pr93b,PrSt88} to work with  
identities among relations.  Although the computation of  
$\pi_2(K^2,*)$ is reduced to an algebraic problem on crossed  
modules, this has not previously helped the computation. We shall  
follow the paper \cite{BrSa97} in developing algorithmic methods for  
this computation. For this, we need the language of free crossed  
modules.\\  
  
Let $\mathcal{P}:= grp \lan X | R \ran$ be a group presentation.
An \textbf{identity among relations} is a specified product of  
conjugates of relations  
$$\iota ~ = ~ ({r_1}^{\ep_1})^{u_1}({r_2}^{\ep_2})^{u_2} \cdots  
({r_n}^{\ep_n})^{u_n}$$  
where $ r_i \in R, \ep_i =\pm 1, u_i \in F(X)$  
such that $\iota$ equals the identity in $F(X)$.  
  
\begin{example}
\emph{Let $grp\lan X|R\ran$ be a group presentation.  
Then for any elements $r,s \in R$ we have the identities}  
\begin{center}  
\begin{tabular}{lll}  
$r^{-1} s^{-1} r s^r$ & $=$ & $\id$,\\  
$r s^{-1} r^{-1} s^{r^{-1}}$ & $=$ & $\id$.  
\end{tabular}  
\end{center}  
\end{example}  
When a group has a Cayley graph which forms a simply connected region 
comprised of cells whose boundaries correspond to relators, an identity $\iota$
may be obtained by the following procedure:
\begin{enumerate}[$\bullet$]
\item
Order the cells as $\gamma_1, \ldots, \gamma_m$ in such a way that for all 
$i=1,\ldots,m$ the first $i$ cells form a simply connected sub-region 
$\Lambda_i$.
\item
Choose to transverse each cell in an anti-clockwise direction.
\item
Form a product of of conjugates of relators $v_1 \cdots v_m$ where $v_i$ is 
determined as cell $\gamma_i$ is added to $\Lambda_{i-1}$. To add $\gamma_i$,
start from the vertex $\id$ and move clockwise around the boundary of 
$\Lambda_{i-1}$ until a suitable start vertex on the boundary of $\gamma_i$ is 
reached. A start vertex is such that the word formed by the anti-clockwise 
boundary of $\gamma_i$ starting at that vertex is either the relator $r_i$ 
or the inverse $r_i^{-1}$ of the relator label corresponding $\gamma_i$.
Let $u_i$ be the word given by the path from $\id$ to the start vertex. 
Then the required term is $v_i:=(r_i^{\ep_i})^{u_i^{-1}}$.
\item
Finally set $\iota:=v_1\cdots v_m r_b^{\ep_b}$ where $r_b^{\ep_b}$ is the 
relator associated to the boundary.
\end{enumerate}

\begin{example}
\emph{In the case of a specific group presentation,  
$S_3 = grp\lan x,y \, | \, x^3, y^2, xyxy \ran $, 
label the relators in $S_3$ as $r,s,t$ respectively, and order the 
cells of the Cayley graph as shown below:}  

{\Large
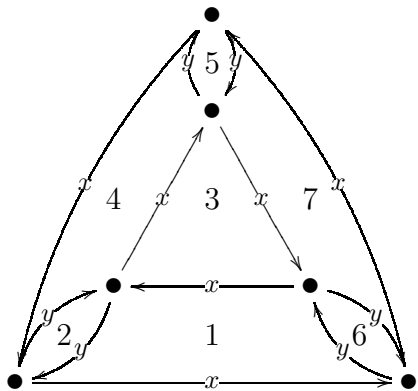
$$  
\xymatrix{&  
          & \bullet \ar@/_1pc/[ddddll]|{x}  
                       \ar@/^0.75pc/[d]|y  
                       \ar@{}[d]|{\text{\large$5$}}
          &  
          &  
          \\  
          & \ar@{}[dd]|{\text{\large$4$}}
          & \bullet \ar[ddr]|x  
                    \ar@/^0.75pc/[u]|y  
                    \ar@{}[dd]|{\text{\large$3$}}
          & \ar@{}[dd]|{\text{\large$7$}}
          &  
          \\  
          &  
          &  
          &  
          &  
          \\  
          &  \bullet  \ar[uur]|x  \ar@/^0.75pc/[dl]|y  \ar@{}[dl]|{\text{\large$2$}}
          &           \ar@{}[d]|{\text{\large$1$}} 
          &  \bullet  \ar[ll]|x   \ar@/^0.75pc/[dr]|y  \ar@{}[dr]|{\text{\large$6$}}
          &  
          \\  
             \bullet  \ar[rrrr]|x \ar@/^0.75pc/[ur]|y  
          &  
          &  
          &  
          &  \bullet  \ar@/_1pc/[uuuull]|x \ar@/^0.75pc/[ul]|y \\}
$$}  
  
\emph{Here cells 1,4,7 (traversed in an anti-clockwise direction) 
correspond to $t$; cells 2,5,6 correspond to $s^{-1}$ 
while cell 3 and the outer boundary (considered as the boundary of the 
``outside cell'') correspond to $r^{-1}$.
We obtain}
\begin{align*} 
\iota ~ :=  & ~  t (s^{-1}) (r^{-1})^{y^{-1}} t^{y^{-1}} (s^{-1})^x  
  (s^{-1})^{x^{-1}y^{-1}x} t^{y^{-1}x} r^{-1}. 
\intertext{\emph{We can verify algebraically that $\iota$ is an 
identity:}}  
\mapsto& ~  
  (xyxy) (y^{-2}) (x^{-3})^{y^{-1}} (xyxy)^{y^{-1}} (y^{-2})^x  
  (y^{-2})^{x^{-1}y^{-1}x} (xyxy)^{y^{-1}x} x^{-3}\\  
=& ~ (xyxy) (y^{-2}) (yx^{-3}y^{-1}) (yxyxyy^{-1}) (x^{-1}y^{-2}x)  
  (x^{-1}yxy^{-2}x^{-1}y^{-1}x) (x^{-1}yxyxyy^{-1}x) (x^{-3})\\  
=& ~ \id.  
\end{align*}  
\end{example}

\section{The Module of Identities Among Relations}  
  
To discuss relations among generators of $G$ we use free groups. To  
discuss  
identities among the relations of $G$  we need free crossed modules.  
The precise idea of a consequence of the relations, and in particular  
of an identity is similar to that of specifying a relator as an  
element of the free group, but takes the action of $F$ into account. \\ 
  
Peiffer and Reidemeister were the first to detail the construction in  
\cite{Peiffer,Reide} in 1949.  
Reidemeister sets up the necessary group action by associating each  
element  
of a first group with an automorphism of a second group, defining a  
homomorphism between the two groups, requiring that it fulfills CM1.  
He looks at the class of Peiffer relations  
of the kernel of this homomorphism, and factors the first  
group by the congruence generated by the Peiffer relations. The  
construction is the same as that detailed below, but he does not mention  
the terms ``group action'' or ``crossed module''. Given that ``crossed  
module'' had only been defined by  
Whitehead in 1946, this is not so surprising.  
It was not until 1982 that perhaps the first paper \cite{BrHu} to  
recognise and name the structures that Reidemeister defined was  
published.\\  
  
Formally, given a group $F$, a \textbf{pre-crossed $F$-module}  
is a pair $(C,\dd)$ where $\dd:C \to F$ is a group  
morphism with an action of $F$ on $C$ denoted  
$c^u$ $(u \in F)$ so that: 
\begin{align*}  
& \text{CM1) } \quad \dd(c^u) ~ = ~ u^{-1}(\dd c) u \quad \text{ for all } \ c
\in  
C, u \in F.\\
\intertext{A \textbf{crossed $F$-module} is a pre-crossed $F$-module
that also 
satisfies  the \textbf{Peiffer relation}:}   
& \text{CM2) } \quad c^{-1}c_1c = c_1^{\dd c} \qquad \text{ for all } \
c,  
c_1 \in C.
\end{align*}  
When $(\dd,C,F)$ is a crossed module it is also
common to refer to it as the crossed $F(X)$-module $(\dd,C)$.
For more information on crossed modules see \cite{BrWe95,BrWe96,BrWe97,LaLu}.\\

The following exposition is a combination of ideas in 
\cite{BrHu,Cremanns,Reide}. It details the construction of the  
module of identities among relations. The construction is not exactly 
the same as that in the references, since it is in terms of rewriting  
systems on a free monoid rather than normal subgroups of a free group.\\

Let \ $\mathcal{P} := grp\lan X, \mathcal{R} ,w\ran $ \ be a  
presentation of a group $G$ where $\mathcal{R}$ is a set of labels for  
the relators identified by the (not necessarily injective function)  
$w:\mathcal{R}\to F(X)$ and $R:=w(\mathcal{R})$.\\
  
A crossed $F(X)$-module $(C,\dd)$ is \textbf{free} on the function 
$w:\mathcal{R} \to F(X)$  
if, given any other crossed $F(X)$-module $(D,\gamma)$ with a map  
$\beta:\mathcal{R} \to D$,  
there exists a unique morphism of crossed $F(X)$-modules $\phi:C \to  
D$ which satisfies $\alpha \circ \phi = \beta$.\\

Define $Y := \mathcal{R} \times F(X)$,  
and write elements of $Y$ in the form $(\rho,u)$, where $\rho \in  
\mathcal{R}, u \in F(X)$.  
  
Put $Y^+:=\{y^+:y \in Y\}$ and $Y^-:=\{y^-:y \in Y\}$. Elements of the  
free monoid  
$(Y^+ \sqcup Y^-)^*$ are called \textbf{Y-sequences} and have the form  
$$(\rho_1,u_1)^{\ep_1} \cdots (\rho_n,u_n)^{\ep_n}.$$  
Define an action of $F(X)$ on $Y$ by  
$$ (\rho,u)^x:=(\rho,ux) \text{ for } x \in F(X).$$  
This induces an action of $F(X)$ on  
$(Y^+ \sqcup Y^-)^*$.  
Define a monoid morphism $\dd: (Y^+ \sqcup Y^-)^* \to F(X)$ to be that
induced by  
$$  
\dd ( \ (\rho,u)^{\ep}) = u^{-1} (w \rho)^{\ep} u  \text{ where } \ep =  
\pm.  
$$  
Define  
\begin{eqnarray*}  
R_P & := & \{(y^-z^+y^+,z^{+\dd y^+}) : y,z \in Y \}\\  
    & \cup & \{(y^+z^-y^-,z^{-\dd y^-}) : y,z \in Y \}\\  
    & \cup & \{(y^-y^+,\id) : y \in Y \}\\  
    & \cup & \{(y^+y^-,\id) : y \in Y \}  
\end{eqnarray*}  
and define $\to_{R_P}$ to be the reduction relation generated by $R_P$  
on $(Y^+ \sqcup Y^-)^*$.  
For $a,b \in (Y^+ \sqcup Y^-)^*$ if $a \stackrel{*}{\lra}_{R_P} b$ then  
$a$ and $b$ are said to be \textbf{Peiffer Equivalent}.\\  
  
\begin{Def}  
The \textbf{Peiffer Problem} is as follows:  
\begin{center}  
\begin{tabular}{lll}  
INPUT: & $a,b \in (Y^+ \sqcup Y^-)^*$ & two elements of the free  
monoid,\\  
QUESTION: & $a \stackrel{*}{\lra}_{R_P} b$? & are they Peiffer  
Equivalent?  
\end{tabular}  
\end{center}  
\end{Def}  
  
The motivation for solving this Peiffer Problem comes from the fact that  
we wish to construct a particular free crossed module, whose kernel will  
be the module of identities among relations.  
Define $$C(R):=  \frac{(Y^+ \sqcup Y^-)^*}{\stackrel{*}{\lra}_{R_P}}.$$  
  
\begin{lem}  
$C(R)$ is a group.  
\end{lem}  
\begin{proof}  
Let $a,b \in (Y^+ \sqcup Y^-)^*$. 
The congruence $\stackrel{*}{\lra}_{R_P}$ preserves the composition of 
Y-sequences so we define $[a]_{R_P}[b]_{R_P}:=[ab]_{R_P}$.  
The identity is $[\id]_{R_P}$, and if $a=y_1^{\ep_1} \cdots y_n^{\ep_n}$  
for $y_1, \ldots, y_n \in Y$, $\ep_1, \ldots, \ep_n= \pm$ then  
$[a]_{R_P}^{-1} := [y_n^{-\ep_n} \cdots y_1^{-\ep_1}]_{R_P}$ is the  
inverse.  
\end{proof}  
  
\begin{lem}  
There is an action of $F(X)$ on $C(R)$ defined by  
$$[a]^x:=[a^x] \text{ for } x \in F(X).$$  
\end{lem}  
\begin{proof}  
Let $y,z \in Y$, $x \in F(X)$  
then $y=(\rho,u)$ and $z=(\sigma,v)$ for some $u,v \in F(X), \rho,  
\sigma \in R.$  
\begin{eqnarray*}  
(y^-z^+y^+)^x & = & (\rho,ux)^- (\sigma,vx)^+ (\rho,ux)^+\\  
              & = & y_1^- z_1^+ y_1^+ \quad \text{ where }  
y_1=(\rho,ux),z_1=(\sigma,vx) \in Y\\  
     &\stackrel{*}{\lra}_P& {z_1}^{+\dd y_1^+}\\  
             & = & (\sigma,vx)^{+\dd(\rho,ux)^+} \quad \text{ by  
definition of } y_1,z_1\\  
             & = & (\sigma,vx(x^{-1}\dd(\rho,u)^+ x)^+ \quad \text{ by  
definition of the action on } (Y^+ \sqcup Y^-)^*\\  
             & = & (\sigma,v \dd(\rho,u)^+ x)^+\\  
             & = & ((\sigma,v)^{+\dd(\rho,u)^+} )^x\\  
             & = & (z^{+\dd y^+})^x \quad \text{ by definition of } y,z  
\end{eqnarray*}  
Similarly $(y^+z^-y^-)^x \stackrel{*}{\lra}_{R_P} (z^{-\dd y^-})^x$, and
it  
is also clear that  
$(y^+y^-)^x \stackrel{*}{\lra}_{R_P} (\id)^x = \id$ and  
$(y^-y^+)^x \stackrel{*}{\lra}_{R_P} (\id)^x = \id$.  
Therefore the action of $F(X)$ on $C(R)$ is well-defined by $[a]^x :=
[a^x]$.  
\end{proof}  
  
\begin{lem}  
There is a group homomorphism $\dd_2:C(R) \to F(X)$ defined by  
$$ \dd_2[a]_{R_P}:=\dd(a) \text{ for } a \in (Y^+ \sqcup Y^-)^*.$$  
\end{lem}  
\begin{proof}  
Let $a,b \in (Y^+ \sqcup Y^-)^*$.  
We require to prove that if $a \stackrel{*}{\lra}_{R_P}  b$ then  
$\dd(a)=\dd(b)$.  
It is therefore sufficient to prove,  
for all $y,z \in Y$, that  
$\dd(y^-z^+y^+)=\dd(z^{+\dd y^+})$,  
$\dd(y^+z^-y^-)=\dd(z^{-\dd y^-})$ and  
$\dd(y^+y^-)=\dd(y^-y^+)=\id_{F(X)}$.  
Let $y=(\rho,u),z=(\sigma,v) \in Y$.  
Then
\begin{align*}  
\dd(y^-z^+y^+) & =  \dd(\rho,u)^- \dd(\sigma,v)^+ \dd(\rho,u)^+,\\  
&= u^{-1} w(\rho)^{-1} u v^{-1} w(\sigma) v u^{-1} w(\rho) u,\\  
&= \dd(\sigma, vu^{-1} w(\rho) u)^+,\\  
&= \dd(\sigma, v \dd(\rho,u)^+)^+,\\  
&= \dd((\sigma,v)^{+\dd(\rho,u)^+}),\\  
&= \dd(z^{+\dd y^+}),
\intertext{and}  
\dd(y^+y^-) &=  \dd(\rho,u)^+\dd(\rho,u)^-,\\  
&= u^{-1} w(\rho) u u^{-1} w(\rho)^{-1} u,\\  
&= \id_{F(X)}. 
\end{align*}
\vspace{-0.5cm}
The other two cases can be proved in the same way,  
therefore $\dd_2$ is well-defined.  
\end{proof}  
  
\begin{thm}  
$(C(R),\dd_2)$ is the free crossed $F(X)$-module on $w:\mathcal{R} \to  
F(X)$.  
\end{thm}  
\begin{proof}  
First we verify the crossed module axioms. \\
CM1:
Let $a=(\rho_1,u_1)^{\ep_1} \cdots (\rho_n,u_n)^{\ep_n}$ for  
$(\rho_1,u_1), \ldots, (\rho_n,u_n) \in Y$, $\ep_1, \ldots, \ep_n=\pm$  
and let $x \in F(X)$.  
Then  
\begin{align*}  
\dd_2([a]_{R_P}^x) &= \dd([(\rho_1,u_1)^{\ep_1}]^x) \cdots  
\dd([(\rho_n ,u_n)^{\ep_n}]^x)\\ &= x^{-1} u_1^{-1}  
w(\rho_1)^{\ep_1 (1)} u_1 x \cdots  
 x^{-1} u_n^{-1} w(\rho_n)^{\ep_n (1)} u_n x,\\  
&= x^{-1} ( u_1^{-1} w(\rho_1)^{\ep_1 (1)} u_1 \cdots u_n^{-1}  
w(\rho_n)^{\ep_n (1)} u_n) x,\\ &= x^{-1} \dd((\rho_1, u_1)^{\ep_1}  
\cdots (\rho_n, u_n)^{\ep_n}) x,\\ &= x^{-1} \dd_2[(\rho_1,  
u_1)^{\ep_1} \cdots (\rho_n, u_n)^{\ep_n}]_{R_P} x,\\ &= x^{-1}  
\dd_2[a]_{R_P} x.\\
\end{align*}

CM2: Let $y,z \in Y$.  We first use the basic rules of $R_P$ to  
verify that $y^+ z^+ y^- \stackrel{*}{\lra}_{R_P} z^{+ \dd y^-}$  
and \linebreak $y^- z^- y^+ \stackrel{*}{\lra}_{R_P} z^{- \dd y^+}$.  
\begin{alignat*}{2}  
z^{+ \dd y^+ y^-} 
&\stackrel{*}{\lra}_{R_P} &\;& (y^+ y^-)^- z^+ (y^+ y^-),\\ 
& =        && y^+ y^- z^+ y^+ y^- \\ 
&\to_{R_P} && y^+ z^{+ \dd y} y^-.
\intertext{Therefore}  
y^+ z^{+ \dd y^+} y^- &\stackrel{*}{\lra}_{R_P}&&(z^{+ \dd  
y^+})^{\dd y_-}. 
\intertext{So for all $z_1 \in Y$} 
y^+ z_1^+ y^- &\stackrel{*}{\lra}_{R_P} &&z_1^{+ \dd y^-}.  
\end{alignat*}  
  
The other case may be proved in the same way but using the basic  
rule $y^+ z^- y^- \to_{R_P} z^{- \dd y^-}$. Therefore the Peiffer  
relation $y^{-\ep} z^\eta y^\ep  \stackrel{*}{\lra}_{R_P} z^{\eta  
\dd y^\ep}$ holds for all $y^\ep, z^\eta \in (Y^+ \sqcup Y^-)^*$. \\ 
  
Let $a=y_1^{\ep_1} \cdots y_n^{\ep_n}$, $b=z_1^{\eta_1} \cdots  
z_m^{\eta_m}$.  
We prove that $[a]_{R_P}^{-1} [b]_{R_P} [a]_{R_P} = [b^{\dd(a)}]_{R_P}$.  
First note that $[a]_{R_P}^{-1}=[y_n^{-\ep_n} \cdots  
y_1^{-\ep_1}]_{R_P}$.  
Now  
\begin{alignat*}{2}  
 y_n^{-\ep_n} \cdots y_1^{-\ep_1} z_1^{\eta_1} \cdots  
z_m^{\eta_m} y_1^{\ep_1} \cdots y_n^{\ep_n} &=&& y_n^{-\ep_n}  
\cdots y_2^{-\ep_2} (y_1^{-\ep_1} z_1^{\eta_1} y_1^{\ep_1} )  
\cdots  ( y_1^{-\ep_1} z_m^{\eta_m} y_1^{\ep_1} ) y_2^{\ep_2}  
\cdots y_n^{\ep_n}, \\ 
& \stackrel{*}{\lra}_{R_P}&\;& y_n^{-\ep_n}  
\cdots y_2^{-\ep_2}  z_1^{\eta_1 \dd y_1^{\ep_1} } \cdots  
z_m^{\eta_m \dd y_1^{\ep_1} } y_2^{\ep_2} \cdots y_n^{\ep_n}.  
\intertext{Repeating the procedure we obtain}  
&\stackrel{*}{\lra}_{R_P} && z_1^{\eta_1 \dd y_1^{\ep_1} \cdots \dd  
y_n^{\ep_n} } \cdots  z_m^{\eta_m \dd y_1^{\ep_1} \cdots \dd  
y_n^{\ep_n} },\\
&= &&(z_1^{\eta_1} \cdots z_m^{\eta_m} )^{\dd(  
y_1^{\ep_1} \cdots y_n^{\ep_n} ) }.
\end{alignat*}
Therefore we have verified CM2:-
$$ 
[y_1^{\ep_1} \cdots y_n^{\ep_n} ]_{R_P}^{-1}  
[z_1^{\eta_1} \cdots z_m^{\eta_m}]_{R_P} [y_1^{\ep_1} \cdots  
y_n^{\ep_n} ]_{R_P} = [(z_1^{\eta_1} \cdots z_m^{\eta_m})^{ \dd  
( y _1^{\ep_1} \cdots y_n^{\ep_n} )}]_{R_P}.  
$$

Finally we show that $(C(R),\dd_2)$ is free on $w:\mathcal{R} \to F(X)$.  
Recall that $F(X)$ acts on $Y$ by $(\rho,u)^x=(\rho,ux)$.  
Define $\alpha:\mathcal{R} \to C(R)$ by $\alpha(\rho):=[(\rho,\id)]_{R_P}$.  
Then let $(D,\gamma)$ be any other crossed $F(X)$-module with a map  
$\beta:\mathcal{R} \to D$.  
We can define a unique morphism of crossed modules $\phi:C(R) \to  
D$ which satisfies $\alpha \circ \phi = \beta$ by  putting 
$\phi([(\rho,u)]_{R_P}):=\beta(\rho)$.\\  
  
Therefore we have proved that $(C(R),\dd_2)$, as defined on $(Y^+ \sqcup  
Y^-)^*$ using ${R_P}$, is the free crossed $F(X)$-module generated by  
$w:\mathcal{R} \to F(X)$.  
\end{proof}  
  
\begin{rem}\emph{  
The usual method of construction of $C(R)$ does not use rewriting  
systems but factors the free precrossed module $(F(Y),\dd')$ by  
the congruence $=_P$ generated by the set of all Peiffer relations  
$P$ on $F(Y)$. Detail of this construction are found in  
\cite{BrHu}. It may be  
verified that the natural map 
$\theta:(Y^+ \sqcup Y^-)^* \to F(Y)$ induces an isomorphism 
$$\theta': \frac{(Y^+ \sqcup  
Y^-)^*}{\stackrel{*}{\lra}_{R_P}} \longrightarrow \frac{F(Y)}{=_P}.$$  
The motivation for this section is to give an exposition of the
construction 
of $C(R)$. Since this thesis is concerned with rewriting, we've
presented the 
exposition in terms of rewriting. 
It is simply an alternative exposition of standard work that is
necessary 
background for what is to follow.} 
\end{rem}  
  
The Peiffer Problem that we have identified is that of determining  
whether two Y-sequences represent the same element of $C(R)$. 
If $a \in (Y^+ \sqcup Y^-)^*$ and  
$\dd_2(a) = \id$ then $[a]_{R_P} \in ker \dd_2$, the module of
identities  
among relations, and $a$ is called an \textbf{identity Y-sequence}.  
There is a special property which will allow us to convert the Peiffer  
Problem for identity sequences into a Gr\"obner basis problem, and this  
will be discussed in Section 6.  
In general there is no procedure for solving the Peiffer Problem. As a  
result the example here is a simple one, included to demonstrate the rewriting  
procedure.  
  
\begin{example}\emph{
The result of the following example is proved in \cite{BrWe95}.}\\

\emph{The multiplicative cyclic group $\bC_n$ of order $n$ has a
presentation  
 $grp\lan x \ | \ x^n\ran $.
Let $r$ represent the relator $x^n$, then  
\ $Y:=\{(r,x^i): i \in \bZ \}$.  
with $\dd:(Y^+ \sqcup Y^-)^* \to F(X)$ defined by $\dd(r) = x^n$ so } 
\begin{eqnarray*}  
\dd_2(r,x^i)^+ &=& x^{-i} \dd(r) x^i ~ = ~ x^{-i} (x^n) x^i ~ = ~ x^n.\\  
\dd_2(r,x^i)^- &=& x^{-i} \dd(r)^{-1} x^i ~ = ~ x^{-i} (x^n) x^i  ~ = ~ x^{-n}.  
\end{eqnarray*}  
\emph{
The action of $F(X)$ on $(Y^+ \sqcup Y^-)^*$ is given by }  
$$(r,x^i)^x = (r,x^{i+1}).$$  

\emph{
The elements of $Y^+ \sqcup Y^-$ can be denoted $a_i, A_i \ i \in \bZ$  
where $a_i:=(r,x^i)^+$, $A_i:=(r,x_i)^-$.  
We consider the rewriting system $R_P$ on $(Y^+ \sqcup Y^-)^*$ given by:}  
$$\{(A_i a_j a_i, a_j^{\dd a_i}) : i, j \in \bZ \}   
\, \cup \, \{(a_i A_j A_i, A_j^{\dd A_i}) : i, j \in \bZ \}   
\, \cup \, \{(a_i A_i, \id) : i \in \bZ \} 
\, \cup \, \{(A_i a_i, \id) : i \in \bZ \}
$$
\emph{The rewriting system is clearly infinite.  
Put $i=j$ in the above rules and we obtain  
$A_i a_i a_i \lra_{R_P} a_{i+n}$ and $a_iA_iA_i \lra_{R_P} A_{i-n}$.  
So $a_{i+n} \to_{R_P} a_i$ and $A_{i} \to_{R_P} A_{i-n}$ for all $i \in  
\bZ$.
It follows immediately from these rules that $\{a_0, \ldots ,a_{n-1},  
A_0, \ldots, A_{n-1} \}$ is a complete set of  
generators for $C(R)$ as a monoid. The now finite set of relations is  
$\{ (a_i A_i , \id), (A_i a_i, \id), (A_i a_j a_i, a_j), (a_i A_j A_i,  
A_j) \}$  
Therefore $C(R)$ for $\bC_n$ is the free abelian  
group on $n$ generators $a_0, a_1, \ldots, a_{n-1}$. 
Further, we find that $a_i^x=a_{i+1}$ for $i=0, \ldots, n-1$ and
$a^x_{n-1}=a_0$. Thus the $C(R)$, which is a $\bC_n$-module is
isomorphic to $\bZ[\bC_n]$, the free $\bC_n$-module on one generator.}  
\end{example}  
  
\begin{rem}  
\emph{  
The Peiffer Problem (of deciding when two sequences are Peiffer
equivalent) does not arise only in crossed modules. When a 2-category is  
constructed, by factoring a sesquicategory (see \cite{Stell,Street})  
by the interchange law, the pairs arising from that interchange law are  
relations among the two cells involving the whiskering action of the  
category  
morphisms. Tim Porter identified this in \cite{Tim} calling them Peiffer  
pairs.  
Thus the Peiffer Problem is not restricted to the construction of  
crossed modules.}  
\end{rem}

\section{Free Crossed Resolutions of Groups}  
  
The following exposition was constructed with Ronnie Brown.\\  
  
The notion of resolution of $\bZ G$-modules for $G$ a group is a  
standard part of homological algebra and the cohomology of groups  
\cite{Cartan,KBrown}.  
It has been shown in \cite{BrWe95,BrPo,BrSa97} 
that there are computational advantages in considering free  
\emph{crossed} resolutions of groups. This will be confirmed by bringing
these  
calculations into the context of rewriting procedures. For this we need  
to give some basic definitions in the form we require.\\  
  
An important aspect of the calculation in \cite{BrSa97} is the use  
of the Cayley graph, being seen here as data for a free crossed  
resolution of the universal covering groupoid $\Gt$ of the group  
$G$. This groupoid corresponds to the action of $G$ on itself by  
right multiplication. That is, the objects of $\Gt$ are the  
elements of $G$ and an arrow of $\Gt$ is a pair $(g_1,g_2):g_1 \to g_1
g_2$,  
with the obvious composition. We have the covering morphism of groupoids
$p_0:  
\Gt \to G:(g_1,g_2) \mapsto g_2$. \\  
  
If $X$ is a set of generators of the group $G$, we have a standard  
morphism $\theta:F(X) \to G$. We also have a standard morphism  
$\widetilde{\theta}: F(\Xt) \to \Gt$. Here  
\begin{enumerate}[i)]  
\item  
$\Xt$ is the Cayley graph of $(X,G)$ with arrows $[g,x]:g \to
g\theta(x)$ for $x \in X, g \in G$.  
\item  
$F(\Xt)$ is the groupoid with objects again the elements of $G$ and
arrows pairs  
$[g,u]:g \to g( \theta u )$ for $g \in G$, $u \in F(X)$,
with composition defined by $[g, u][g(\theta u), v] := [g, uv]$.  In
fact $F(\Xt)$ is the free groupoid on the graph $\Xt$, so that a
morphism $f$ from $F(\Xt)$ to a groupoid is determined by the graph
morphism $f|_{\Xt}$.
\end{enumerate}  
Then $\widetilde{\theta} : F(\Xt) \to \Gt$ is given on arrows by  
$\widetilde{\theta}[g,u] := [g, \theta(u)]$. There is also the  
covering morphism $p_1: F(\Xt) \to F(X)$ given by $p_1 [g,u]:=u$.  
This gives the commutative diagram of morphisms of groupoids
\vspace{-.3cm}  
\begin{equation}
\label{one} \xymatrix{ F(\Xt) \ar[r]^-{\widetilde{\theta}}  
\ar[d]_-{p_1}  
           & \Gt \ar[d]^-{p_0} \\ F(X) \ar[r]_-\theta & G\\}  
\end{equation}  
In fact this diagram is a pullback in the category of  
groupoids. Also, $p_1$ maps $F(\Xt)(1,1)$ isomorphically to $ker
\theta$, and $F(\Xt)$ is the free groupoid on the graph $\Xt$.  
  
Now let $\mathcal{P}= grp \lan X | R \ran$ be a presentation of  
$G$. As explained in the previous section, this gives rise to a  
free crossed $F(X)$-module $\dd_2:C(R) \to F(X)$, whose kernel is  
$\pi_2(\mathcal{P})$, the $\bZ G$-module of identities among  
relations. The aim is to compute a presentation for this module in  
terms of information on the Cayley graph. For this we extend  
diagram \ref{one} in the first instance to 
\begin{equation}
\label{two}  
\xymatrix{ C(\Rt) \ar[d]_-{p_2} \ar[r]^-{\db_2} &  
           F(\Xt) \ar[r]^-{\widetilde{\theta}}
\ar[d]_-{p_1} &   
           \Gt \ar[d]^-{p_0} \\  
         C(R) \ar[r]^-{\dd_2} & 
         F(X) \ar[r]^-\theta  & G\\}  
\end{equation}  
Here $\db_2: C(\Rt) \to F(\Xt)$ is a free crossed module of groupoids.  
For details, we refer the reader to \cite{BrSa97}.  
All the reader needs to know for now is that  
\begin{enumerate}[i)]  
\item  
$C(\Rt)$ is a disjoint union of groups $C(\Rt)(g)$ for $g \in G$  
and $\db_2$ maps $C(\Rt)(g)$ to $F(\Xt)(g,g)$.  
\item for each $g \in G, \;  
p_2$ maps the group $C(\Rt)(g)$ isomorphically to $C(R)$, so that  
elements of $C(\Rt)(g)$ are specified by pairs $[g,c]$ where $c  
\in  C(R)$.  
\item  
$F(\Xt)$ operates on  
$C(\Rt)$ by $[g, c]^{[g,u]} := [g \theta(u), c^u]$  
for $g \in G$, $c \in C(R)$, $u \in F(X)$.  
\item  
The morphisms $\db_2$, $p_2$ are given by  
$\db_2[g,c] := [g, \dd_2 c]$ and $p_2[g, c]:= c$.  
\end{enumerate}  
A proof that $\db_2: C(\Rt) \to F(\Xt)$ is the free crossed  
$F(\Xt)$-module on $\Rt:=G \times R$ is given in \cite{BrSa97}. 
This implies that morphisms and homotopies on $C(\Rt)$ can be defined by
their values on the elements $[g,r]$ 
for $g \in G$, $r \in R$.  \\
  
The key feature of this construction is that $\Gt$ is a  
contractible groupoid, i.e. it is connected and has trivial vertex  
groups. We are going to construct a partial contracting homotopy  
of $\db_2: C(\Rt) \to F(\Xt)$. This is a key part of the procedure  
of constructing generators (and then relations) for  
$\pi_2(\mathcal{P})$. The philosophy as stated in \cite{BrSa97} is  
to construct a ``home'' for a contracting homotopy -- this will be  
explained later. The point is that this leads to a  
``tautological'' proof that the generators constructed do in fact  
generate $\pi_2(\mathcal{P})$.\\  
  
Such a partial contracting homotopy consists of functions  
$$ 
h_0  : G  \to F(\Xt) \quad \text{ and  } \quad 
h_1  : F(\Xt)  \to C(\Rt)  
$$  
with the properties that  
\begin{enumerate} [i)]  
\item  
$h_0(g):g \mapsto \id$ in $F(\Xt)$, $g \in G$.  
\item  
$h_1$ is a morphism (from a groupoid to a group).  
\item  
$\db_2 h_1 [g,u]= (h_0 g)^{-1} [g, u] h_0(g ( \theta u ))$ for all
$[g,u]  
\in F(\Xt)$.  
 \end{enumerate}  
We always assume that $h_0(\id)= \id \in F(\Xt)(\id)$ 
  
\begin{rem} \emph{  
$h_0$ and $h_1$ are related to what are commonly called $0$- and  
$1$-combings of the Cayley graph \cite{Hermiller3}. We hope to pursue
this  
elsewhere. }  
\end{rem}  
  
The choice of $h_0$ is equivalent to choosing a section $\sigma$ of
$\theta: F(X)  
\to G$, i.e. a representative word for each element of $G$, by $h_0(g) =
[g, \sigma(g)^{-1}]$, for $g \in G$.  
What $h_1$ does is provide for each word $u \in F(X)$ a representation  
$$u=\dd_2(\mathtt{proc}_R(u))N_R(u)$$  
where $\mathtt{proc}(u)=p_2 h_1 [ \id, u] \in C(R)$ -- the procedure through
which the normal form $N_R(u) := (\sigma \theta (u))^{-1}$ is reached.
To verify this 
consider (iii), assuming $h_0(\id)=\id$, we have  
$$\db_2 h_1 [ \id,u] = [\id, u] h_0 ( \theta u ).$$  
Then  
\begin{eqnarray*}  
\dd_2(proc(u)) &=& \dd_2 p_2 h_1 [\id, u]\\  
               &=& p_1 \db_2 h_1 [\id, u]\\  
               &=& p_1([\id, u] h_0( \theta u ) )\\  
               &=& u p_1 h_0 ( \theta u )
\end{eqnarray*}  
  
Thus $\mathtt{proc}(u)$ shows how to write $u(N_R(u))^{-1} \in \dd_2 C(R)$ as
a consequence of the relators $R$.  
Conversely, a rewriting procedure to be given later will allow us to  
determine $h_1$ given $h_0$ and a complete rewriting system for  
$\mathcal{P}= grp \lan X | R \ran$.\\  
  
We can now state  
  
\begin{prop}  
Given $h_0$, $h_1$ as above, the module $\pi_2( \mathcal{P} )$ is  
generated by the (separation) elements 
\begin{equation} \label{identityformula}
sep(g,r) := p_2 (h_1 \db_2[g, r])^{-1} r^{\sigma(g)^{-1}} 
\end{equation} 
for all $g \in G$, $r \in R$.  
\end{prop}  
  
\begin{out}
The fact that the elements $sep(g,r)$ of \ref{identityformula} are
identities among  
relations is easily checked, as follows:
\begin{align*} 
\dd_2(p_2 (h_1 \db_2[g, r])^{-1} r^{\sigma(g)^{-1}} ) 
& = \dd_2(p_2 (h_1 [\id, \dd_2(r^g)])^{-1} r^{\sigma(g)^{-1}})\\  
& = \dd_2(p_2 ( [\id,c])^{-1} r^{\sigma(g)^{-1}}) 
\text{ where $c$ satisfies } \dd_2(c)=\dd_2(r^{\sigma(g)^{-1}}),\\ 
& = \dd_2(c)^{-1} \dd_2(r^{\sigma(g)^{-1}}) \\
& = \id.
\end{align*}
The important point is that these elements $sep(g,r)$ generate  
the \emph{module} of identities. The proof of this can be made  
tautologous by taking the construction one step further, i.e.  
$$ 
\xymatrix{ \Cb_3 \ar[r]^-{\db_3} \ar[d]_-{p_3}  
& C(\Rt) \ar[d]^-{p_2} \ar[r]^-{\db_2}  
& F(\Xt) \ar[r]^-{\widetilde{\theta}} \ar[d]^-{p_1}  
& \Gt \ar[d]^-{p_0} \\  
  C_3 \ar[r]^-{\dd_3}  
& C(R)\ar[r]^-{\dd_2}  
& F(X) \ar[r]^-\theta  
& G\\}  
$$ 
Here $C_3$ is the free $\bZ G$-module on $(g, r) \in \bar{R}$ where  
$\bar{R}:=G \times R$ -- we use round brackets to distinguish elements of  
$\bar{R}$ from those of $\Rt$. The morphism $\dd_3$ is defined by  
\begin{align*}
\dd_3(g,r) &:= p_2(( h_1 \db_2[g,r])^{-1})r^{\sigma(g)^{-1}}.  
\intertext{The definition is verified by checking that 
$\dd_2 \dd_3 (g,r)=\id$ \ i.e. } 
\dd_2 \dd_3(g,r)  
& = \dd_2 p_2((h_1 \db_2[g,r])^{-1} r^{\sigma(g)^{-1}}) \\
& = \dd_2(c^{-1} r^{\sigma(g)^{-1}}) 
\text{ where $c$ satisfies $\dd_2(c)=\dd_2(r^{\sigma(g)^{-1}})$}\\ 
& = \id.
\intertext{(Mapping a free $\bZ G$-module into a free crossed $G$-module, 
is acceptable because the image lies in $ker \dd_2$ which is a $\bZ G$-module.)  
In fact we define $\Cb_3$, $h_2$ and $\db_3$ as follows}  
\Cb_3(g)           & :=  \{ g \} \times C_3, \\  
h_2[ g, r ]        & :=  (\id,(g,r)),\\  
\db_3(g_2,[g_1,r]) & :=  (g_2, \dd_3 (g_1,r) ).  
\intertext{We now check directly that}  
\db_3 h_2[g,r] & = [ \id, \dd_3(g,r)],\\  
               & = [ \id, p_2(( h_1 \db_2[g,r])^{-1}) r^{\sigma(g)^{-1}} ],  
\intertext{so that}  
               & = (h_1( \dd_1 [g,r] ))^{-1} r^{\sigma(g)^{-1}}.\\ 
\intertext{In the partial resolution of $\Gt$ we have, for any $c \in C(\Rt)$,}
\db h_2(c)     & = (h_1(\db_2 c))^{-1} c^{h_0 \id},
\end{align*}
since this holds for all $c=[g,r] \in \Rt$. So 
$$\db_2(c)=0 \text{ implies that } c= \db_3( (h_2 c)^{(h_0 \id)^{-1}}.$$
Hence $ker \db_2  \subseteq im \db_3$, so $ker \db_2  = im \db_3$. 
Therefore $ker \dd_2  = im \dd_3$. 
\end{out}

To summarise: the problem of constructing a crossed resolution of a  
group given a particular presentation has been reduced to the problem of  
constructing a contracting homotopy and a covering crossed complex that  
begins with a groupoid defined on the Cayley graph.

\section{Completion Procedure and Contracting Homotopies}  
  
In this section we define what we call an ``extra information completion  
procedure''. The implementation may be found in $kb2.g$.  
Input to the procedure is a set of relators for a group.  
If the procedure terminates then the output is a set of ``extra  
information'' rules. These rules will not only reduce any word in the  
free group to a unique irreducible but will express the actual reduction  
in terms of the original relators.  
  
\begin{Def}  
An \textbf{extra information rewriting system} for a group presentation
$grp\lan X|R \ran$ is a set of triples $R2 := \{ (l_1, c_1, r_1),  
\ldots,(l_n,c_n,r_n) \}$,  
where  $R1 := \{ (l_1,r_1), \ldots, (l_n,r_n)\}$ 
is a rewriting system on $F(X)$ and $c_1, \ldots, c_n \in C(R)$, 
such that $l_i = \dd_2(c_i) r_i$ for $i=1,\ldots,n$. 
We say $R2$ is \textbf{complete} if $R1$ is complete.  
\end{Def}  
  
\begin{lem}  
Let $R2$ be a complete EIRS for $grp\lan  
X|R \ran$.  
Then for any $w \in F(X)$ there exists $(c,z)$, $c \in C(R)$, $z \in
F(X)$ such that $z$ is  
irreducible with respect to $\to_{R1}$, and $w= (\dd_2 c) z$.  
\end{lem}
  
\begin{proof}  
If $w$ is irreducible then we take $z=w$ and $c=\id_{C(R)}$.  
Otherwise there is a sequence of reductions
\vspace{-.5cm}
\begin{align*}  
w & =  u_1l_1v_1\\  
u_1r_1v_1 & =  u_2l_2v_2\\  
\cdots &  \qquad \cdots\\  
u_nr_nv_n & =  z  
\intertext{where $n \geq 1$, and for $i=1, \ldots, n$, $u_i,v_i \in F(X)$ and there
exists $c_i \in C(R)$ such that $(l_i,c_i,r_i) \in R2 $.  
Then since $l_i= (\dd_2 c_i)\, r_i$ for $i=1, \ldots, n$}
\vspace{-.5cm}  
w & =  u_1 \, (\dd_2 c_1) \, r_1 v_1\\  
u_1r_1v_1 & =  u_2 \, (\dd_2 c_2) \, r_2 v_2\\  
\cdots &  \qquad \cdots\\  
u_nr_nv_n & =  z.  
\end{align*}  
Hence $w = ( (\dd_2 c_1)^{u_1^{-1}} \cdots (\dd_2 c_n)^{u_n^{-1}} ) z$.  
\end{proof}  
  
This defines the function $\mathtt{ReduceWord2}$, which accepts as input  
$(w,R2)$ and returns as output $(c,z)$. 
We will write $w \to_{R2} (c,z)$.

\begin{lem}  
Let $grp\lan X|R \ran$ be a finite group presentation which is  
completable with respect to an ordering $>$.  
Then there exists a procedure $\mathtt{KB2}$ which will return the  
complete EIRS for the group.  
\end{lem}  

\begin{proof}  
Define $R2:=\{ (\dd \rho, (\rho,\id),\id) : \rho \in R\}$.  
It is clear that this defines an EIRS since  
$\dd \rho = \dd_2(\rho,\id) \id$.  
  
If $R1$ is complete then $R2$ is complete. If $R1$ is not complete then  
there is an overlap between a pair of rules $(l_1,r_1), (l_2,r_2)$ of  
$R1$ where $(l_1,c_1,r_1), (l_2,c_2,r_2) \in R2$.  
There are two cases to consider.  
  
For the first case suppose $u_1l_1v_1=l_2$ for some $u_1, v_1 \in F(X)$.  
Then the critical pair resulting from the overlap is $(u_1r_1v_1,r_2)$.  
Reduce each side of the pair using $\mathtt{ReduceWord2}$, so  
$u_1r_1v_1 \to_{R2} (d_1, z_1)$ and $r_2 \to_{R2} (d_2, z_2)$. Then if  
$z_1>z_2$ add the extra information rule  
$(z_1, d_1^{-1} {c_1}^{-u_1^{-1}} c_2 d_2, z_2)$ or if $z_2>z_1$ add  
$(z_2, d_2^{-1} c_2^{-1} c_1^{u_1^{-1}} d_1,z_1)$.  
  
For the second case suppose $u_1l_1=l_2v_2$ for some $u_1, v_2 \in  
F(X)$.  
Then the critical pair resulting from the overlap is $(u_1r_1,r_2v_2)$.  
Reduce each side of the pair by $R2$ as before, so that  
$u_1r_1 \to_{R2} (d_1, z_1)$ and $r_2v_2 \to_{R2} (d_2, z_2)$. Then if  
$z_1>z_2$ add the extra information rule  
$(z_1, d_1^{-1}  c_1^{-u_1^{-1}} c_2 d_2, z_2)$ or if $z_2>z_1$ add  
$(z_2, d_2^{-1}  c_2^{-1} c_1^{u_1^{-1}} d_1,z_1)$.  
  
It can be seen immediately from the above that the effect on $R1$ is a  
standard completion of the rewriting system, and that the triples  
$(l,c,r)$ added to $R2$ satisfy the requirement $l=\dd_2(c) r$, so that  
when the completion procedure terminates $R2$ will be a complete extra  
information rewriting system.
\end{proof}
  
This defines the procedure $\mathtt{KB2}$.

\begin{example}\emph{
$Q8$ is presented by $grp\lan a,b \ | \ a^4, b^4, abab^{-1}, a^2b^2\ran $.  
Let $r,s,t$ and $u$ denote the relators i.e. $\dd(r)=a^4,\dd(s)=b^3, \ldots$.  
We begin with the EIRS 
$$R2:=\{(a^4,r,\id), \ (b^4,s,\id), \ (aba,t,b), \ (a^2b^2,u,\id)\}.$$  
As explained before, all the extra information rules are triples $(l,c,r)$
such that $l= (\dd_2 c) r$ and we write $l \to_{R2} (c, r)$, 
thinking of the $(c)$ part as the record of the procedure by which $r$ is 
obtained from $l$ using the original group relators. For example 
$aba \to_{R2} (t, b)$ -- we have to work with a monoid presentation and choose 
to make use of the fact that $Q8$ is finite, rather than introduce generators 
for the inverses, which is what the computer program does.
We look for overlaps between the left hand sides of the
rules. The first overlap we examine is between the first and third rules:  
$$  
\xymatrix{& a^4ba \ar[dl]|{a^4 \to \id} \ar[dr]|{aba \to b} & \\  
          ba && a^3b \ar@{-->}[ll] \\}  
$$  
Without the extra information the critical pair is $(a^3b,ba)$ and the
new rule is $a^3b \to ba$. 
For the EIRS rule we need $c$ so that  $a^3b= \dd_c(c)ba$  
where $c$ is a product of conjugates of relators. 
The new EIRS rule as defined in the proof (second case)
is $(a^3b,t^{-a^{-3}}r,ba)$.  This is checked by:}
$$
a^4ba=(a^4)ba \to_{R2} (r,\id)ba = (r,ba) \text{ and }
a^4ba=a^3(aba) \to_{R2} a^3(t,b) = (t^{a^{-3}},a^3b).
$$  
\emph{Therefore $\dd_2(r) ba = \dd_2(t^{a^{-3}}) a^3b$,  
so $a^3b = \dd_2(t^{a^{-3}})^{-1}(r) ba = \dd_2(t^{-a^{-3}}r) ba.$  
so $c=t^{-a^{-3}}r$.
If we continue this ``extra information completion'' for $Q8$ we end up  
with the EIRS}  
\begin{align*}    
b^2 & \to_{R2} (r^{-1}u^{a^{-2}}, a^2),\\  
aba & \to_{R2} (t, b),\\  
ba^2 & \to_{R2} (t^{-1}t^{-a^{-1}}r^{b^{-1}a^{-2}}, a^2b),\\  
bab & \to_{R2} (r^{-b^{-1}a^{-1}b^{-1}} t^{b^{-1}} r^{-1} u^{-a^{-2}} r, a),\\  
a^4 & \to_{R2} (r, \id),\\  
a^3b & \to_{R2} (t^{-a^{-3}} r, ba).  
\end{align*}  
\emph{ So, for example, $a^5ba^3$ reduces to $a^2b$ and  
$a^5ba^3 = \dd_2(rt^{-a^{-1}}r^{b^{-1}a^{-2}})a^2b$.}
\end{example}

The ``extra information'' Knuth-Bendix procedure $\mathtt{KB2}$ 
results in a rewriting system  
with information on where the rules came from. This extra information  
is in no way unique.\\
  
Let $grp\lan X|R\ran$ be a presentation of a group $G$.  
Let $\Xt$ denote the Cayley graph. Edges of the graph are recorded as  
pairs  
$[g,x]$, where $g$ is the group element identified with the source  
vertex, and $x$ is a group generator identified with the edge label.
  
\begin{lem}[Complete Rewriting Systems Determine $h_0$]\mbox{ }\\ 
Let $G$ be a finite group, finitely presented by $grp\lan X|R\ran$, with 
quotient morphism $\theta:F(X) \to G$.  
Then a complete rewriting system $R1$ for the presentation determines  
$h_0:G \to F(\Xt)$.  
\end{lem}  
\begin{proof}  
Let $N$ be the normal form function defined by $\to_{R1}$ on $F(X)$.  
Define $h_0(g) := [\id, N(g)]^{-1}$.  
Then $h_0(g):g \to \id$ in $F(\Xt)$ as required.  
\end{proof}  
  
\begin{thm}[Complete EIRS's Determine $h_1$]\mbox{ }\\
Let $G$ be a finite group, finitely presented by $grp\lan X|R\ran$, with
quotient morphism $\theta:F(X) \to G$.  
Then a complete EIRS $R2$ for the  
presentation determines $h_1: F(\Xt) \to C(\Rt)(\id)$.  
\end{thm}  
\begin{proof}  
Recall that $\Xt$ is the Cayley graph of $G$. 
Let $[g,x] \in \Xt$.  
Define  
$$h_1[g,x] := [\id, \mathtt{ReduceWord2}(N(g)xN(g\theta x)^{-1},R2)[1]].$$  
Then clearly $h_1[g,x] \in C(\Rt)(\id)$ and  
$\db_2 h_1[g,x]=[\id,N(g)] [g,x] [\id,N(g\theta x)]^{-1} 
               = h_0(g)^{-1} [g,x] h_0(g \theta x) $.
Extending this definition of $h_1$ on $\Xt$ therefore gives the
morphism $h_1$ of the groupoid $F(\Xt)$ to the group $C(\Rt)(\id)$ satisfying
the required conditions.  
\end{proof}  
  
\begin{cor}  
There exists an algorithm for defining $h_0,h_1$ for any finite  
completable group presentation $grp\lan X|R \ran$.  
\end{cor}  
\begin{proof}  
Calculate $R2$, using $\mathtt{KB2}$. Let $N$ be the normal form  
function defined by $\to_{R1}$ (recall $R1$ is part of $R2$).  
Put $h_0(g):=[\id,N(g)]^{-1}$.  
Put $h_1[g,x]:=[\id,\mathtt{ReduceWord2}(N(g) x N(g \theta x)^{-1}, R2)[1]]$.  
\end{proof}  
  
\begin{example}
\emph{Below is the Cayley graph for $Q8$. The double edges indicate the tree  
defined by the length-lex ordering.}  
$$  
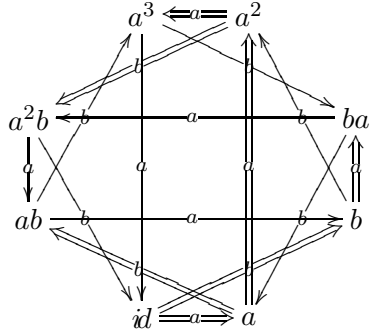
\xymatrix{& a^3 \ar[ddd]|a \ar[drr]|b  
          & a^2 \ar@{=>}[l]|a \ar@{=>}[dll]|b &\\  
          a^2b  \ar[d]|a \ar[ddr]|b &  
          && ba \ar[lll]|a \ar[ddl]|b \\  
          ab    \ar[rrr]|a \ar[uur]|b &  
          && b  \ar@{=>}[u]|a \ar[uul]|b \\  
          & \id \ar@{=>}[r]|a \ar@{=>}[urr]|b  
          & a   \ar@{=>}[uuu]|a \ar@{=>}[ull]|b &\\}  
$$  
\emph{A typical relator cycle is $[\id,b] [b,a] [ba,a] [a^2b,a] [ab,a]  
[\id,b]^{-1}$  
this is equivalent to $ba^4b^{-1}$ or $r^{b^{-1}}$, the cycles represent  
conjugates of relators in the graph.\\}

\emph{The extra information in our rewriting system may be used to express the  
cycle created by adding an edge $\alpha$ to the tree as such a product, or in  
fact to  
express its retraction as a product of conjugates of relators.  
For example, add the edge $[a^3,a]$ and the cycle  
$[\id,a][a,a][a^2,a][a^3,a]$ is created. The retraction is $a^4$. 
We know that $\dd_2(r)=a^4$, so $[\id,r]$ is the cycle as a product of 
relator cycles.}\\

\emph{That was an easy example. 
If we add the edge $[a^3,b]$ then the retraction is  
$a^3 b a^{-1} b^{-1}$ or $a^3 b (a^3) (a^2b)$, (since the rewriting system is 
defined on the monoid presentation we replace inverse elements by their 
normal forms). It is more difficult to see how  
this word may be written as a product of conjugates of relators. 
In fact we just reduce it using the extra information rules :}  
\begin{align*}  
a^3ba^5b &\to_{R2} a^3ba(r,\id)b \\  
         &\to_{R2} a^3bab(r^b,\id)\\  
         &\to_{R2} a^3((r^{-b^{-1}a^{-1}b^{-1}} t^{b^{-1}} r^{-1}  
                         u^{-a^{-2}} r, a) (r^b,\id)\\  
         &\to_{R2} a^4(r^{ -b^{-1}a^{-1}b^{-1}} t^{b^{-1}} r^{-1}  
                        u^{ -a^{-2} } r)^a, \id) (r^b,\id)\\  
         &\to_{R2} (r (r^{ -b^{-1} a^{-1} b^{-1} } t^{ b^{-1} } r^{-1}  
                           u^{ -a^{-2} } r )^a r^b , \id)
\end{align*}
\emph{The order in which the rules are applied does not matter for our
purposes -- it does affect the answer but we only wish to find a representation
of the word as a product of conjugates of relators, which representation it is
is not important -- though smaller ones are preferable for efficiency reasons.
The list below gives the cycles created by adding in non-tree 
edges as products of relator cycles.}
\begin{center}
\begin{tabular}{llll}   
$\underline{[g,x]}$  & $\mapsto \underline{h_1[g,x]}$ &&\\  
$[b,b] $ & $\mapsto ~ bb(a^2)^{-1}$     & $\to ~ b^2a^{-2}$ & $\to ~ su^{-1}$,\\  
$[ab,a]$ & $\mapsto ~ ab(b)^{-1}$       & $\to ~ aba^2b$    & $\to ~ t$,\\  
$[ab,b]$ & $\mapsto ~ ab^2(a^3)^{-1}$   & $\to ~ ab^2a$     & $\to ~ u^a r^{-1}$,\\  
$[ba,a]$ & $\mapsto ~ ba^2(a^2b)^{-1}$  & $\to ~ ba^4ba^6$  & $\to ~ s u^{-1} u^{ba^{-2}} s^{-a^{-2}}$,\\  
$[ba,b]$ & $\mapsto ~ bab(a)^{-1}$      & $\to ~ baba^3$    & $\to ~ t^{-a^{-1}b^{-1}} u^{b^{-1}}$,\\  
$[a^3,a]$ & $\mapsto ~ a^4(\id)^{-1}$   & $\to ~ a^4$       & $\to ~ r$,\\  
$[a^3,b]$ & $\mapsto ~ a^3b(ba)^{-1}$   & $\to ~ a^3ba^5b$  & $\to ~ rt^{-a}$,\\  
$[a^2b,a]$ & $\mapsto ~ a^2ba(ab)^{-1}$ & $\to ~ a^2ba^3ba^3$ & $\to ~ t^{a^{-1}}$,\\  
$[a^2b,b]$ & $\mapsto ~ a^2b^2(\id)^{-1}$ & $\to ~ a^2b^2$ & $\to u$.  
\end{tabular}
\end{center}
\emph{This example gives 32 generators for the module of identities. 
In fact this can be reduced to 7 but the reduction requires methods not 
dealt with in this thesis. }
\end{example}  
  
\section{Algorithm for Computing a Set of Generators for $\pi_2$}

Section 5.3 described how the problem of specifying a free crossed  
resolution of a group reduced to the problem of defining 
a contracting homotopy of a covering crossed complex.\\
  
The computation of a complete rewriting system for the group is used to  
define the first part of the 
contracting homotopy $h_1$ on the edges of the Cayley graph.  
The formulae from the definition of the covering crossed complex are  
used to  
find a complete set of generators for the kernel of $\dd_2$ (the  
identities  
among relations).  
The pre-images of these elements generate $C_3$ as a $\bZ  
G$-module.  
By reducing this set of generators and writing each of the reducible  
generators  
in terms of the irreducible ones we define $h_2$ on the generators of $C_2$.  
This is made clear in the example, and is the part which corresponds to  
the Gr\"obner basis computation, though we do it by inspection.\\
  
Now the crossed complex formulae with $h_2$ are used to find a complete  
set of  
generators for the kernel of $\dd_3$ (the identities among identities).  
Again,  
we reduce the set of identities, so that their pre-images freely  
generate $C_4$  
as a $\bZ G$-module. The process of reduction of the identities defines  
the next  
contracting homotopy $h_3$, and again we use the formulae to find a  
complete  
set of generators for $ker\dd_4$, and reduction to a set whose pre-image  
freely  
generates $C_5$ as a $\bZ G$-module.\\
  
This procedure may in theory be repeated as much as is wished, in order  
to  
compute the resolution of the group up to any level. The limitations are  
ones of  
practicality: in our example the reduction of the set of identities is  
done by  
inspection (involving a lot of trial and error) this takes time (weeks).  
A Gr\"obner basis procedure (over the group ring) would provide a  
computerisable  
method for defining $h_n$, and this would mean that the computation of  
the  
resolution was limited only by the computer's capacity.  
The correspondence between the homotopy definition and the Gr\"obner  
basis  
computation (for reduction) is explained more fully in the next section.  
  
\subsection{Specification of the Program}  

A collection of $\mathsf{GAP3}$ functions has been written to perform these 
calculations and will be rewritten in $\mathsf{GAP4}$ and submitted as a share package.  
The function $\mathtt{IdRel1}$ accepts as input a free group and a list of relators.  
It goes through a number of calculations, including an ``extra  
information''  
Knuth-Bendix completion procedure and returns a complete set of  
generators  
for the module of identities among relations.  
The structure of the program $\mathtt{idrel.g}$ is outlined below.\\
  
Preliminary functions necessary are:
\begin{enumerate}[]
\item  
$\mathtt{ReduceWord}( word, R1 )$:  
reduces a word with respect to a rewriting system $R1$, in the standard way.  
\item
$\mathtt{ReduceWord2}( word, R2 )$:  
applies an EIRS $R2$  
to a word and reduces it as far as possible within that system.  
Output is a pair $[c, w]$ where $word = \dd_2(c) w$,  
where $c$ is a Y-sequence.  
\item
$\mathtt{InverseYsequence}( a )$:  
Y-sequences are represented by lists  
$a=[s_1,u_1],\ldots, [s_n, u_n]$  
where $u_i \in F$ and $s_i$ is a relator or an inverse of 
a relator.  
This function inverts such a sequence to  
$[s_n,u_n], \ldots,[s_1,u_1]$.  
This is used to invert products of conjugates of relators which are  
represented as Y-sequences.  
\item
$\mathtt{KB2}( R2 )$:  
is an implementation of the ``extra information'' Knuth-Bendix  
procedure described in Section 4. The input rules are  
in the form of lists of length three where the middle entry represents  
the product of conjugates of relators $(r_1,u_1)^{\ep_1} \cdots  
(r_n,u_n)^{\ep_n}$  
as a Y-sequence $[[r_1^{\ep_1},u_1],\ldots,[r_n^{\ep_n},u_n]]$. The  
output rules  
will have the same form. If $[l,c,r]$ is a rule in such a system then  
$l \to r$ and $l= \dd_2(c) r$.
\end{enumerate}  
  
Given a presentation $grp \lan X | relts \ran$, define $F:=F(X)$.
The main function is:

\begin{enumerate}[]  
\item
$\mathtt{IdRel1}(F,relts)$.  
First $G$ is defined to be the quotient of the free group $F$ by the  
relators $relts$. Let $\theta:F \to G$ be the quotient morphism.
It is necessary to keep track of whether an element is in $G$ or $F$.
The next step is to construct the initial EIRS from the relators. 
The program uses the monoid  
presentation of the group to enable it to accept relators containing  
inverses without changing them.  
The resulting EIRS is then completed using $\mathtt{KB2}$ to obtain  
$K2$. The analogous ordinary system is $K1$.  
The Cayley graph is represented by a list of edges, which are pairs  
$[g,x]$  
where $g$ is an irreducible in $F$ and $x$ is a generator.  
The so-called alpha-edges are the edges not in the spanning tree given by the  
length-lex order. The map $h_1$ is defined on these alpha-edges by 
$h_1[g,x]:=  
[\id,\mathtt{ReduceWord2}( N(g)x, K2 )]$  
and we apply $p_2$ immediately, so recording only the second part of  
this pair.  
To obtain the identities among relations all relator cycles in the  
Cayley  
graph must be considered. These are recorded as pairs $[g,r]$  
where $g$ is a vertex and $r$ is a relator.  
The boundary $\db_2$ of the cycle is basically found by splitting up the  
relator $r$ to obtain a list of edges. Non-alpha edges are removed since  
$h_1$ maps any edge of the tree to $\id$. The remaining edges of each  
cycle are identified with their images under $p_2h_1$.  
The identities are calculated by manipulating the information held so as  
to  
obtain a representation  of $p_2 (h_1\db_2[g,r])^{-1} r^{\sigma(g)^{-1}}$ for
each $[g,r]$ pair.
\end{enumerate}  
  
The output is in the form of a record $\mathtt{id1}$ (say) with the
following  
fields:
\begin{enumerate}[]
\item  
$\mathtt{id1.free}$ the free group $F$;
\item  
$\mathtt{id1.rels}$ the relators $relts$;
\item  
$\mathtt{id1.elF}$ the normal forms of the group elements;
\item  
$\mathtt{id1.K}$ the (ordinary) completed rewriting system;
\item  
$\mathtt{id1.idents}$ the generating set of identities among relations;
\item  
$\mathtt{id1.isIdsRecord}$ true -- a check that the identities generated  
all have the image $\id$.
\end{enumerate}
  
A small example is printed here -- others are on disk in files  
$\mathtt{idreleg1.g}$ to  
$\mathtt{idreleg3.g}$.  
If $\mathtt{IdRelPrintLevel}$ is set to be greater than 1 (up to 3)  
information on the progression through the program is printed to the  
screen.  
{\small  
\begin{verbatim}  
gap> Read("idrel.g");  
gap> IdRelPrintLevel:=1;;  
gap> F:=FreeGroup("a","b");;  
gap> a:=F.1;;b:=F.2;;  
gap> R:=[a^3,b^2,a*b*a*b];;  
gap> id1:=IdRel1(F,R);;  
gap> id1.idents;  
[ [ [ r1-1, IdWord ], [ r1^-1, IdWord ] ], 
  [ [ r1^-1, IdWord ], [ r1, a^-1 ] ],  
  [ [ r3^-1, IdWord ], [ r2, a^-1*b^-1*a^-1 ],  
    [ r2^-1, IdWord ], [ r1^-1, b^-1 ], [ r3, a^-2*b^-1 ], [ r1, b^-1 ] ],  
  [ [ r1^-1, IdWord ], [ r1, a^-2 ] ],  
  [ [ r2^-1, IdWord ], [ r1^-1, b^-1 ], [ r3, a^-2*b^-1 ],  
      [ r3^-1, IdWord ], [ r2, a^-1*b^-1*a^-1 ], [ r1, b^-1*a^-1 ] ],  
  [ [ r3^-1, IdWord ], [ r2, a^-1*b^-1*a^-1 ],  
      [ r2^-1, IdWord ], [ r1^-1, b^-1 ], [ r3, a^-2*b^-1 ],  
      [ r1, a^-1*b^-1 ] ], [ [ r2^-1, IdWord ], [ r2, IdWord ] ],  
  [ [ r2^-1, a^-1 ], [ r2, a^-1 ] ], 
  [ [ r2^-1, IdWord ], [ r2, b^-1 ] ],  
  [ [ r3^-1, a^-2 ], [ r1, IdWord ], [ r2^-1, a^-1*b^-1 ], 
  [ r1^-1, IdWord ], [ r3, a^-2 ], [ r2, a^-2 ] ],  
  [ [ r2^-1, a^-1 ], [ r2, b^-1*a^-1 ] ],  
  [ [ r2^-1, a^-1*b^-1 ], [ r1^-1, IdWord ], [ r3, a^-2 ],  
      [ r3^-1, a^-2 ], [ r1, IdWord ], [ r2, a^-1*b^-1 ] ],  
  [ [ r2^-1, IdWord ], [ r3^-1, IdWord ], [ r2, a^-1*b^-1*a^-1 ], [ r3, IdWord ] ],  
  [ [ r2^-1, a^-1 ], [ r2^-1, IdWord ], [ r1^-1, b^-1 ], [ r3, a^-2*b^-1 ],  
      [ r2^-1, a^-1*b^-1 ], [ r1^-1, IdWord ], [ r3, a^-2 ], [ r3, a^-1 ] ],  
  [ [ r1^-1, IdWord ], [ r3^-1, a^-2 ], [ r1, IdWord ], [ r3, b^-1 ] ],  
  [ [ r3^-1, a^-2 ], [ r1, IdWord ], [ r1^-1, IdWord ], [ r3, a^-2 ] ],  
  [ [ r2^-1, IdWord ], [ r3^-1, IdWord ], [ r2, a^-1*b^-1*a^-1 ], [ r3, b^-1*a^-1 ] ],  
  [ [ r2^-1, a^-1*b^-1 ], [ r1^-1, IdWord ], [ r3, a^-2 ], [ r2^-1, a^-1 ],  
    [ r2^-1, IdWord ], [ r1^-1, b^-1 ], [ r3, a^-2*b^-1 ], [ r3, a^-1*b^-1 ] ] ]  
gap> 
\end{verbatim} } 
  
The program returns a set of 18 generators for $ker \dd_2$, these are  
the images under $\dd_3$ of a set of generators for $\Cb_3$.  
For the output of  
higher stages to be useful implementation of some Gr\"obner basis  
procedures will be necessary. This is discussed in Section 6.  
  
\begin{example}\emph{  
We now present the results obtained for $S_3$  
followed by some of the details of the calculations which can be done by  
hand in this case, beginning with the presentation  
$$  
G:=grp\lan x,y \ | \ x^3,y^2,(xy)^2\ran .  
$$ 
The description of the partial free crossed resolution is as follows.  
Let $X=\{x,y\}$ and  
define $\mathcal{R}$ to be the set of relator labels $\{r,s,t\}$  
whose images under $w$ are  
$$\{ x^3, y^2, (xy)^2\}.$$  
$C_2$ is the free crossed $F(X)$-module on $w:\mathcal{R}\to F(X)$.\\  
$C_3$ is the free $\bZ G$-module generated by four elements  
$\{\iota_1,\ldots ,\iota_4\}$ whose images under $\dd_3$ generate $ker  
\dd_2$  
and are  
\begin{center}  
$\{ r^{-1}r^{x^{-1}}, \ s^{-1}s^{y^{-1}}, \ t^{-1}t^{y^{-1}x^{-1}}, \  
 \  ts^{-xy}r^{-y}s^{-1}t^xs^{-x}r^{-1}t^{x^{-1}} \}.$  
\end{center}  
$C_4$ is the free $\bZ G$-module generated by five elements  
$\{\eta_1, \ldots ,\eta_5\}$ whose images under $\dd_4$ generate $ker  
\dd_3$  
and are  
\begin{center}  
$  
\{ \iota_1(\id+x+x^2), \ \iota_2(\id+y), \ \iota_3(x+y),  
 \ \iota_4(x^2-\id)-\iota_2(yx+x^2)-\iota_1(xy-\id),  
 \ \iota_4(y-1)-\iota_3(x-yx+\id)+\iota_2 \}.  
$  
\end{center}  
$C_5$ is the free $\bZ G$-module generated by six elements  
$\{\mu_1, \ldots,\mu_6\}$ whose images under $\dd_5$ generate $ker  
\dd_4$ and  
are  
\begin{center}  
$  
\{ \eta_1(x-\id), \ \eta_2(y-\id), \ \eta_3(x^2-y),  
 \ \eta_4(\id+x+x^2)+ \eta_2(\id+x+x^2)- \eta_1(\id-y),$\\  
$ \eta_5(\id+yx)+ \eta_4(x+y)+\eta_3+ \eta_2(x^2),  
 \ \eta_5(\id+y)+ \eta_3(\id-x+y)- \eta_2 \}.  
$  
\end{center}  
$C_6$ is the free $\bZ G$-module generated by seven elements  
$\{\nu_1,\ldots,\nu_7\}$ whose images under $\dd_6$ generate $ker \dd_5$  
and  
are  
\begin{center}  
$  
\{ \mu_1(\id+x+x^2), \  
   \mu_2(\id+y), \  
   \mu_3(x+y),  
 \ \mu_4(x^2-\id) - \mu_1(x^2+y),$\\  
$ \mu_6(\id+x+x^2) - \mu_5(\id+y+xy) + \mu_4(\id+y) -\mu_3(y) -  
\mu_2(x^2),$\\  
$ \mu_5(x^2-y)+ \mu_2(x)- \mu_3, \  
   \mu_6(yx-x) -\mu_3(\id+x+y) \}.  
$  
\end{center}  
This defines the resolution of the group ($C_0$) up to the sixth level  
$C_6$.  
If identities among relations $\iota_i$ are equivalent to first order  
syzygies  
then the $\nu_i$ are like the fourth order syzygies.\\}  
  
\emph{  
The calculations proceeded as follows:\\  
First of all we computed an ``extra information'' complete rewriting  
system for  
the group (GAP output):}
{\small  
\begin{verbatim}  
gap> R:=[x^3,y^2,x*y*x*y];  
[ x^3, y^2, x*y*x*y ]  
gap> R2:=List( R, r -> [ r, [ [ r, IdWord ] ], IdWord ] );  
[ [ x^3, [ [ x^3, IdWord ] ], IdWord ], [ y^2, [ [ y^2, IdWord ] ],
IdWord ],  
  [ x*y*x*y, [ [ x*y*x*y, IdWord ] ], IdWord ] ]  
gap> KB2(R2);  
[ [ y^2, [ [ y^2, IdWord ] ], IdWord ], 
  [ x^3, [ [ x^3, IdWord ] ], IdWord ],  
  [ x^2*y, [ [ y^-1*x^-1*y^-1*x^-1, x^-2 ], 
            [ y^2, x^-1*y^-1*x^-3 ],  [ x^3, IdWord ] ],  y*x ],  
  [ x*y*x, [ [ y^-2, x^-1*y^-1*x^-1 ], [ x*y*x*y, IdWord ] ], y ],  
  [ y*x^2, [ [ y^-1*x^-1*y^-1*x^-1, x^-2*y^-1 ], 
             [ x^3, y^-1 ], [ y^2, IdWord ] ], x*y ],  
  [ y*x*y, [ [ x^-3, IdWord ], [ x*y*x*y, x^-2 ] ], x^2 ] ]  
\end{verbatim} } 
\emph{ The six rules may be translated as follows:}  
\begin{center}  
\begin{tabular}{ll}  
$y^2 \to_{R2}  (s, \id)$ & $x^3 \to_{R2}  (r, \id)$\\  
$x^2y \to_{R2} (t^{-x^{-2}} s^{ x^{-1} y^{-1} x^{-3} } r , yx)$  
                   & $xyx \to_{R2}  (s^{-x^{-1} y^{-1} x^{-1} } t, y)$ \\  
$yxy \to_{R2} (r^{-1} t^{x^{-2}}, x^2)$  
                   & $yx^2 \to_{R2} (t^{-x^{-2} y^{-1}} r^{y^{-1}} s, xy)$\\  
\end{tabular}  
\end{center}  
\emph{The word on the left hand side reduces to the word at the right
hand  
end, and  
is equal to the boundary of the entry in brackets multiplied by that  
reduced word. $N(g)$  
denotes the normal form (unique reduced word) in $F(X)$ representing  
the element $g$ and $\theta$ is the quotient map : $F(X) \to G$.  
The homotopy $h_1$ is defined on the edges $[g,x]$ of the Cayley graph  
($G \times X$) by finding products of conjugates of the relators ($R$)  
whose  
images under $\delta_2$ are $N(g)xN(g\theta(x))^{-1}$. (For small groups  
like  
this one it is possible to do this quite efficiently by inspection.)  
In general one defines $h_1$ algorithmically by using the ``extra  
information'' rewriting system introduced in the previous section.  
The definition of $h_1$ in this example is as follows: (I have chosen to  
use a  
more efficient definition than that suggested by the computer program  
because  
it simplifies the manual calculations to follow. The only loss by using  
the computer generated definition is that of space. With  groups even  
a little larger or more complex there is no option but to use the  
computer  
 generated definition.)}  
\begin{center}  
\begin{tabular}{l|l|l}  
edge $[g,x]$ & $h_1[g,x]$ & $p_2h_1[g,x]$ \\  
{\em in} $\Cb_1$ & {\em in} $\Cb_2$ & {\em in} $C_2$ \\  
\hline  
   $[\id,x]$ & $1$                   & $1$ \\  
   $[\id,x]$ & $1$                   & $1$ \\  
   $[x,x]$   & $1$                   & $1$ \\  
   $[x,y]$   & $1$                   & $1$ \\  
   $[y,x]$   & $1$                   & $1$ \\  
   $[y,y]$   & $[\id,s]$             & $s$ \\  
   $[x^2,x]$ & $[\id,r]$             & $r$ \\  
   $[x^2,y]$ & $[\id,rs^xt^{-x}]$    & $rs^xt^{-x}$ \\  
   $[xy,x]$  & $[\id,ts^{-1}]$       & $ts^{-1}$ \\  
   $[xy,y]$  & $[\id,ts^{xy}t^{-1}]$ & $ts^{xy}t^{-1}$ \\  
   $[yx,x]$  & $[\id,sr^yt^{-1}]$    & $sr^yt^{-1}$ \\  
   $[yx,x]$  & $[\id,st^ys^{-1}r^{-1}]$ & $st^ys^{-1}r^{-1}$ \\  
\end{tabular}  
\end{center}  
$$\text{\em Table 1: Defining } h_1$$  
\emph{ The formulae for the crossed complex give us a complete set of  
generators for  $ker \dd_2$. }  

\begin{center}  
\begin{tabular}{l|l|l|l}  
$[g,r]$          & $\db_2[g,r]  $         & $p_2((h_1\db_2[g,r])^{-1}[g,  
r]^{[g,g^{-1}]})$ & $p_3h_2[g,r]$ \\  
{\em in} $\Cb_2$   & {\em in} $\Cb_1$         & {\em in} $C_2$ & {\em in} $C_3$ \\  
\hline  
$[\id,r]$ & $[1,x][x,x][x^2,x]$        & $1$  
& $0$ \\  
$[x,r]$   & $[x,x][x^2,x][1,x]$        & $r^{-1}r^{x^{-1}}$  
& $\iota_1$ \\  
$[y,r]$   & $[y,x][yx,x][xy,x]$        & $1$  
& $0$ \\  
$[x^2,r]$ & $[x^2,x][1,x][x,x]$        & $r^{-1}r^{x^{-2}}$  
& $\iota_1(1+x^2)$ \\  
$[xy,r]$   & $[xy,x][y,x][yx,x]$       & $r^{-x^{-1}y^{-1}x^{-1}}r^{y^{-  
1}x^{-1}}$              & $-\iota_1(xy)$ \\  
$[yx,r]$  & $[yx,x][xy,x][y,x]$        & $r^{-y^{-1}}r^{x^{-1}y^{-1}}$  
& $\iota_1(y)$ \\  
$[\id,s]$ & $[1,y][y,y]$               & $1$  
& $0$ \\  
$[x,s]$   & $[x,y][xy,y]$              & $s^{-y^{-1}x^{-1}}s^{x^{-1}}$  
& $-\iota_2(x^2)$ \\  
$[y,s]$   & $[y,y][1,y]$               & $s^{-1}s^{y^{-1}}$  
& $\iota_2$ \\  
$[x^2,s]$ & $[x^2,y][yx,y]$            & $t^{y^{-1}x^{-3}}t^{-x^{-2}}$  
& $-\iota_3(x)$ \\  
$[xy,s]$  & $[xy,y][x,y]$              & $1$  
& $0$ \\  
$[yx,s]$  & $[yx,y][x^2,y]$            & $t^xs^{-x}t^{-y^{-1}}s^{-x^{-1}  
y^{-1}}$                    & $\iota_3(y)-\iota_2(yx)$ \\  
$[\id,t]$ & $[1,x][x,y][xy,x][y,y]$    & $1$  
& $0$ \\  
$[x,t]$   & $[x,x][x^2,y][yx,x][xy,y]$ & $ts^{-xy}r^{-y}s^{-1}t^x s^{-x}  
r^{-1} t^{x^{-1}}$        & $\iota_4$ \\  
$[y,t]$   & $[y,x][yx,y][x^2,x][1,y]$  & $1$  
& $0$ \\  
$[x^2,t]$ & $[x^2,x][1,y][y,x][yx,y]$  & $t^{y^{-1}x^{-3}} t^{-x^{-2}}$  
& $-\iota_3(x)$ \\  
$[xy,t]$  & $[xy,x][y,y][1,x][x,y]$    & $t^{-1}t^{y^{-1}x^{-1}}$  
& $\iota_3$ \\  
$[yx,t]$  & $[yx,x][xy,y][x,x][x^2,y]$ & $t^xs^{-x}r^{-1}ts^{-xy}r^{-y}  
s^{-1} t^{x^{-1}y^{-1}}$  & $\iota_4(1)-\iota_3(yx)$ \\  
\end{tabular}  
\end{center}  
$$\text{\em Table 2: Calculating $ker\dd_2$ and defining } h_2$$  
  
\emph{ The last column shows how the other identities found may be
expressed  (in $C_3$)  
in terms of the four generating ones. The main result so far is that the  
module of identities among relations for this group presentation is  
generated by four elements. This result can be obtained by other methods. 
However, we  
now use the results of that last column to calculate a set of generators  
for the  
module of identities among identities. This last column defines $h_2$ on  
the free  
generators of $\Cb_2$ (listed in the second column of the table) so that it 
 annihilates the action of $\Cb_1$ as required.}\\  

{\em The elements $p_3(-h_2\db_3[g,\iota]+[g,\iota]^{h_0(g)})$ for $[g,\iota]  
\in \Cb_3$  
are a generating set of identities among the identities.  
The table below gives the identity resulting from each generator  
$[g,\iota]$ of $\Cb_3$.  
These were obtained by first calculating the images under $\db_3$.  
This effectively gives us the boundary of the generator.}\\ 
 
{\em For example, $\db_3[\id,\iota_1]$ is $[\id,r]^{-1}[x,r]^{[x,{x^{-1}}]}$,  
This is because $\dd_3(\iota_1)=r^{-1}r^{x^{-1}}$, and  
$\db_n(g,\gamma):=[g,\dd_n(\gamma)]$ and we then write  
$[g,\dd_n(\gamma)]$ as a product of the  
generators of $C_{n-1}$  as a $C_1$-module as $h_2$ will be defined on  
these generators.  
Similarly, $\db_3[x^2,\iota_4]$ is  
$[x^2,t][y,s]^{-[y,{xy}]} [yx,r]^{-[yx,{y}]} [x^2,s]^{-1}[x,t]^{[x,x]}  
[x,s]^{-[x,x]} [x^2,r]^{-1}[\id,t]^{[\id,x^{-1}]}$.
(Recall that the action is defined as $[g,\gamma]^{[g,y]}=[g\theta y,  
\gamma^y]$.) \\}

{\em When we have turned the $[g,\iota]$ into such a product of $\Cb_2$  
generators,  
we can calculate $h_2(\db_3[g,\iota])$ using the last table.  
Note that a property of $h_2$ is that it must annihilate the action of  
$\Cb_1$,  
it is also a morphism, in that it preserves the multiplication of the  
elements of $\Cb_2$.  
Therefore \ $h_2\db_3[\id,\iota_1]$ is $h_2[\id,r]^{-1} = h_2[x,r]$ \\  
and \ $h_2\db_3[x^2,\iota_4]$ is  
$h_2[x^2,t] - h_2[y,s] - h_2[yx,r] - h_2[x^2,s] + h_2[x,t] - h_2[x,s] -  
h_2[x^2,r] + h_2[\id,t]$.  
We can read these values off the previous table, as we have defined  
$h_2$ on all the  
elements $[g,r]$.  
So \ $h_2\db_3[\id,\iota_1]$ is $[\id,-0+\iota_1]=[\id,\iota_1]$\\  
and \ $h_2\db_3[x^2,\iota_4]$ is  
$[\id, \iota_4-\iota_2-\iota_1(y)-(-\iota_3(x))+\iota_4-(-\iota_2(x^2))-  
\iota_1(1+x^2)+0]$.}\\  
  
\emph{ To obtain the identities we negate the above
$h_2\db_3[g,\iota]$'s and  
add $[g,\iota]^{h_0(g)}$  
which is effectively $[\id,\iota(g)]$.  
We finally project this sum down to $C_3$:  
$p_2h_2\db_3[\id,\iota_1]$ is $-\iota_1+\iota_1=0$ and  
$p_2h_2\db_3[x^2,\iota_4]$ is $\iota_4(x-1)-\iota_2(x^2-\id)+\iota_1(\id  
+x^2+y)$.}\\  
  
\emph{The following table gives the identities resulting from all the  
generators.}  
  
\begin{center}  
\begin{tabular}{l|l|l}  
$[g,\iota]$      & $p_3(-h_2\db_3[g,\iota]+[g,\iota]^{h_0(g)})$ &  
$p_4h_3[g,\iota]$\\  
in $\bar{C_3}$      & {\em in} $C_3$ & {\em in} $C_4$\\  
\hline  
$[\id,\iota_1]$  & $0$ & $0$\\  
$[x,\iota_1]$    & $0$ & $0$\\  
$[y,\iota_1]$    & $0$ & $0$\\  
$[x^2,\iota_1]$  & $\iota_1(\id+x+x^2)$ & $\eta_1$\\  
$[xy,\iota_1]$   & $0$ & $0$\\  
$[yx,\iota_1]$   & $\iota_1(y+xy+yx)$ & $\eta_1(y)$\\  
$[\id,\iota_2]$  & $0$ & $0$\\  
$[x,\iota_2]$    & $0$ & $0$\\  
$[y,\iota_2]$    & $\iota_2(\id+y)$ & $\eta_2$\\  
$[x^2,\iota_2]$  & $\iota_2(x+yx)-\iota_3(x+y)$ & $\eta_2(x)-\eta_3$\\  
$[xy,\iota_2]$   & $\iota_2(x^2+xy)$ & $\eta_2(x^2)$\\  
$[yx,\iota_2]$   & $\iota_3(x+y)$ & $\eta_3$\\  
$[\id,\iota_3]$  & $0$ & $0$\\  
$[x,\iota_3]$    & $\iota_3(x^2+yx)$ & $\eta_3(x)$\\  
$[y,\iota_3]$    & $\iota_3(x+y)$ & $\eta_3$\\  
$[x^2,\iota_3]$  & $0$ & $0$\\  
$[xy,\iota_3]$   & $\iota_3(\id+xy)$ & $\eta_3(y)$\\  
$[yx,\iota_3]$   & $0$ & $0$\\  
$[\id,\iota_4]$  & $0$ & $0$\\  
$[x,\iota_4]$ &  $\iota_4(x^2\!-\!\id)\!+\!\iota_3(x\!+\!y)\!-\!\iota_2(  
yx\!+\!x^2) \! - \! \iota_1(xy \! - \! 1)$  
 & $\eta_4+\eta_3$\\  
$[y,\iota_4]$    & $\iota_4(y-1)\!-\!\iota_3(x\!-\!yx\!+\!\id)\!+\!  
\iota_2$  
 & $\eta_5$\\  
$[x^2,\iota_4]$  & $\iota_4(x\!-\!1)-\iota_2(x^2\!-\!\id)\!+\!\iota_1  
(\id \!+\!x^2 \! + \! y)$  
 & $-\eta_4(x)-\eta_2(x^2)-\eta_1$\\  
$[xy,\iota_4]$   & $\iota_4(xy-1)-\iota_3(\id\!-\!yx\!-\!y) \! -  
\!\iota_2(yx)\!- \! \iota_1(xy \! - \! \id)$  
 & $\eta_5(x^2)+\eta_4+\eta_3(x)$\\  
$[yx,\iota_4]$   & $\iota_4(yx\!-\!\id)\!-\!\iota_3(\id\!-\!y\!-\!yx)\!-  
\!\iota_2(x^2 \! + \! yx \! - \! \id) \! + \! \iota_1(x^2 \! + \! \id \!  
+ \! y)$  
 & $-\eta_5(yx)-\eta_4(x)-\eta_2(x^2)+\eta_1$\\  
\end{tabular}  
\end{center}  
$$\text{\em Table 3: Calculating $ker\dd_3$ and defining }h_3$$  
  
\emph{  
The images of the $\eta_i$ generate the kernel as a $\bZ G$-module, the  
$\eta_i$ themselves provide a set of generators for $\bar{C_4}$.  
We use the formula $p_4(-h_3\db_4[g,\eta]+[g,\eta]^{h_0(g)})$ to  
calculate  
a generating set of 30 elements for $ker\dd_4$, which we can reduce to  
six.  
The last table defines $h_3$ (``in $\Cb_4$'' column) on the generators of $C_3$
($[g,\iota]$ column).}  
  
\begin{center}  
\begin{tabular}{l|l|l}  
$[g,\eta]$ &  $p_4(-h_3\db_4[g,\eta]+[g,\eta]^{h_0(g)})$ &  
$p_5h_4[g,\eta]$\\  
{\em in} $\bar{C_4}$ &  {\em in} $C_4$ & {\em in} $C_5$\\  
\hline  
$[\id,\eta_1]$  & $-\eta_1+\eta_1$ &            $0$\\  
$[x,\eta_1]$    & $-\eta_1+\eta_1(x^2)$ &       $-\mu_1(x^2)$\\  
$[y,\eta_1]$    & $-\eta_1(y)+\eta_1(y)$ &      $0$\\  
$[x^2,\eta_1]$  & $-\eta_1+\eta_1(x)$ &         $\mu_1$\\  
$[xy,\eta_1]$   & $-\eta_1(y)+\eta_1(xy)$ &     $\mu_1(y)$\\  
$[yx,\eta_1]$   & $-\eta_1(y)+\eta_1(yx)$ &     $-\mu_1(yx)$\\  
$[\id,\eta_2]$  & $-\eta_2+\eta_2$ &            $0$\\  
$[x,\eta_2]$    & $-\eta_2(x^2)+\eta_2(x^2)$ &  $0$\\  
$[y,\eta_2]$    & $-\eta_2+\eta_2(y)$ &         $\mu_2$\\  
$[x^2,\eta_2]$  & $-\eta_2(x)+\eta_2(x)$ &      $0$\\  
$[xy,\eta_2]$   & $-\eta_2(x^2)+\eta_2(xy)$ &   $\mu_2(x^2)$\\  
$[yx,\eta_2]$   & $-\eta_2(x)+\eta_2(yx)$ &     $\mu_2(x)$\\  
$[\id,\eta_3]$  & $-\eta_3+\eta_3$ &            $0$\\  
$[x,\eta_3]$    & $-\eta_3(y)+\eta_3(x^2)$ &    $\mu_3$\\  
$[y,\eta_3]$    & $-\eta_3(y)+\eta_3(y)$ &      $0$\\  
$[x^2,\eta_3]$  & $-\eta_3(x)+\eta_3(x)$ &      $0$\\  
$[xy,\eta_3]$   & $-\eta_3(x)+\eta_3(xy)$ &     $\mu_3(yx)$\\  
$[yx,\eta_3]$   & $-\eta_3+\eta_3(yx)$ &        $\mu_3(y)$\\  
$[\id,\eta_4]$  & $-\eta_4+\eta_4$ &            $0$\\  
$[x,\eta_4]$    & $\eta_4(x+\id)+\eta_2(x^2+x+\id)-\eta_1(\id-  
y)+\eta_4(x^2)$ &        $\mu_4$\\  
$[y,\eta_4]$    & $\eta_5(\id+yx)+\eta_4(x)+\eta_3+\eta_2(x^2)+\eta_4(y)  
$ &            $\mu_5$\\  
$[x^2,\eta_4]$  & $-\eta_4(x)+\eta_4(x)$ &  
$0$\\  
$[xy,\eta_4]$   & $\eta_5(x^2+\id)+\eta_4+\eta_3(x-  
\id)+\eta_2(x+\id)+\eta_4(xy)$ &    $-\mu_6+\mu_5(y)-\mu_2(x)$\\  
$[yx,\eta_4]$   & $-\eta_5(x^2+yx)-\eta_4(x+\id)-\eta_3(x)+\eta_1(y-  
\id)+\eta_4(yx)$ & $-\mu_6(x^2+x)+\mu_5(xy)$\\  
                &  
& \quad $-\mu_4+\mu_3(y)-\mu_2$\\  
$[\id,\eta_5]$  & $-\eta_5+\eta_5$ &  
$0$\\  
$[x,\eta_5]$    & $-\eta_5(x^2)+\eta_5(x^2)$ &  
$0$\\  
$[y,\eta_5]$    & $\eta_5+\eta_3(\id-x+y)-\eta_2+\eta_5(y)$ &  
$\mu_6$\\  
$[x^2,\eta_5]$  & $\eta_5(yx)+\eta_3(x-y+\id)-\eta_2(x)+\eta_5(x)$ &  
$\mu_6(x)+\mu_3(\id+x)$\\  
$[xy,\eta_5]$   & $\eta_5(x^2)+\eta_3(y-\id+x)-\eta_2(x^2)+\eta_5(xy)$ &  
$\mu_6(x)+\mu_3(x^2-\id)$\\  
$[yx,\eta_5]$   & $-\eta_5(yx)+\eta_5(yx)$ &  
$0$\\  
\end{tabular}  
\end{center}  
$$\text{\em Table 4: Calculating $ker\dd_4$ and defining }h_4$$  
  
\emph{ So now we have six generators for $\Cb_5$ : $\{\mu_1,\ldots,\mu_6\}$  
and their images ~ 
$
\{\eta_1(x-\id), \ \eta_2(y-\id), \ \eta_3(x^2-y), \  
\eta_4(\id+x+x^2)+\eta_2(\id+x+x^2)-\eta_1(\id-y), \  
\eta_5(\id+yx)+\eta_4(x+y)+\eta_3+\eta_2(x^2), \  
\eta_5(\id+y)+\eta_3(\id-x+y)-\eta_2\}
$ ~
generate the module of identities among  
the identities among identities ($ker\dd_4$).  
The last column defines $h_4$. }\\  
  
\begin{center}  
\begin{tabular}{l|l|l}  
$[g,\mu]$ &  $p_5(-h_4\db_5[g,\mu]+[g,\mu]^{h_0(g)})$ & {\em in} $C_6$\\  
\hline  
$[\id,\mu_1]$  & $0$                &            $0$\\  
$[x,\mu_1]$    & $0$                &            $0$\\  
$[y,\mu_1]$    & $0$                &            $0$\\  
$[x^2,\mu_1]$  & $\mu_1(\id+x+x^2)$ &            $\nu_1$\\  
$[xy,\mu_1]$   & $\mu_1(y+xy+yx)$   &            $\nu_1(y)$\\  
$[yx,\mu_1]$   & $0$                &            $0$\\  
  
$[\id,\mu_2]$  & $0$                &            $0$\\  
$[x,\mu_2]$    & $0$                &            $0$\\  
$[y,\mu_2]$    & $\mu_2(\id+y)$     &            $\nu_2$\\  
$[x^2,\mu_2]$  & $0$                &            $0$\\  
$[xy,\mu_2]$   & $\mu_2(x^2+xy)$    &            $\nu_2(x^2)$\\  
$[yx,\mu_2]$   & $\mu_2(x+yx)$      &            $\nu_2(x)$\\  
  
$[\id,\mu_3]$  & $0$                &            $0$\\  
$[x,\mu_3]$    & $\mu_3(xy+x^2)$    &            $\nu_3(x)$\\  
$[y,\mu_3]$    & $0$                &            $0$\\  
$[x^2,\mu_3]$  & $\mu_3(x+y)$       &            $\nu_3$\\  
$[xy,\mu_3]$   & $\mu_3(\id+xy)$    &            $\nu_3(y)$\\  
$[yx,\mu_3]$   & $0$                &            $0$\\  
  
$[\id,\mu_4]$  & $0$                &            $0$\\  
$[x,\mu_4]$    & $\mu_4(x^2\!-\!\id)-\mu_1(x^2\!+\!y)$ &  $\nu_4$\\  
$[y,\mu_4]$    & $\mu_6(\id\!+\!x\!+\!x^2)-\mu_5(1\!+\!y\!+\!xy)+\mu_4(1  
\!+\!y)-\mu_3(y)-\mu_2(x^2)$ &   $\nu_5$\\  
$[x^2,\mu_4]$  & $\mu_4(x-\id)+\mu_1(yx+\id)$ &    $-\nu_4(x)$\\  
$[xy,\mu_4]$   & $\mu_6(\id\!+\!x\!+\!x^2)-\mu_5(1\!+\!y\!+\!xy)+\mu_4(1  
\!+\!xy)-\mu_3(y)-\mu_2(x^2)+\mu_1(y\!+\!x^2)$ &   $\nu_5-\nu_4(xy)$\\  
$[yx,\mu_4]$   & $\mu_6(\id\!+\!x\!+\!x^2)-\mu_5(1\!+\!y\!+\!xy)+\mu_4(1  
\!+\!yx)-\mu_3(y)-\mu_2(x^2)-\mu_1(yx\!+\!\id)$ &    $\nu_5+\nu_4(y)$\\  
  
$[\id,\mu_5]$  & $0$                &            $0$\\  
$[x,\mu_5]$    & $\mu_5(x^2-y)+\mu_2(x)-\mu_3$ & $\nu_6$\\  
$[y,\mu_5]$    & $0$                &            $0$\\  
$[x^2,\mu_5]$  & $-\mu_6(x^2\!+\!x\!+\!\id)+\mu_5(\id\!+\!y\!+\!x)-\mu_4  
(\id\!+\!y)-\mu_3(x\!+\!x^2)+\mu_2(\id\!+\!x^2)$ &  
$\nu_6(x^2)-\nu_5\!+\!\nu_3$\\  
$[xy,\mu_5]$   & $-\mu_6(x^2\!+\!x\!+\!\id)+\mu_5(\id\!+\!y\!+\!yx)-  
\mu_4(\id\!+\!y)-\mu_3(x\!+\!x^2\!+\!yx)+\mu_2(x^2)$ &  
$-\nu_5-\nu_3(x\!+\!\id)$\\  
$[yx,\mu_5]$   & $\mu_5(yx\!-\!\id)-\mu_3(y)-\mu_2(x^2)$ &  
$\nu_6(y)-\nu_2(x^2)$\\  
  
$[\id,\mu_6]$  & $0$                &            $0$\\  
$[x,\mu_6]$    & $-\mu_3(x^2+yx)$   &            $-\nu_3(x^2)$\\  
$[y,\mu_6]$    & $\mu_6(y\!-\!\id)+\mu_3(yx)$  &     $\nu_7$\\  
$[x^2,\mu_6]$  & $-\mu_3(x+y)$      &            $-\nu_3$\\  
$[xy,\mu_6]$   & $\mu_6(xy\!-\!x^2)+\mu_3(y\!-\!x^2\!+\!yx)+\mu_2(x^2)$  
& $\nu_7(x^2)+\nu_3(x)$\\  
$[yx,\mu_6]$   & $\mu_6(yx\!-\!x)-\mu_3(x\!+\!y+\!\id)+\mu_2(x)$ &  
$\nu_7(x)-\nu_3(y+\id)$\\  
\end{tabular}  
\end{center}  
$$\text{\em Table 5: Calculating $ker\dd_5$ and defining }h_5$$  
  
\emph{We could calculate the identities for the next level,  
using the last table as a definition  
for $h_5$, computing a set of 42 generators for $ker\dd_6$  
(using  $p_6(-h_5\db_6[g,\nu]+[g,\nu]^{h_0(g)})$ ) and reducing them  
as before. It does not get more complicated: for $n \! \geq \! 3$  \  
$C_n$ is a  
$\bZ G$-module and the expression  
$p_n(-h_{n-1}\db_n[g,\gamma]+[g,\gamma]^{h_0(g)})$,  
where $\gamma$ is a generator of $C_n$,  
gives a set of generators for $C_{n+1}$ as a $\bZ G$-module  
(which may be reduced over the $\bZ G$-module).  
It is in principle possible to continue this exercise further, but  
it is  
not of value to do so here. The obvious conjecture it that $C_n$ will be  
the  
free $\bZ G$-module generated by $n+1$ elements.}  
\end{example}  
  
Notice that every time we are choosing a set of  
independent generators for the $\bZ G$-submodule; the set is not unique,  
and we do  
not have an algorithm for determining which generator is expressible in  
terms of  
the others or how to express it in this way. The method used is no more  
than  
inspection and trial and error. The purpose of including this example is  
that it  
best shows what may be achieved using the covering groupoids and  
homotopies  
methods, the complexity of even a very small example, and thus  
illustrates  
the necessity for a computer algorithm to extract such information as  
was  
summarised at the beginning of this example.  
The next section shows that these problems can be expressed in terms of  
noncommutative Gr\"obner bases over group rings. New work is being  
developed  
\cite{Birgit} on algorithms for such problems, and so expressing the  
problem  
of devising an algorithm for obtaining \emph{reduced} sets of identities  
and  
higher identities is a step forward, and until such Gr\"obner basis  
algorithms become  
available we cannot expect to be able to have algorithms for reducing  
the  
sets of generating identities.

\section{The Submodule Problem}

The previous sections have shown that a variation of the  
noncommutative Buchberger algorithm (Knuth-Bendix algorithm)  
may be applied to a group presentation to obtain the contracting  
homotopy $h_1$,  
and a set of generators for the module of identities among relations for  
the  
group presentation.  
This much has been implemented in the program $\mathtt{idrel.g}$ for 
$\mathsf{GAP}$.  
The remaining problem is that of reducing the set of generators with  
respect to the action of $\bZ G$ on the module.\\
  
We discussed earlier the Peiffer Problem which occurs at the first level  
(identities among relations: $ker\dd_2 \subseteq C(R)$).  
This problem is difficult because we need to test for equality  
in the free crossed $F(X)$-module, in other words, to test for Peiffer  
equivalence of two sequences  
(recall that the Peiffer rules imply that  
$[s,v][r,u]=[r,u][s,v \dd(r)^u]=[r,u \dd(s)^v]$).  
In this case we essentially wish to be able to reduce the set of  
generating  
identities to a set $\{\iota_1, \ldots, \iota_k\}$ that is in some sense
minimal over  
$\bZ G$ i.e. no $\iota_j$ can be written as a sum of $\bZ G$-multiples  
of the  
other identities. To summarise -- there are great difficulties in  
reducing the set of generators of the module of identities among  
relations.  
Furthermore, unless we can express each of the original generators in terms 
of those in the reduced set it is not practical to define 
$h_2$ on such a large set.\\  
  
We will now use a property which converts the Peiffer Problem into a  
Gr\"obner basis problem.  
This property is fully explained in \cite{BrHu}.  
First, recall that the crossed module is defined by taking the Peiffer  
equivalence classes of the free group $F(R \times F(X))$. This is the  
same as  
looking at the free monoid $(Y^+ \sqcup Y^-)^*$ factored by the  
relations needed for the group as well as by the Peiffer relations.  
Elements of  
$(Y^+ \sqcup Y^-)^*$ are called \textbf{Y-sequences}.\\  
  
An \textbf{identity Y-sequence} is one whose image under $\dd_2$ is the  
identity  
in $F(X)$.\\  
  
The identity property uses a result on the abelianisation of $C(R)$ to
describe a useful way of determining whether an identity  
$Y$-sequence  
(i.e. one identified with an element of the kernel of $\dd_2$, which is
abelian) is Peiffer equivalent to the empty sequence.\\  
  
An identity $Y$-sequence $a=(r_1,u_1)^{\ep_1},\ldots,(r_k,u_k)^{\ep_1}$  
has the \textbf{Primary Identity Property} if the indexing numbers  
$1, \ldots,k$ of the sequence $y$  
can be paired $(i,j)$ so that $r_i=r_j$, $\theta(u_j) = \theta(u_j)$ and  
$\ep_i = -\ep_j$.  
  
\begin{lem}[\cite{BrHu}]  
Let $a \in (Y^+ \sqcup Y^-)^*$. Then $a$ has the Primary Identity  
Property if and only if it is Peiffer equivalent to the empty sequence.  
\end{lem}

Let $X$ be a set and let $K$ be a ring. Recall that the \textbf{free
right  
$K$-module} $K[X]$ on $X$ has as elements all formal sums $x_1k_1+  
\cdots +x_nk_n$ where $x_1, \ldots, x_n \in X$ and $k_1, \ldots, k_n \in  
K$. Right multiplication by elements of $K$ and addition of elements of  
$K[X]$ are defined, with a zero and inverses, and  
$(x_1+x_2)k=x_1k+x_2k$.  
  
Let $P:=\{p_1, \ldots ,p_n\} \subseteq K[X]$.  
Recall that the \textbf{sub $\bZ G$-module generated by $P$} is  
$$\lan P \ran:= \{ p_1\zeta_1+ \cdots +p_n\zeta_n :  \zeta_1, \ldots,  
\zeta_n \in K\}$$

Let $grp\lan X|R \ran$ be a presentation of a group $G$.  
The group ring $\bZ G$ is the free right $\bZ$-module on $G$ together  
with a composition, making it an algebra over the ring $\bZ$. The free  
right $\bZ G$-module $\bZ G[R]$ on the set $R$ has elements of the form  
$r_1 \zeta_1 + \cdots + r_n \zeta_n$ where $r_1, \ldots, r_n \in R$ and  
$\zeta_1, \ldots, \zeta_n \in \bZ G$.  
  
\begin{lem}  
Let $grp\lan X|R \ran$ be a presentation of a group $G$, with quotient  
morphism $\theta:F(X) \to G$.  
Let $\iota=(r_1,u_1)^{\ep_1} \cdots (r_n,u_n)^{\ep_n}$ be an identity Y-  
sequence and let $\lambda$ denote the empty sequence.  
Define $\alpha:(Y^+ \sqcup Y^-)^* \to \bZ G[R]$ by  
$\alpha((r,u)^\ep):=r(\theta u\ep)$ with $\alpha(\lambda)=0$.  
Then $\iota \stackrel{*}{\lra}_{R_P} \lambda$ if and only if  
$\alpha(\iota)=0$.  
\end{lem}
  
\begin{proof}  
We verify that $\alpha$ preserves the $G$-action: $\alpha(((r,u)^\ep)^v)  
= \alpha((r,uv)^\ep)= r(\theta(uv)\ep) = (\alpha(r,u)^\ep)^{\theta v}$.  
The result now follows immediately from the definition of $\alpha$, the  
Primary Identity Property and the previous lemma.  
\end{proof}

\begin{cor}  
Let $\iota_1, \iota_2$ be identity Y-sequences. Then $\iota_1  
\stackrel{*}{\lra}_{R_P} \iota_2$ if and only if $\lan \iota_1 \ran =  
\lan \iota_2 \ran$ in $\bZ G [R]$.  
\end{cor}  
  
\begin{Def}  
Let $K[X]$ be a right $K$-module and let $a,b \in K[X]$. The Submodule  
Problem is  
\begin{center}  
\begin{tabular}{lll}  
INPUT & $a,b \in K[X]$ & (two elements of the right $K$-module,)\\  
QUESTION & $\lan a \ran = \lan b \ran$? & (do they generate the same  
submodule?)  
\end{tabular}  
\end{center}  
\end{Def}  
  
So we have shown that the Peiffer Problem for identity Y-sequences  
simplifies to the Submodule Problem. If the Submodule Problem can be  
solved then it is possible to reduce the set of generators of $ker  
\dd_2$ to a set  
of generating identities $\{\iota_1, \ldots \iota_t \}$ such that no  
subset of this will generate the same sub $\bZ G$-module. This is in  
some sense a minimal set of generators for $ker \dd_2$ (see later  
note).\\  
  
At the next levels, $ker \dd_n$ for $n \geq 3$, the problem is simpler  
in that we are now working entirely  
in $\bZ G$-modules, and do not encounter the Peiffer Problem.  
The only problem we now encounter is the Submodule Problem.\\  
  
In the $ker\dd_3$ case (Table 3) we have a set of 24 generators as  
elements  
of $C_2$, which here is the free $\bZ G$-module on  $\{\iota_1, \ldots
,\iota_4\}$.  
Some of these generators are zero, others are of the form  
$\iota_1(\id \! + \! x \! + \! x^2)$  
and $\iota_2(x \! + \! yx) - \iota_3(x \! + \! y)$.\\  
  
The problem may be phrased in the terms of a Gr\"obner basis problem.  
This is a reasonable approach, because methods for dealing with  
commutative  
Gr\"obner bases over rings exist \cite{AdLo} (essentially for Principal  
Ideal Domains) and methods for noncommutative Gr\"obner bases over rings  
(specifically group and monoid rings) are being developed \cite{Birgit}.  
Let $P:=\{p_1, \ldots ,p_n\}$ be a set of polynomials with coefficients  
in $\bZ G$ and  
monomials from a set $M$ ~ i.e. $p_1, \ldots ,p_n$ are elements of the  
$\bZ G$-module  
$\bZ G(M)$.  
The task is to find a set $Q:=\{q_1, \ldots ,q_m\}$ that generates the  
same sub  
$\bZ G$-module, but is such that no $q_i$ is a sum of $\bZ G$-multiples  
of the other $q_j$.\\  
  
Bases for modules are not in general unique or of the same rank.
So it is possible that there are two such  
sets $Q$ and $Q'$ and that these are of different sizes. We are  
concerned not  
with finding the generating set with smallest cardinality but with
finding a set which contains no subset which would generate the same
submodule.\\  
  
If $Q$ is a Gr\"obner basis for $P$ then by definition 
$\lan P\ran =\lan Q\ran $. If $Q$ is a reduced Gr\"obner basis then it  
is such that no element $q_i$ of $Q$ is a sum of
 $\bZ G$-multiples of the other elements $q_j$ of $Q$. This puts the problem  
of finding a reduced set of sub-module generators in terms of a Gr\"obner  
basis problem.

\section{Concluding Remarks}  
  
The purpose of this chapter was to make algorithmic the methods given in  
\cite{BrSa97}. In fact we have computerised the initial part of the  
construction, using rewriting  
theory and the Knuth-Bendix completion procedure to algorithmically  
define  
the first contracting homotopies $h_0$ and $h_1$. The program  
$idrels.g$ will compute, from a group presentation, a complete  
generating set  
for the module of identities among relations.\\

Unfortunately we cannot yet  
produce an algorithm for the minimalisation of this set of generators.  
Two major barriers to a reduction procedure have been identified.  
Firstly, the Peiffer Problem, a particularly difficult  
word problem encountered in crossed modules and 2-categories as a result  
of the Peiffer rules or interchange law.  
This has been reduced, using a property defined in \cite{BrHu} to  
the Submodule Problem, which is also encountered at higher levels, and  
indicates that methods for noncommutative Gr\"obner bases over group  
rings are required.  
Methods for solving this problem are progressing, thanks to
collaboration with Birgit Reinert (Kaiserslautern). A program for
reducing the first generating set of identities exists. This work will
continue  
with the aim of extending the program so that it will compute  
minimal generating sets for the $\bZ G$-modules $C_n$ for any given $n$.\\
  
Investigation of whether the completion of a monoid presentation  
yields something useful for the construction of a resolution of the  
monoid  
would also be an interesting area of work. We do not know whether the  
covering  
groupoids methods of \cite{BrSa97} might generalise to a covering  
categories of  
monoids method for calculating something corresponding to identities  
among  
relations for monoids. This looks like the beginnings of a  
noncommutative  
syzygy theory, and would definitely be worth investigating.  

\begin{appendix}

\textbf{\Large File 1: knuth.g}\\

The first program is an implementation of the standard Knuth-Bendix procedure
which may be applied to string rewriting.
A rewrite system $R$ is input in the form of a list $\mathtt{R}$ 
of pairs of words. The important subroutines are:\\

$\bullet \ \mathtt{OnePass(word,R)}$: reduces $\mathtt{word}$ (if possible) 
by applying one rule from $\mathtt{R}$.
This procedure involves searching to see if the left side of a rule in 
$\mathtt{R}$ 
is a subword of $\mathtt{word}$ and then replacing that part of $\mathtt{word}$ 
with the right side of the rule.

$\bullet \ \mathtt{ReduceWord(word,R)}$ reduces $\mathtt{word}$ as far as 
possible with respect to $\mathtt{R}$ by the
repeated application of the previous function. (Note that the reduced form
can only be guaranteed to be unique if $\mathtt{R}$ is complete.)

$\bullet \ \mathtt{CriticalPairs(R)}$: overlaps between the left hand sides 
of the rules in $\mathtt{R}$ are found, and the resulting critical pairs are 
found and reduced with respect to $\mathtt{R}$.

$\bullet \ \mathtt{OnePassKB(R)}$: this function computes the critical pairs 
of a rewrite system $\mathtt{R}$ and then resolves these critical pairs by 
adding then to $\mathtt{R}$.

$\bullet \ \mathtt{SystemReduce(R)}$: is an efficiency measure rather 
than theoretically essential. It normalises an ordinary rewrite system by 
reducing the rules (both sides of each rule are reduced by the other rules 
and the rules implied by other rules within the system are hence removed).\\

The main function of the program is $\mathtt{KB}$.\\

$\bullet \ \mathtt{KB(R)}$: attempts to complete the rewrite system (with respect to the length-lex order).
If it achieves the completion it returns the complete (reduced) rewrite system 
as a list of ordered pairs.\\

When the rewriting system is for a monoid there are further functions which
will enumerate the elements of the monoid.\\

$\bullet \ \mathtt{NextWords(F,Words)}$: creates new words of length 
$n\!+\!1$ by composing single generators from (the free group) $\mathtt{F}$ 
with irreducible words of length $n$.
 
$\bullet \ \mathtt{Enumerate(F,R)}$: uses the previous function and 
$\mathtt{reduce(word,R)}$ to
build up blocks of words of the same length (on the irreducibles one unit
shorter) and then to reduce these words as far as possible. When a whole block
of new words is reducible, there are no more irreducible words to be found.\\

\textbf{\Large File 2: kan.g}\\

The main function of the program is called $\mathtt{Kan}$. 
The input, functions and output are fully described in Chapter Two.

$\bullet \ \mathtt{InitialRules(KAN)}$:
The first sub-routine constructs the initial rewrite system of mixed one-sided 
and two-sided rules. All the rules of the form $(x\iota Fa,Xa(x))$ for $a \in \bA$
are added to the relations of the category $B$. This establishes an initial 
rewriting system for the group.

$\bullet \ \mathtt{Kan(KAN)}$:
This completes the rewriting system with respect to length-lex (where possible)
by calling  $\mathtt{knuth.g}$. 
It then enumerates the elements of the sets which make up the Kan
extension. The action of $\bB$ on the resulting elements 
can easily be computed.\\

\textbf{\Large File 3: ncpoly.g}\\

This file provides definitions and some operations for polynomials with
rational coefficients and non-commutative monomials in a semigroup.\\

$\bullet \ \mathtt{PolyFromTerms( [[k_1,m_1],..,[k_n,m_n]] )}$: creates a 
(noncommutative) polynomial from a list of terms. 
A polynomial is stored as a record but printed nicely
as a polynomial $\mathtt{k1 ~ m1 + \cdots + kn ~ mn}$. 
There are a number of operations:\\

$\bullet \ \mathtt{IsNonCommPoly( poly )}$: tests whether a record is a 
polynomial.

$\bullet \ \mathtt{LengthPoly( poly )}$: returns the number of terms.

$\bullet \ \mathtt{LeadTerm( poly )}$: extracts the leading term 
(which consists of the monomial of greatest size with respect to the 
length-lex order and its coefficient).

$\bullet \ \mathtt{LeadCoeff( poly )}$: returns the coefficient of the 
leading term.

$\bullet \ \mathtt{LeadMonom( poly )}$: returns the monomial part of 
the leading term.

$\bullet \ \mathtt{MakeMonic( poly )}$: divides a a polynomial by its 
leading coefficient to return a monic polynomial.

$\bullet \ \mathtt{NeatenPoly( poly )}$: adds like terms (non-destructive).

$\bullet \ \mathtt{poly_1=poly_2}$: equality between polynomials is well defined.

$\bullet \ \mathtt{AreEquivPolys(poly_1,poly_2)}$: polynomials are equivalent 
if one is a multiple of the other. 

$\bullet \ \mathtt{AddPoly(poly_1,poly_2)}:$ returns the `neatened' sum of two 
`neat' polynomials.
 
$\bullet \ \mathtt{SubtractPoly(poly_1,poly_2)}:$ returns the `neatened' 
difference of two `neat' polynomials.\\

To summarise: a polynomial record $\mathtt{poly}$ has the following fields:
$\mathtt{poly.IsNonComPoly}$ is true;
$\mathtt{poly.terms}$ is a list of terms $\mathtt{[c, m]}$ 
                      where $\mathtt{c}$ is a rational and $\mathtt{m}$ is a word;
$\mathtt{poly.isNeat}$ is either true or false;
$\mathtt{poly.operations}$ will be $\mathtt{NonCommPolyOps}$;
$\mathtt{poly.lead}$ is a term $\mathtt{[c, m]}$;
$\mathtt{poly.leadmon}$ is $\mathtt{poly.lead[2]}$;
$\mathtt{poly.isMonic}$ is either true or false.\\

All these functions are required for the noncommutative Gr\"obner basis
program.\\

\textbf{\Large File 4: grobner.g}\\

This is a program for computing the noncommutative Gr\"obner basis of a set of 
polynomials. It consists of a number of functions:\\

$\bullet \ \mathtt{ReducePoly( poly, POL )}$: reduces a polynomial $poly$ by subtracting 
multiples of polynomials in $\mathtt{POL}$. The reduced form can only be guaranteed 
to be unique with a Gr\"obner basis.

$\bullet \ \mathtt{OrderSystem( POL )}$: orders a set of polynomials with respect to their leading
monomials.

$\bullet \ \mathtt{PolySystemReduce( POL )}$: Removes polynomials which are sums of 
multiples of other polynomials in the system.

$\bullet \ \mathtt{SPolys(ALL,NEW)}$: compares two lists of polynomials for matches (if the lists
are equal then this is the standard procedure and finds all matches in the
system) and calculates the resulting S-polynomials.

$\bullet \ \mathtt{GB(POL)}$: returns (where possible) a Gr\"obner basis 
for a system of noncommutative polynomials over the rationals 
(with respect to the length-lex order).\\

\textbf{\Large File 5: idrel.g}\\

This program accepts as input a free group and a list of relators.
It goes through a number of calculations, including an ``extra information''
Knuth-Bendix completion procedure and returns a complete set of generators 
for the module of identities among relations. The input, functions and output
are fully described in Chapter Five, with examples.
\end{appendix}

\end{document}